\documentclass[oneside,english,american]{amsart}
\usepackage[T1]{fontenc}
\usepackage[latin9]{inputenc}
\usepackage{textcomp}
\usepackage{mathrsfs}
\usepackage{amstext}
\usepackage{amsthm}
\usepackage{amssymb}
\usepackage{graphicx}

\makeatletter
\numberwithin{equation}{section}
\numberwithin{figure}{section}
\theoremstyle{plain}
\newtheorem{thm}{\protect\theoremname}[subsection]
\theoremstyle{plain}
\newtheorem{lem}[thm]{\protect\lemmaname}
\theoremstyle{plain}
\newtheorem{prop}[thm]{\protect\propositionname}
\theoremstyle{plain}
\newtheorem{cor}[thm]{\protect\corollaryname}
\theoremstyle{remark}
\newtheorem*{rem*}{\protect\remarkname}
\theoremstyle{plain}
\newtheorem*{prop*}{\protect\propositionname}
\theoremstyle{definition}
\newtheorem{defn}[thm]{\protect\definitionname}
\theoremstyle{plain}
\newtheorem*{thm*}{\protect\theoremname}
\theoremstyle{remark}
\newtheorem{rem}[thm]{\protect\remarkname}

\usepackage[all]{xy}

\makeatother

\usepackage{babel}
\addto\captionsamerican{\renewcommand{\corollaryname}{Corollary}}
\addto\captionsamerican{\renewcommand{\definitionname}{Definition}}
\addto\captionsamerican{\renewcommand{\lemmaname}{Lemma}}
\addto\captionsamerican{\renewcommand{\propositionname}{Proposition}}
\addto\captionsamerican{\renewcommand{\remarkname}{Remark}}
\addto\captionsamerican{\renewcommand{\theoremname}{Theorem}}
\addto\captionsenglish{\renewcommand{\corollaryname}{Corollary}}
\addto\captionsenglish{\renewcommand{\definitionname}{Definition}}
\addto\captionsenglish{\renewcommand{\lemmaname}{Lemma}}
\addto\captionsenglish{\renewcommand{\propositionname}{Proposition}}
\addto\captionsenglish{\renewcommand{\remarkname}{Remark}}
\addto\captionsenglish{\renewcommand{\theoremname}{Theorem}}
\providecommand{\corollaryname}{Corollary}
\providecommand{\definitionname}{Definition}
\providecommand{\lemmaname}{Lemma}
\providecommand{\propositionname}{Proposition}
\providecommand{\remarkname}{Remark}
\providecommand{\theoremname}{Theorem}

\usepackage[table]{xcolor}

\usepackage{MnSymbol}

\begin{document}
\title{Foliations on Shimura varieties in positive characteristic}
\author{Eyal Z. Goren and Ehud de Shalit}
\keywords{\selectlanguage{english}%
Shimura variety, foliation}
\subjclass[2000]{\selectlanguage{english}%
11G18, 14G35}
\address{\selectlanguage{english}%
Eyal Z. Goren, McGill University, Montréal, Canada}
\address{\selectlanguage{english}%
eyal.goren@mcgill.ca}
\address{\selectlanguage{english}%
Ehud de Shalit, Hebrew University of Jerusalem, Israel}
\address{\selectlanguage{english}%
ehud.deshalit@mail.huji.ac.il}

\selectlanguage{american}%
\begin{abstract}
This paper is a continuation of \cite{G-dS1}. We study foliations
of two types on Shimura varieties $S$ in characteristic $p$. The
first, which we call ``tautological foliations'', are defined on
Hilbert modular varieties, and lift to characteristic $0$. The second,
the ``$V$-foliations'', are defined on unitary Shimura varieties
in characteristic $p$ only, and generalize the foliations studied
by us before, when the CM field in question was quadratic imaginary.
We determine when these foliations are $p$-closed, and the locus
where they are smooth. Where not smooth, we construct a ``successive
blow up'' of our Shimura variety to which they extend as smooth foliations.
We discuss some integral varieties of the foliations. We relate the
quotient of $S$ by the foliation to a purely inseparable map from
a certain component of another Shimura variety of the same type, with
parahoric level structure at $p$, to $S.$

\tableofcontents{}
\end{abstract}

\maketitle

\vspace{-1cm}
\section{Introduction}

Let $S$ be a non-singular variety over a field $k$, and $\mathcal{T}$
its tangent bundle. A (smooth) foliation on $S$ is a vector sub-bundle
of $\mathcal{T}$ that is closed under the Lie bracket. If $\text{\rm char}(k)=p>0$,
a $p$-foliation is a foliation that is, in addition, closed under
the operation $\xi\mapsto\xi^{p}.$ As explained below, such $p$-foliations
play an important role in studying purely inseparable morphisms $S\to S'.$ 

Foliations have been studied in many contexts, and our
purpose here is to explore certain $p$-foliations on Shimura varieties
of PEL type, which bear a relation to the underlying moduli problem.
The connection between the two topics can go either way. It can be
seen as using the rich geometry of Shimura varieties to produce interesting
examples of foliations, or, in the other direction, as harnessing
a new tool to shed light on some geometrical aspects of these Shimura
varieties, especially in characteristic $p$.

We study two types of foliations, that could well turn out to be particular
cases of a more general theory. The first lie on Hilbert modular varieties.
Let $S$ be a Hilbert modular variety associated with a totally real
field $L,$ $[L:\mathbb{Q}]=g$ (see the text for details). Let $\Sigma$
be any proper non-empty subset of $\mathscr{I}=\mathrm{Hom}(L,\mathbb{R}).$ Fixing
an embedding of $\bar{\mathbb{Q}}$ in $\bar{\mathbb{Q}}_{p}$ we
view $\mathscr{I}$ also as the set of embeddings of $L$ in $\mathbb{\bar{Q}}_{p}.$
If $p$ is unramified in $L$ these embeddings end up in $\mathbb{Q}_{p}^\text{\rm nr}$
and the Frobenius automorphism $\phi\in {Gal}(\mathbb{Q}_{p}^\text{\rm nr}/\mathbb{Q}_{p})$
permutes them. 

Consider the uniformization $\Gamma\setminus\mathfrak{H}^{\mathscr{I}}\simeq S(\mathbb{C}),$
describing the complex points of $S$ as a quotient of the product
of $g$ copies of the upper half plane (indexed by $\mathscr{I}$)
by an arithmetic subgroup of $\text{\rm SL}_{2}(L).$ The foliation $\mathscr{F}_{\Sigma}$,
labelled by the subset $\Sigma$, is defined most easily complex analytically:
at any point $x\in S(\mathbb{C})$ its fiber is spanned by $\partial/\partial z_{i}$
for $i\in\Sigma\subset\mathscr{I}.$ Our main results about it are
the following:
\begin{itemize}
\item Using the Kodaira-Spencer isomorphism, $\mathscr{F}_{\Sigma}$ can
be defined algebraically, hence also in the characteristic $p$ fiber
of $S$, for any good unramified rational prime $p.$ Its ``reduction
modulo $p$'' is a smooth $p$-foliation if and only if the subset
$\Sigma$ is invariant under Frobenius.
\item If the singleton $\{\sigma\}$ is not Frobenius-invariant (i.e. the
corresponding prime of $L$ is not of absolute degree 1), the obstruction
to $\mathscr{F}_{\{\sigma\}}$ being $p$-closed can be identified
with the square of a partial Hasse invariant.
\item Let $\mathfrak{p}$ be a prime of $L$ dividing $p$. Let $S_0(\mathfrak{p})$ be the special fiber of the integral model of the Hilbert modular variety with $\Gamma_0(\mathfrak{p})$-level structure studied in \cite{Pa}. 
When $\Sigma$ consists of all the
embeddings \textit{not} inducing $\mathfrak{p}$, the (purely inseparable) quotient of $S$ by the foliation $\mathscr{F}_{\Sigma}$
can be identified with a certain irreducible component of $S_0(\mathfrak{p})$.
\item It is easy to see that $\mathscr{F}_{\Sigma}$ does not have
any integral varieties in characteristic~$0$. In contrast, as we show below, any $p$-foliation in characteristic $p$ admits a plentiful supply of integral varieties.  
In our case, we show that certain Goren-Oort
strata in the reduction modulo $p$ of $S$ are integral varieties
of the foliation $\mathscr{F}_{\Sigma}$. 
\end{itemize}
The second class of foliations studied in this paper lives on unitary
Shimura varieties $M$ of arbitrary signature, associated with a CM
field $K$. They generalize the foliations studied in \cite{G-dS1}
when $K$ was quadratic imaginary, and are again labelled by subsets
$\Sigma$ of $\mathscr{I}^{+}=\mathrm{Hom}(L,\mathbb{R})$ where $L$
is the totally real subfield of~$K$. Unlike the foliations of the
first type, they are particular to the characteristic $p$ fibers,
where~$p$ is again a good unramified prime, and are of different
genesis. They are defined using the Verschiebung isogeny of the universal
abelian scheme over the $\mu$-ordinary locus $M^\text{\rm ord}$ of $M$.
Their study relies to a great extent on the work of Wedhorn and Moonen
cited in the bibliography. We refer to the text for the precise definition
of the foliation denoted by $\mathscr{F}_{\Sigma}$, as it is a bit
technical. A rough description is this: In the $p$-kernel of the
universal abelian scheme over $M^\text{\rm ord}$ lives an important
subgroup scheme, which played a special role in the work of Oort and
his school, namely the maximal subgroup scheme of $\alpha_{p}$-type.
Its cotangent space serves to define the $\mathscr{F}_{\Sigma}$ via
the Kodaira-Spencer isomorphism. The main results concerning $\mathscr{F}_{\Sigma}$
are the following:
\begin{itemize}
\item $\mathscr{F}_{\Sigma}$ is a smooth $p$-foliation on $M^\text{\rm ord},$
regardless of what $\Sigma$ is. Involutivity follows from the flatness
of the Gauss-Manin connection, but being $p$-closed is more delicate,
and is a consequence of a theorem of Cartier on the $p$-curvature
of that connection.
\item Although in general more complicated than the relation found in \cite{G-dS1},
one can work out explicitly the relation between the foliation $\mathscr{F}_{\Sigma}$
and the ``cascade structure'' defined by Moonen \cite{Mo} on the
formal completion $\widehat{M}_{x}$ of $M$ at a $\mu$-ordinary
point $x$. While the cascade structure does not globalize, its ``trace''
on the tangent space, constructed from the foliations $\mathscr{F}_{\Sigma},$
does globalize neatly. 
\item There is a maximal open subset $M_{\Sigma}\subset M$ to which $\mathscr{F}_{\Sigma}$
extends as a smooth $p$-foliation with the same definition used to
define it on $M^\text{\rm ord}.$ This $M_{\Sigma}$ is a union of Ekedahl-Oort
(EO) strata, and in fact consists of all the strata containing in
their closure a smallest one, denoted $M_{\Sigma}^\text{\rm fol}.$ The description
of which EO strata participate in $M_{\Sigma}$ is given combinatorially
in terms of ``shuffles'' in the Weyl group.
\item Outside $M_{\Sigma}$ the foliation $\mathscr{F}_{\Sigma}$ acquires
singularities, but we construct a ``successive blow up'' $\beta:M^{\Sigma}\to M$,
which is an isomorphism over $M_{\Sigma}$, to which the lifting of
$\mathscr{F}_{\Sigma}$ extends as a smooth $p$-foliation. This $M^{\Sigma}$
is an interesting (characteristic $p$) moduli problem in its own
right. It is non-singular, and the extension of $\mathscr{F}_{\Sigma}$
to it is transversal to the fibers of $\beta$.
\item When $K$ is quadratic imaginary, the EO stratum $M_{\Sigma}^\text{\rm fol}$
was proved to be an integral variety (in the sense of foliations) of $\mathscr{F}_{\Sigma}$. A
similar result is expected here when $\Sigma=\mathscr{I}^{+}.$ This
can be probably proved via elaborate Dieudonné module computations,
as in \cite{G-dS1}, but in this paper we content ourselves with checking
that the dimensions match. 
\item A natural interesting question is to identify the purely inseparable
quotient of $M^{\Sigma}$ by (the extended) $\mathscr{F}_{\Sigma}$
with a certain irreducible component of the special fiber of the Rapoport-Zink
model of a unitary Shimura variety with parahoric level structure
at $p$. This was done in our earlier paper when~$K$
was quadratic imaginary, and was used to obtain some new results on
the geometry of that particular irreducible component. In the general
case treated in this paper, we know of a natural candidate with which
we would like to identify the quotient of $M^{\Sigma}$ by $\mathscr{F}_{\Sigma}.$
However, following the path set in \cite{G-dS1} for a general CM
field $K$ would require a significant amount of work, and we leave
this question for a future paper.
\end{itemize}

\

Section 2 is a brief review of general results on foliations, especially
in characteristic~$p$. The main two sections of the paper, Sections
3 and 4, are devoted to the two types of foliations, respectively.

\subsubsection{Notation}
\begin{itemize}
\item For any commutative $\mathbb{F}_{p}$-algebra $R$ we let $\phi:R\to R$
be the homomorphism $\phi(x)=x^{p}.$
\item If $S$ is a scheme over $\mathbb{F}_{p},$ $\Phi_{S}$ denotes its
absolute Frobenius morphism. It is given by the identity on the underlying
topological space of $S$, and by the map $\phi$ on its structure
sheaf. If $\mathcal{H}$ is an $\mathcal{O}_{S}$-module then we write
$\mathcal{H}^{(p)}$ (or $\mathcal{H}^{(p)/S}$) for $\Phi_{S}^{*}\mathcal{H}=\mathcal{O}_{S}\otimes_{\phi,\mathcal{O}_{S}}\mathcal{H}$.
\item If $T\to S$ is a morphism of schemes over $\mathbb{F}_{p},$ $T^{(p)}$
(or $T^{(p)/S}$) is $S\times_{\Phi_{S},S}T$ and $\text{\rm Fr}_{T/S}:T\to T^{(p)}$
is the relative Frobenius morphism, characterized by the relation
$pr_{2}\circ \text{\rm Fr}_{T/S}=\Phi_{T}.$
\item If $A\to S$ is an abelian scheme and $A^{t}\to S$ is its dual, then
$\text{\rm Fr}_{A/S}$ is an isogeny and Verschiebung $\text{\rm Ver}_{A/S}:A^{(p)}\to A$
is the dual isogeny of $\text{\rm Fr}_{A^{t}/S}.$
\item If $\mathcal{H}$ is an $\mathcal{O}_{S}$-module with $\mathcal{O}$
action (for some ring $\mathcal{O}$), and $\tau:\mathcal{O}\to\mathcal{O}_{S}$
is a homomorphism, then
\[
\mathcal{H}[\tau]=\{\alpha\in\mathcal{H}|\,\forall a\in\mathcal{O}\,\,a.\alpha=\tau(a)\alpha\}.
\]
If $T:\mathcal{H\to\mathcal{G}}$ is a homomorphism of sheaves of
modules, we denote $\ker T=\mathcal{H}[T].$
\item If $x\in S$, the fiber of $\mathcal{H}$ at $x$ is denoted $\mathcal{H}_{x}.$
This is a vector space over the residue field $k(x).$ The same notation
is used for the fiber $\mathcal{H}_{x}=x^{*}\mathcal{H}$ at a geometric
point $x:\text{\rm Spec}(k)\to S$.
\item If $\mathcal{H}^{\vee}$ is the dual of a locally free $\mathcal{O}_{S}$-module
$\mathcal{H}$ we denote the pairing $\mathcal{H}^{\vee}\times\mathcal{H}\to\mathcal{O}_{S}$
by $\left\langle ,\right\rangle .$
\item By the Dieudonné module of a $p$-divisible group over a perfect field
$k$ in characteristic $p$, or of a finite commutative $p$-torsion
group scheme over $k,$ we understand its \emph{contravariant} Dieudonné
module.
\end{itemize}

\subsubsection{Acknowledgements}

We would like to thank N. Shepherd-Barron and R. Taylor for sharing
their unpublished manuscript \cite{E-SB-T} with us. This research
was supported by ISF grant 276.17 and NSERC grant 223148 and the hospitality
of McGill University and the Hebrew University.

\section{Generalities on Foliations}

\subsection{Smooth foliations}

Let $k$ be a perfect field and $S$ a $d$-dimensional smooth $k$-variety, $d>0$.
Let $\mathcal{T}$ denote the tangent bundle of $S$. If $U\subset S$
is Zariski open, then sections $\xi\in\mathcal{T}(U)$ are \emph{vector
fields} on $U$ and act on $\mathcal{O}_{S}(U)$ as derivations. The
space $\mathcal{T}(U)$ has a structure of a Lie algebra (of infinite
dimension) over $k,$ when we define
\[
[\xi,\eta](f)=\xi(\eta(f))-\eta(\xi(f)).
\]
A \emph{foliation $\mathscr{\mathcal{F}}$} on $S$ is a saturated
subsheaf $\mathcal{F}\subset\mathcal{T}$ closed under the Lie bracket.
For every Zariski open set $U,$ the vector fields on $U$ along the
foliation form a saturated $\mathcal{O}_{S}(U)$-submodule $\mathcal{F}(U)\subset\mathcal{T}(U)$
closed under the Lie bracket, i.e. if $f\in\mathcal{O}_{S}(U),$ $\xi\in\mathcal{T}(U)$
and $f\xi\in\mathcal{F}(U)$ then $\xi\in\mathcal{F}(U),$ and if
$\xi,\eta\in\mathcal{F}(U)$ then $[\xi,\eta]\in\mathcal{F}(U).$

The foliation $\mathcal{F}$ is called \emph{smooth} if it is a vector
sub-bundle of $\mathcal{T},$ namely if both~$\mathcal{F}$ and $\mathcal{T}/\mathcal{F}$
are locally free sheaves. Since $\mathcal{F}$ is assumed to be saturated,
and since a torsion-free finite module over a discrete valuation ring
is free, the locus $Sing(\mathcal{F})$ where $\mathcal{F}$ is \emph{not}
smooth is a closed subset of $S$ of codimension $\ge 2$. As an example,
the vector field $x\partial/\partial x+y\partial/\partial y$ generates
a rank-1 foliation on $\mathbb{A}^{2},$ whose singular set is the
origin.

If $k=\mathbb{C}$ then by a well known theorem of Frobenius every
$x\in S-Sing(\mathcal{F})$ has a \emph{classical} open neighborhood
$V\subset S(\mathbb{C})$ which can be decomposed into a disjoint
union of parallel leaves $L$ of the foliation. Each leaf $L$ is a smooth
complex submanifold of $V$ and if $y\in L$ then the tangent space
$T_{y}L\subset T_{y}S=\mathcal{T}_{y}$ is just the fiber of $\mathcal{F}$
at $y$. One says that $L$ is an \emph{integral subvariety} of the
foliation. Whether a smooth foliation has an \emph{algebraic} integral
subvariety passing through a given point $x$, and the classification
of all such integral subvarieties, is in general a hard problem, see
\cite{Bo}.

There is a rich literature on foliations on algebraic varieties. See,
for example, the book \cite{Br}.

\subsection{Foliations in positive characteristic}

\subsubsection{$p$-foliations and purely inseparable morphisms of height 1}

If $\text{\rm char}(k)=p$ is positive then $\mathcal{T}$ is a $p$-Lie algebra,
namely if $\xi\in\mathcal{T}(U)$ then $\xi^{p}=\xi\circ\cdots\circ\xi$
(composition $p$ times) is also a derivation, hence lies in $\mathcal{T}(U).$
A foliation $\mathcal{F}$ is called a $p$\emph{-foliation} if it
is $p$-\emph{closed}: whenever $\xi\in\mathcal{F}(U)$, then $\xi^{p}\in\mathcal{F}(U)$
as well.

The interest in $p$-foliations in characteristic $p$ stems from
their relation to purely inseperable morphisms of height 1. The following
theorem has its origin in Jacobson's \emph{inseparable Galois theory
}for field extensions (\cite[\S8]{Jac}). We denote by $\phi:k\to k$
the $p$-power map, by $S^{(p)}=k\times_{\phi,k}S$ the $p$-transform
of $S$, and by
\[
\text{\rm Fr}_{S/k}:S\to S^{(p)}
\]
the relative Frobenius morphism. We denote by $\Phi_{S}:S\to S$ the
absolute Frobenius of $S.$ Thus,
\[
\Phi_{S}=pr_{2}\circ \text{\rm Fr}_{S/k},
\]
where $pr_{2}:S^{(p)}=k\times_{\phi,k}S\to S$ is the projection onto
the second factor.
\begin{thm}
\label{thm:Quotient by foliation}\cite{Ek} Let $k$ be a perfect
field, $\text{\rm char}(k)=p.$ Let $S$ be a smooth $k$-variety and denote
by $\mathcal{T}$ its tangent bundle. There exists a one-to-one correspondence
between smooth $p$-foliations $\mathcal{F}\subset\mathcal{T}$ and
factorizations of the relative Frobenius morphism $\text{\rm Fr}_{S/k}=g\circ f,$
\[
S\overset{f}{\to}T\overset{g}{\to}S^{(p)},
\]
where $T$ is a smooth $k$-variety (equivalently, where $f$ and $g$
are finite and flat). We call $T$ the quotient of $S$ by the foliation  $\mathcal{F}$.

Given $\mathcal{F}$, if (locally) $S=\text{\rm Spec}(A)$, then $T=\text{\rm Spec}(B)$
where $B=A^{\mathcal{F}=0}$ is the subring annihilated by $\mathcal{F}$,
and $f$ is induced by the inclusion $B\subset A.$ Conversely, given
a factorization as above, then $\mathcal{F}=\ker(df)$ where $df$
is the map induced by~$f$ on the tangent bundle.

Furthermore, if $r=\mathrm{rk}(\mathcal{F})$ then $\deg(f)=p^{r}.$
\end{thm}

As mentioned above, the birational version of this theorem is due
to Jacobson. From this version one deduces rather easily a correspondence
as in the theorem, when $T$ is only assumed to be \emph{normal, }and
$\mathcal{F}$ is saturated, closed under the Lie bracket and $p$-closed,
but not assumed to be smooth. The main difficulty is in showing that
$T$ is smooth if and only if $\mathcal{F}$ is smooth, i.e. locally
a direct summand everywhere. The reference \cite{Ek} only cites the
work of Yuan \cite{Yuan} and of Kimura and Nitsuma \cite{Ki-Ni},
but does not give the details. The proof in the book \cite{Mi-Pe}
seems to be wrong. A full account may be found in \cite{Li}.

\subsubsection{The obstruction to being $p$-closed\label{subsec:obstruction p closed}}
\label{subsubsec: obstruction}

If $\mathcal{F}\subset\mathcal{T}$ is a smooth foliation, the map
$\xi\mapsto\xi^{p}\mod\mathcal{F}$ induces an $\mathcal{O}_{S}$-linear
map of vector bundles
\[
\kappa_{\mathcal{F}}:\Phi_{S}^{*}\mathcal{F}\to\mathcal{T}/\mathcal{F}
\]
which is identically zero if and only if $\mathcal{F}$ is $p$-closed.
See \cite{Ek}, Lemma 4.2(ii). We call the map $\kappa_{\mathcal{F}}$
the \emph{obstruction to being} $p$-\emph{closed}.

\subsection{Integral varieties in positive characteristic}
\label{sec:general integral}
In contrast to the situation in characteristic 0, integral varieties
of $p$-foliations in characteristic $p$ abound, and are easily described.
The goal of this section is to clarify their construction.\footnote{This construction shows that the conjecture ``of André-Oort type'',
suggested in §5.3 of \cite{G-dS1}, is far from being true.} As in the previous section we let $S$ be a smooth $d$-dimensional quasi-projective
variety over a perfect field $k$ of characteristic $p,$ $\mathcal{T}=TS$
its tangent bundle, and $\mathcal{F}$ a smooth $p$-foliation of
rank $r$. We denote by $S\overset{f}{\to}T$ the quotient of $S$
by $\mathcal{F}$, as in Theorem~\ref{thm:Quotient by foliation}.

Let $\mathcal{G}=\mathcal{T}/\mathcal{F}=\mathrm{Im}(df).$ This is
a smooth $p$-foliation of rank $d-r$ on $T$ and the quotient of
$T$ by $\mathcal{G}$ is $T\overset{g}{\to}S^{(p)}.$ The factorizations
of the relevant Frobenii are $Fr_{S/k}=g\circ f$ and $Fr_{T/k}=f^{(p)}\circ g.$
\begin{defn}
Let $\iota:W\hookrightarrow S$ be a closed subvariety of $S$ and
$W^{\mathrm{sm}}$ the (open dense) smooth locus in $W$. We say that
$W$ is an \emph{integral variety} of $\mathcal{F}$ if at every $x\in W^{\mathrm{sm}}$
we have $T_{x}W=\iota^{*}\mathcal{F}_{x}$ (as subspaces of $\iota^{*}T_{x}S$).
In this case $\dim W=r.$ We say that $W$ is \emph{transversal} to
$\mathcal{F}$ at $x\in W^{\mathrm{sm}}$ if $T_{x}W\cap\iota^{*}\mathcal{F}_{x}=0,$
and that it is \emph{generically transversal to $\mathcal{F}$ }if
the set of points where it is transversal is a dense open set of $W.$
\end{defn}

\begin{rem}
Unlike the case of characteristic 0, an integral subvariety of a smooth
$p$-foliation need not be smooth. Consider, for example, the foliation
generated by $\partial/\partial v$ on $\mathbb{A}^{2}=Spec(k[u,v])$.
The irreducible curve
\[
u(u+v^{p})+v^{2p}=0
\]
is an integral curve of the foliation, but is singular at $x=(0,0)$.
The curve $u-v^{2}=0$ is generically transversal to the same foliation,
although it is not transversal to it at $x$.

If $S=Spec(A)$ and $W=Spec(A/I)$ for a prime ideal $I$, then regarding
$\mathcal{F}$ as a submodule of the module of derivations of $A$
over $k,$ $W$ is an integral variety of $\mathcal{F}$ if and only
if $\mathcal{F}(I)\subset I.$
\end{rem}

\begin{prop}
Let $\iota:W\hookrightarrow S$ be a closed $r$-dimensional subvariety
of $S$ and $Z=f(W)\hookrightarrow T$ the corresponding subvariety
of $T$ (also $r$-dimensional). Then the following are equivalent:
\end{prop}

\begin{enumerate}
\item $f_{W}:W\to Z$ is purely inseparable of degree $p^{r}$;
\item $W$ is an integral variety of $\mathcal{F}$;
\item $g_{Z}:Z\to W^{(p)}$ is a birational isomorphism;
\item $Z$ is generically transversal to $\mathcal{G}.$
\end{enumerate}
\begin{proof}
Remark first that since $g_{Z}\circ f_{W}=Fr_{W/k}$ induces a purely
inseparable field extension $k(W^{(p)})\subset k(W)$ of degree $p^{r},$
(1) and (3) are equivalent, and in fact are equivalent to $g_{Z}$
being separable (generically étale). 

Let $x\in W^{\mathrm{sm}}$ be such that $y=f(x)\in Z^{\mathrm{sm}}.$
The commutative diagram
\[\xymatrix@M=0.3cm{T_{x}W \ar[d]_{df_{W,x}} \ar@{^{(}->}[r] & T_{x}S\ar[d]^{df_x}\\
T_{y}Z \ar@{^{(}->}[r] & T_{y}T
}
\]
and the fact that $\ker(df)=\mathcal{F}$ tell us that
\[
\ker(df_{W})=TW^{\mathrm{sm}}\cap\iota^{*}\mathcal{F}.
\]
It follows that $f_{W}$ is purely inseparable of degree $p^{r}$
($df_{W}=0$) if and only if $TW^{\mathrm{sm}}$ and $\iota^{*}\mathcal{F}$,
both rank-$r$ vector sub-bundles of $\iota^{*}TS,$ coincide along
$W^{\mathrm{sm}}.$ This shows the equivalence of $(1)$ and $(2)$.
It also follows that $f_{W}$ is separable (generically étale) if
and only for generic $x$ we have $T_{x}W^{\mathrm{sm}}\cap\iota^{*}\mathcal{F}_{x}=0.$
When applied to $g$ and~$Z$, instead of $f$ and $W$, this gives
the equivalence of $(3)$ and $(4)$.
\end{proof}
\begin{thm}
Notation as above, any two points $x_{1},x_{2}$ of $S$ lie on an
integral variety of $\mathcal{F}.$ 
\end{thm}

\begin{proof} One may assume that $0< {\rm rank}(\mathcal{F}) < d = {\rm dim}(S)$. We have the diagram $\xymatrix@C=0.8cm{S \ar[r]^f & T \ar[r]^g& S^{(p)}}$, where the first arrow is dividing by the foliation $\mathcal{F}$ and the second arrow by the foliation $\mathcal{G}$. 
Let $y_{i}=f(x_{i})\in T.$ A standard application of Bertini's theorem \cite[Theorem II 8.18]{H} shows that there is a subvariety
$Z\subset T$ of dimension $r$ passing through the points~$y_{i}$ which
is generically transversal to $\mathcal{G}.$ Choosing $W$ so that
$g(Z)=W^{(p)},$ and hence $f(W)=Z,$ we conclude from the previous Proposition
that $W$ is an integral variety of $\mathcal{F}$ passing through
$x_{1}$ and~$x_{2}.$
\end{proof}
Thus integral varieties in characteristic $p$ abound, and are easy
to classify. Nevertheless, given a particular subvariety $W$, it
is still interesting to decide whether it is an integral variety of
$\mathcal{F}$ or not.

\section{Tautological foliations on Hilbert modular varieties}

\subsection{Hilbert modular schemes}

Let $L$ be a totally real field, $[L:\mathbb{Q}]=g\ge2$, $N\ge4$
an integer, and $\mathfrak{c}$ a fractional ideal of $L$, relatively
prime to $N$, called the \emph{polarization module}. We denote by
$\mathfrak{d}$ the different ideal of $L/\mathbb{Q}$. Let $D=\mathrm{disc}_{L/\mathbb{Q}}$
be the absolute discriminant of $L$.

Consider the moduli problem over $\mathbb{Z}[(ND)^{-1}],$ attaching
to a scheme $S$ over $\mathbb{Z}[(ND)^{-1}]$ the set $\mathscr{M}(S)$
of isomorphism classes of four-tuples
\[
\underline{A}=(A,\iota,\lambda,\eta),
\]
where
\begin{itemize}
\item $A$ is an abelian scheme over $S$ of relative dimension $g$.
\item $\iota:\mathcal{O}_{L}\hookrightarrow\mathrm{End}(A/S)$ is an injective
ring homomorphism making the tangent bundle $TA$ a locally free sheaf
of rank $1$ over $\mathcal{O}_{L}\otimes\mathcal{O}_{S}$. We denote
by $A^{t}$ the dual abelian scheme, and by $\iota(a)^{t}$ the dual
endomorphism induced by $\iota(a).$
\item $\lambda:\mathfrak{c}\otimes_{\mathcal{O}_{L}}A\simeq A^{t}$ is a
$\mathfrak{c}$-polarization of $A$ in the sense of \cite{Ka} (1.0.7).
This means that $\lambda$ is an isomorphism of abelian schemes compatible
with the $\mathcal{O}_{L}$ action, where $a\in\mathcal{O}_{L}$ acts
on the left hand side via $a\otimes1=1\otimes\iota(a)$ and on the
right hand side via $\iota(a)^{t}.$ Furthermore, under the identification
\[
\mathrm{Hom}_{\mathcal{O}_{L}}(A,A^{t})\simeq\mathrm{Hom}_{\mathcal{O}_{L}}(A,\mathfrak{c}\otimes_{\mathcal{O}_{L}}A)
\]
induced by $\lambda,$ the symmetric elements on the left hand side
(the elements $\alpha$ satisfying $\alpha^{t}=\alpha$ after we canonically
identify $A^{tt}$ with $A$) correspond precisely to the elements
of $\mathfrak{c}$, and those arising from an ample line bundle~$\mathcal{L}$
on $A$ via $\alpha_{\mathcal{L}}(x)=[\tau_{x}^{*}\mathcal{L}\otimes\mathcal{L}^{-1}]$
correspond to the totally positive cone in~$\mathfrak{c}.$ Note that
the Rosati involution induced by $\lambda$ on $\mathcal{O}_{L}$
is the identity.
\item $\eta$ is a $\Gamma_{00}(N)$-level structure on $A$ in the sense
of \cite{Ka} (1.0.8), i.e. a closed immersion of $S$-group schemes
$\eta:\mathfrak{d}^{-1}\otimes\mu_{N}\hookrightarrow A[N]$ compatible
with $\iota$.
\end{itemize}
The moduli problem $\mathscr{M}$ is representable by a smooth scheme
over $\mathbb{Z}[(ND)^{-1}],$ of relative dimension $g$, which we
denote by the same letter $\mathscr{M}$, and call the \emph{Hilbert
modular scheme} (see \cite{Ra}). If we want to remember the dependence on $\mathfrak{c}$
we use the notation $\mathscr{M}^{\mathfrak{c}}$ for $\mathscr{M}.$
The Hilbert moduli scheme admits smooth toroidal compactifications,
depending on some extra data. See \cite{Lan, Ra}.

The complex points $\mathscr{M}(\mathbb{C})$ of $\mathscr{M}$ may
be identified, as a complex manifold, with $\Gamma\backslash\mathfrak{H}^{g}$,
where $\mathfrak{H}$ is the upper half plane, and $\Gamma\subset \text{\rm SL}(\mathcal{O}_{L}\oplus\mathfrak{d}^{-1}\mathfrak{c}^{-1})$
is some congruence subgroup.

If $\mathfrak{c}=\gamma\mathfrak{c'}$, where $\gamma\gg0$, i.e. $\gamma$ is a totally positive element of $L$,  then $\lambda\circ(\gamma\otimes1):\mathfrak{c}'\otimes_{\mathcal{O}_{L}}A\simeq A^{t}$
is a $\mathfrak{c}'$-polarization of~$A$. Thus only the strict ideal
class of $\mathfrak{c}$ in $Cl^{+}(L)$ matters: the moduli schemes
$\mathscr{M}^{\mathfrak{c}}$ and $\mathscr{M}^{\mathfrak{c}'}$ are
isomorphic, via an isomorphism depending on the choice of $\gamma.$

Fix a prime number $p$ which is unramified in $L$ and relatively
prime to $N$. By the last remark, we may (and do) assume that $p$
is also relatively prime to $\mathfrak{c}.$ Write
\[
p\mathcal{O}_{L}=\mathfrak{p}_{1}\dots\mathfrak{p}_{r},\,\,\,f_{i}=f(\mathfrak{p}_{i}/p)
\]
the inertia degree of $\mathfrak{p}_{i}$, and let $\kappa$ be a
finite field into which all the $\kappa(\mathfrak{p}_{i})=\mathcal{O}_{L}/\mathfrak{p}_{i}$
embed, i.e. $\kappa=\mathbb{F}_{p^{n}}$ where $n$ is divisible by
$\mathrm{lcm}\{f_{1},\dots,f_{r}\}.$ Let $W(\kappa)$ be the Witt
vectors of $\kappa$. We have a decomposition
\[
\mathbb{B}=\mathrm{Hom}(L,W(\kappa)[p^{-1}])=\coprod_{\mathfrak{p}|p}\mathbb{B}_{\mathfrak{p}}
\]
indexed by the $r$ primes of $\mathcal{O}_{L}$ dividing $p,$ where
$\sigma:L\hookrightarrow W(\kappa)[1/p]$ belongs to $\mathbb{B}_{\mathfrak{p}}$
if $\sigma^{-1}(pW(\kappa))\cap\mathcal{O}_{L}=\mathfrak{p}$. The
Frobenius automorphism of $W(\kappa)$, denoted $\phi$, operates
on $\mathbb{B}$ from the left via $\sigma\mapsto\phi\circ\sigma$,
and permutes each $\mathbb{B}_{\mathfrak{p}}$ cyclically.

We denote by
\[
M=\mathscr{M}\otimes_{\mathbb{Z}[(ND)^{-1}]}\mathbb{F}_{p}
\]
the special fiber of $\mathscr{M}$ at $p$. Note that $M$ is smooth
over $\mathbb{F}_{p}$.

\subsection{The tautological foliations}

\subsubsection{The definition}

Complex analytically, the complex manifold $\Gamma\backslash\mathfrak{H}^{g}$
admits $g$ tautological rank 1 smooth foliations, generated by the
vector fields $\partial/\partial z_{i}$, where $z_{i}$ ($1\le i\le g$)
are the coordinate functions. As the derivations $\partial/\partial z_{i}$
commute with each other, any $r$ of them generate a rank $r$ analytic
foliation on $\mathscr{M}(\mathbb{C})$. We shall now show that these
foliations are in fact algebraic, and address the question which of
them descends modulo $p$ to a smooth $p$-foliation on $M$. We shall
then relate the quotients of $M$ by these foliations to Hilbert modular
varieties with Iwahori level structure at $p$.

\medskip{}

\textbf{Convention}. From now on we denote by $\mathscr{M}$ the base
change of the Hilbert modular scheme from $\mathbb{Z}[(ND)^{-1}]$
to $W(\kappa)$ and by $M$ its special fiber, a smooth variety over
$\kappa$. Recall that $\kappa$ is assumed to be large enough so
that there are $g$ distinct embeddings $\sigma:L\hookrightarrow W(\kappa)[1/p].$

\medskip{}

Let $\pi:A^\text{\rm univ}\to\mathscr{M}$ denote the universal abelian variety
over $\mathscr{M}$ and 
\[
\underline{\omega}=\pi_{*}\Omega_{A^\text{\rm univ}/\mathscr{M}}^{1}
\]
its \emph{Hodge bundle}. Since $p\nmid\mathrm{disc}_{L/\mathbb{Q}}$,
it decomposes under the action of $\mathcal{O}_{L}$ as a direct sum
of $g$ line bundles
\[
\underline{\omega}=\oplus_{\sigma\in\mathbb{B}}\mathscr{L}_{\sigma}
\]
where $\mathscr{L}_{\sigma}=\{\alpha\in\underline{\omega}|\,\iota(a)^{*}(\alpha)=\sigma(a)\alpha\,\,\forall a\in\mathcal{O}_{L}\}.$

Let $\underline{\text{\rm Lie}}=\underline{\text{\rm Lie}}(A^\text{\rm univ}/\mathscr{M})=\underline{\omega}^{\vee}$
be the relative tangent space of the universal abelian variety. The
Kodaira-Spencer isomorphism is an isomorphism of $\mathcal{O}_{\mathscr{M}}$-modules
(\cite{Ka}, (1.0.19)-(1.0.20))
\begin{multline}
\text{\rm KS}:\mathcal{T}_{\mathscr{M}/W(\kappa)}\simeq\mathrm{Hom}_{\mathcal{O}_{L}\otimes\mathcal{O}_{\mathscr{M}}}(\underline{\omega},\underline{\text{\rm Lie}}((A^\text{\rm univ})^{t}/\mathscr{M}))\label{eq:KS}
\\
\simeq\mathrm{Hom}_{\mathcal{O}_{L}\otimes\mathcal{O}_{\mathscr{M}}}(\underline{\omega},\underline{\text{\rm Lie}}\otimes_{\mathcal{O}_{L}}\mathfrak{c})=\underline{\text{\rm Lie}}^{\otimes2}\otimes_{\mathcal{O}_{L}}\mathfrak{dc},
\end{multline}
where the second isomorphism results from the polarization $\lambda^\text{\rm univ}:\mathfrak{c}\otimes_{\mathcal{O}_{L}}A^\text{\rm univ}\simeq(A^\text{\rm univ})^{t}$, and the $\otimes^2$ is the tensor product of $\mathcal{O}_L \otimes \mathcal{O}_\mathscr{M}$-modules.
Since $(p,\mathfrak{dc})=1$ we have 
\[
\underline{\text{\rm Lie}}\otimes_{\mathcal{O}_{L}}\mathfrak{dc}=\underline{\text{\rm Lie}}\simeq\oplus_{\sigma\in\mathbb{B}}\mathscr{L}_{\sigma}^{\vee}.
\]
We therefore get from $\text{\rm KS}$ a canonical decomposition of the tangent
space of $\mathscr{M}$ into a direct sum of $g$ line bundles
\[
\mathcal{T}_{\mathscr{M}/W(\kappa)}\simeq\oplus_{\sigma\in\mathbb{B}}\mathscr{L}_{\sigma}^{-2}.
\]

We denote by $\mathscr{F}_{\sigma}$ the line sub-bundle of $\mathcal{T}_{\mathscr{M}/W(\kappa)}$
corresponding to $\mathscr{L}_{\sigma}^{-2}$ under this isomorphism.
\begin{lem}
Let $\Sigma\subset\mathbb{B}$ be any set of embeddings of $L$ into
$W(\kappa)[1/p].$ Then $\mathscr{F}_{\Sigma}=\oplus_{\sigma\in\Sigma}\mathscr{F}_{\sigma}$
is involutive: if $\xi,\eta$ are sections of $\mathscr{F}_{\Sigma}$,
so is $[\xi,\eta].$
\end{lem}

\begin{proof}
Recall that the Kodaira-Spencer isomorphism is derived from the Gauss-Manin
connection
\[
\nabla:H_{dR}^{1}(A^\text{\rm univ}/\mathscr{M})\to H_{dR}^{1}(A^\text{\rm univ}/\mathscr{M})\otimes_{\mathcal{O}_{\mathscr{M}}}\Omega_{\mathscr{M}/W(\kappa)}^{1}.
\]
If $\xi$ is a section of $\mathcal{T}_{\mathscr{M}/W(\kappa)}$ we
denote by $\nabla_{\xi}:H_{dR}^{1}(A^\text{\rm univ}/\mathscr{M})\to H_{dR}^{1}(A^\text{\rm univ}/\mathscr{M})$
the map obtained by contraction with $\xi$. The Gauss-Manin connection
is well-known to be flat, namely
\[
\nabla_{[\xi,\eta]}=\nabla_{\xi}\circ\nabla_{\eta}-\nabla_{\eta}\circ\nabla_{\xi}.
\]
Now, $\nabla_{\xi}$ commutes with $\iota(a)^{*}$ for $a\in\mathcal{O}_{L}$,
and therefore preserves the $\sigma$-isotypic component $H_{dR}^{1}(A^\text{\rm univ}/\mathscr{M})[\sigma]$
for each $\sigma\in\mathbb{B}$. \emph{By} \emph{definition,} $\xi\in\mathscr{F}_{\sigma}$
if for every $\tau\ne\sigma$ the operator $\nabla_{\xi}$ maps the subspace
\[
\mathscr{L}_{\tau}=\underline{\omega}[\tau]\subset H_{dR}^{1}(A^\text{\rm univ}/\mathscr{M})[\tau]
\]
to itself. Similarly, $\xi\in\mathscr{F}_{\Sigma}$ if the same holds
for every $\tau\notin\Sigma.$ It follows at once from the flatness
of $\nabla$ that if this condition holds for $\xi$ and $\eta$,
it holds for $[\xi,\eta].$
\end{proof}
We conclude that $\mathscr{F}_{\Sigma}$ is a smooth foliation. We
call these foliations \emph{tautological.}

\subsubsection{The main theorem}

We consider now the foliations $\mathscr{F}_{\Sigma}$ in the special
fiber $M=\mathscr{M}\times_{W(\kappa)}\kappa$ only. The following
theorem summarizes the main results in the Hilbert modular case. As we learned from \cite{E-SB-T}, point (i) was also observed there
some years ago.

\begin{thm}
\label{Main Theorem HMV}(i) The smooth foliation $\mathscr{F}_{\Sigma}$
is $p$-closed if and only if $\Sigma$ is invariant under the action
of Frobenius, namely $\phi\circ\Sigma=\Sigma.$ In particular, $\mathscr{F}_{\sigma}$
is $p$-closed if and only if $f(\mathfrak{p}_{\sigma}/p)=1$ where
$\mathfrak{p}_{\sigma}$ is the prime induced by $\sigma.$

(ii) Suppose that $f(\mathfrak{p}_{\sigma}/p)\ne1.$ Then, up to a
unit, the obstruction $\kappa_{\mathscr{F}_{\sigma}}$ to $\mathscr{F}_{\sigma}$
being $p$-closed (\ref{subsec:obstruction p closed}) is equal to
the square of the $\phi\circ\sigma$-partial Hasse invariant $h_{\phi\circ\sigma}$
\cite{Go}.

(iii) Let $\mathfrak{p}$ be a prime of $L$ above $p$ and $\mathscr{F}_{\mathfrak{p}}=\mathscr{F}_{\mathbb{B}_{\mathfrak{p}}}=\oplus_{\sigma\in\mathbb{B}_{\mathfrak{p}}}\mathscr{F}_{\sigma}$
the corresponding $p$-foliation. The quotient of $M$ by the $p$-foliation
$\oplus_{\mathfrak{q}\ne\mathfrak{p}}\mathscr{F}_{\mathfrak{q}}$
may be identified with the étale component of the $\Gamma_{0}(\mathfrak{p})$-moduli
scheme $M_{0}(\mathfrak{p})$ (see details in the proof).
\end{thm}

\subsection{Preliminaries}

\subsubsection{Tate objects\label{subsec:Tate-objects}}

We begin our proof of Theorem \ref{Main Theorem HMV} by recalling
a result of Katz \cite{Ka}, who computed the effect of the Kodaira-Spencer
isomorphism on $q$-expansions.

As in \cite{Ka} (1.1.4), let $S$ be a set of $g$ linearly independent
$\mathbb{Q}$-linear forms $l_{i}:L\to\mathbb{Q}$ preserving (total)
positivity. Let $\mathfrak{a}$ and $\mathfrak{b}$ be fractional
ideals of $L$, relatively prime to $N,$ such that $\mathfrak{c}=\mathfrak{ab}^{-1}$.
Let $\mathcal{R}=W(\kappa)\otimes\mathbb{Z}((\mathfrak{ab},S))$ be
the ring defined in \cite{Ka} (1.1.7), after base change to $W(\kappa).$
Let 
\[
\mathrm{Tate}_{\mathfrak{a,b}}(q)=\mathbb{G}_{m}\otimes\mathfrak{d}^{-1}\mathfrak{a}^{-1}/q(\mathfrak{b})
\]
be the abelian scheme over $\mathcal{R}$ constructed in \cite{Ka}
(1.1.13). For fixed $\mathfrak{a}$ and $\mathfrak{b}$ it is essentially
independent of $S$. Since $\mathfrak{a}$ is relatively prime to
$N$, $\mathrm{Tate}_{\mathfrak{a,b}}(q)$ admits a \emph{canonical}
$\Gamma_{00}(N)$-level structure $\eta_{can}$ (denoted in \cite{Ka}
(1.1.16) by $i_{can}$). It also admits a canonical $\mathfrak{c}$-polarization
$\lambda_{can}$ and a canonical action $\iota_{can}$ of $\mathcal{O}_{L}.$
We thus obtain an object
\[
\underline{\mathrm{Tate}}_{\mathfrak{a,b}}(q)=(\mathrm{Tate}_{\mathfrak{a,b}}(q),\iota_{can},\lambda_{can},\eta_{can})
\]
over the ring $\mathcal{R}$ (for any choice of $S$). For the definition
of $q$-expansions at the cusp labeled by the pair $\mathfrak{(a},\mathfrak{b})$,
recalled below, we assume, in addition, that $\mathfrak{a}$ is relatively
prime to $p.$ This can always be achieved, since only the classes
of $\mathfrak{a}$ and $\mathfrak{b}$ in the strict ideal class group
$Cl^{+}(L)$ matter. 

The Lie algebra of $\mathrm{Tate}_{\mathfrak{a,b}}(q)$ is given by
a canonical isomorphism (\cite{Ka} (1.1.17))
\[
\omega_{\mathfrak{a}}:\underline{\text{\rm Lie}}=\text{\rm Lie}(\mathrm{Tate}_{\mathfrak{a,b}}(q)/\mathcal{R})\simeq \text{\rm Lie}(\mathbb{G}_{m}\otimes\mathfrak{d}^{-1}\mathfrak{a}^{-1}/\mathcal{R})=\mathfrak{d}^{-1}\mathfrak{a}^{-1}\otimes\mathcal{R},
\]
hence the Kodaira-Spencer map $\text{\rm KS}$ (\ref{eq:KS}) induces a map
\begin{equation}
\text{\rm KS}:\mathcal{T}_{\text{\rm Spec}(\mathcal{R})/W(\kappa)}\to\underline{\text{\rm Lie}}^{\otimes2}\otimes_{\mathcal{O}_{L}}\mathfrak{dc}\simeq\mathfrak{d}^{-1}\mathfrak{a}^{-1}\mathfrak{b}^{-1}\otimes\mathcal{R}.\label{eq:KS-q}
\end{equation}
The tangent bundle $\mathcal{T}_{\text{\rm Spec}(\mathcal{R})/W(\kappa)}$ is
the module of $W(\kappa)$-derivations of $\mathcal{R}$. For $\gamma\in\mathfrak{d}^{-1}\mathfrak{a}^{-1}\mathfrak{b}^{-1}$
we may consider the derivation $D(\gamma)$ of $\mathcal{R}$ (analogue
of $q\frac{d}{dq})$ given by
\[
D(\gamma)(\sum_{\alpha\in\mathfrak{ab},\,\,l_{i}(\alpha)\ge-n}a_{\alpha}q^{\alpha})=\sum_{\alpha\in\mathfrak{ab},\,\,l_{i}(\alpha)\ge-n}\text{\rm Tr}_{L/\mathbb{Q}}(\alpha\gamma)a_{\alpha}q^{\alpha}.
\]
We then have the following elegant result.
\begin{lem}
\label{lem:Katz' formula}(\cite{Ka} (1.1.20)) The image of $D(\gamma)$
under the Kodaira-Spencer map $\text{\rm KS}$ in $(\ref{eq:KS-q})$ is $\gamma\otimes1.$
\end{lem}

\subsubsection{Hilbert modular forms and partial Hasse invariants}

In this subsection we set up notation for Hilbert modular forms, and
recall some results due to the first author and to Diamond and Kassaei
on the $q$-expansions of such forms over $\kappa.$ 

By our assumption on $\kappa,$ for any $W(\kappa)$-algebra $R$
we have
\[
\mathcal{O}_{L}\otimes R\simeq\oplus_{\sigma\in\mathbb{B}}R_{\sigma}
\]
where $R_{\sigma}$ is the ring $R$ equipped with the action of $\mathcal{O}_{L}$
via $\sigma:\mathcal{O}_{L}\hookrightarrow W(\kappa)$. A \emph{weight
}is a tuple $k=(k_{\sigma})_{\sigma\in\mathbb{B}}$ with $k_{\sigma}\in\mathbb{Z}.$
We shall often write $k$ also as the element
\[
\sum_{\sigma\in\mathbb{B}}k_{\sigma}[\sigma]\in\mathbb{Z}[\mathbb{B}].
\]
The weight $k$ defines a homomorphism of algebraic groups over $W(\kappa),$
$\chi=\chi_{k}:\text{\rm Res}_{W(\kappa)}^{\mathcal{O}_{L}\otimes W(\kappa)}(\mathbb{G}_{m})\to\mathbb{G}_{m}$,
given on $R$-points by 
\[
\chi:(\mathcal{O}_{L}\otimes R)^{\times}\to R^{\times},\,\,\,\chi(a\otimes x)=\prod_{\sigma\in\mathbb{B}}(\sigma(a)x)^{k_{\sigma}}.
\]
A weight $k$, level $N$ Hilbert modular form (HMF) $f$ over $W(\kappa)$
``à la Katz'' is a rule associating to any $W(\kappa)$-algebra
$R$, any four-tuple $\underline{A}$ over $R$ as above, and any
$\mathcal{O}_{L}\otimes R$-basis $\omega$ of $\underline{\omega}_{A/R}$
(such a basis always exists locally on $R$) an element $f(\underline{A},\omega)\in R$,
which depends only on the isomorphism type of the pair $(\underline{A},\omega)$,
is compatible with base change $R\to R'$ (over $W(\kappa)$), and
satisfies
\[
f(\underline{A},\alpha\omega)=\chi(\alpha)^{-1}f(\underline{A},\omega)
\]
for $\alpha\in(\mathcal{O}_{L}\otimes R)^{\times}.$ The $(\mathfrak{a,b})$-$q$-expansion
of $f$ is the element 
\[
f(\underline{\mathrm{Tate}}_{\mathfrak{a,b}}(q),(\text{\rm Tr}_{L/\mathbb{Q}}\otimes1)\circ\omega_{\mathfrak{a}})\in\mathcal{R},
\]
where $\mathcal{R},$ $\mathfrak{a}$ and $\mathfrak{b}$ are as above.
Note that by our assumption that $\mathfrak{a}$ is relatively prime
to $p$, $\mathfrak{d}^{-1}\mathfrak{a}^{-1}\otimes\mathcal{R}=\mathfrak{d}^{-1}\otimes\mathcal{R}$.
Since $\text{\rm Tr}_{L/\mathbb{Q}}:\mathfrak{d}^{-1}\to\mathbb{Z}$ is an $\mathcal{O}_{L}$-basis
of $\mathrm{Hom}(\mathfrak{d}^{-1},\mathbb{Z}),$ 
\[
(\text{\rm Tr}_{L/\mathbb{Q}}\otimes1)\circ\omega_{\mathfrak{a}}:\underline{\text{\rm \text{\rm Lie}}}\to\mathcal{R}
\]
 is indeed an $\mathcal{O}_{L}\otimes\mathcal{R}$-basis of $\underline{\omega}_{A/\mathcal{R}}$
for $A=\mathrm{Tate}_{\mathfrak{a,b}}(q).$

For a weight $k$ denote by $\mathscr{L}_{\chi}$ the line bundle
\[
\mathscr{L}_{\chi}=\bigotimes_{\sigma\in\mathbb{B}}\mathscr{L}_{\sigma}^{\otimes k_{\sigma}}
\]
on $\mathscr{M}.$ Let $R$ be a $W(\kappa)$-algebra and $\underline{A}$
a four-tuple over $R$ as above, corresponding to a morphism $h:\text{\rm Spec}(R)\to\mathscr{M}$
over $W(\kappa)$. An $\mathcal{O}_{L}\otimes R$-basis $\omega$
of $\omega_{A/R}$ yields $R$-bases $\omega_{\sigma}$ of the line
bundles $\omega_{A/R}[\sigma]=h^{*}\mathscr{L}_{\sigma}$ for every
$\sigma\in\mathbb{B},$ hence a basis $\omega_{\chi}$ of $\bigotimes_{\sigma\in\mathbb{B}}\omega[\sigma]^{\otimes k_{\sigma}}=h^{*}\mathscr{L}_{\chi}$.
If $f$ is a weight $k$ HMF then $f(\underline{A},\omega)\cdot\omega_{\chi}$
is independent of $\omega$. We may therefore regard $f$ as a global
section of $\mathscr{L}_{\chi},$ and vice versa, any global section
of $\mathscr{L}_{\chi}$ over $\mathscr{M}$ is a HMF of weight $k$
and level $N.$

Since we assume that $g\ge2$, by the Köcher principle any HMF $f$
is automatically holomorphic at the cusps, and the $q$-expansions
of $f$ lie in
\[
\mathcal{R}\cap\{a_{0}+\sum_{\alpha\gg0}a_{\alpha}q^{\alpha}|\,a_{0},a_{\alpha}\in W(\kappa)\}.
\]

The same analysis holds if we restrict to $\kappa$-algebras $R$
rather than $W(\kappa)$-algebras, and yields a definition of HMF's
of weight $k$ and level $N$ over $\kappa$, as well as an interpretation
of such modular forms as global sections of $\mathscr{L}_{\chi}$
over $M,$ the special fiber of $\mathscr{M}.$ We also get the mod-$p$
$(\mathfrak{a,b})$-$q$-expansion of a modular form over $\kappa$
as an element of $\mathcal{R}/p\mathcal{R}$ by the same recipe. However,
in general, not every HMF over~$\kappa$ lifts to a HMF over $W(\kappa).$
The exact sequence
\[
0\to\mathscr{L}_{\chi}\overset{\times p}{\to}\mathscr{L}_{\chi}\to\mathscr{L}_{\chi}/p\mathscr{L_{\chi}}\to0
\]
shows that the obstruction to lifting a HMF over $\kappa$ lies in
$H^{1}(\mathscr{M},\mathscr{L}_{\chi}).$

Let $M_{k}(N,W(\kappa))=H^{0}(\mathscr{M},\mathscr{L}_{\chi_{k}})$
denote the space of weight $k,$ level $N$ HMF's over $W(\kappa)$
and similarly $M_{k}(N,\kappa)=H^{0}(M,\mathscr{L}_{\chi_{k}})$ the
space of weight $k$, level $N$ HMF's over $\kappa.$ The $q$-\emph{expansion
principle }says that a modular form over $W(\kappa)$ (or over $\kappa$),
all of whose $q$-expansions, for all choices of $(\mathfrak{a},\mathfrak{b})$
- corresponding to the various cusps of the Hilbert modular scheme
- vanish, is 0.

The space
\[
M_{*}(N,W(\kappa))=\bigoplus_{k\in\mathbb{Z}[\mathbb{B}]}M_{k}(N,W(\kappa))
\]
carries a natural ring structure, and is called the \emph{ring of
modular forms of level} $N$ over $W(\kappa).$ Similar terminology
applies to the ring $M_{*}(N,\kappa)$ of modular forms of level $N$
over $\kappa$. The $q$-expansion homomorphisms extend naturally
to ring homomorphisms from these rings to the rings $\mathcal{R}$
or $\mathcal{R}/p\mathcal{R}$ (depending on the choice of $\mathfrak{a}$
and $\mathfrak{b}$). However, different HMF's over $\kappa$ (of
different weights) may now have the same $q$-expansions.

An important role in the study of $q$-expansions in characteristic
$p$ is played by the $g$ \emph{partial Hasse invariants} $h_{\sigma}$
($\sigma\in\mathbb{B}).$ These are modular forms over $\kappa,$
of weights
\[
k_{\sigma}=p[\phi^{-1}\circ\sigma]-[\sigma],
\]
whose $q$-expansions at every unramified cusp is $1$. The reader is referred to Theorem~2.1 of \cite{Go} for details, but briefly the situation is as follows: 

Let $(\underline{A}, \omega)$ be an object as above, over a $\kappa$-algebra $R$, where $\underline{{\rm Lie}} \,A$ is a free $\mathcal{O}_L\otimes R$ module of rank $1$ and $\omega$ is a basis for the relative differentials of $A$ over $R$. Under the decomposition $\mathcal{O}_L\otimes R = \oplus_{\sigma} R$, where $\sigma$ runs over the homomorphisms $\mathcal{O}_L \rightarrow \kappa$,  we have a corresponding decomposition $\omega = \oplus_\sigma \omega_\sigma$. Let $\eta = \oplus_\sigma \eta_\sigma$ be the decomposition of the dual basis $\eta$ for $R^1\pi_\ast \mathcal{O}_A$ under the polarization. Then, for every $\sigma$, $h_\sigma((\underline{A}, \omega))$ is defined by the identity $h_\sigma(\underline{A}, \omega)\cdot\eta_\sigma = {\rm Fr}(\eta_{\phi^{-1}\circ \sigma})$.

\medskip 

Following \cite{DK} we define the following cones $C^\text{\rm min}\subset C^\text{\rm std}\subset C^\text{\rm hasse}$
in $\mathbb{R}[\mathbb{B}]$:
\begin{itemize}
\item $C^\text{\rm min}=\{\sum a_{\sigma}[\sigma]|\,\,pa_{\sigma}\ge a_{\phi^{-1}\circ\sigma}\,\,\forall\sigma\}$
(the \emph{minimal cone})
\item $C^\text{\rm std}=\{\sum a_{\sigma}[\sigma]|\,\,a_{\sigma}\ge0\,\,\forall\sigma\}$
(the \emph{standard cone})
\item $C^\text{\rm hasse}=\{\sum_{\sigma}a_{\sigma}(p[\phi^{-1}\circ\sigma]-[\sigma])|\,\,a_{\sigma}\ge0\}$
(the \emph{Hasse cone}).
\end{itemize}
For example, when $g=2$ and $p$ is split in $L$ all three cones
coincide. In contrast, when $p$ is inert in $L$ they look as follows:

\bigskip{}

\begin{center} \includegraphics[scale=0.35]{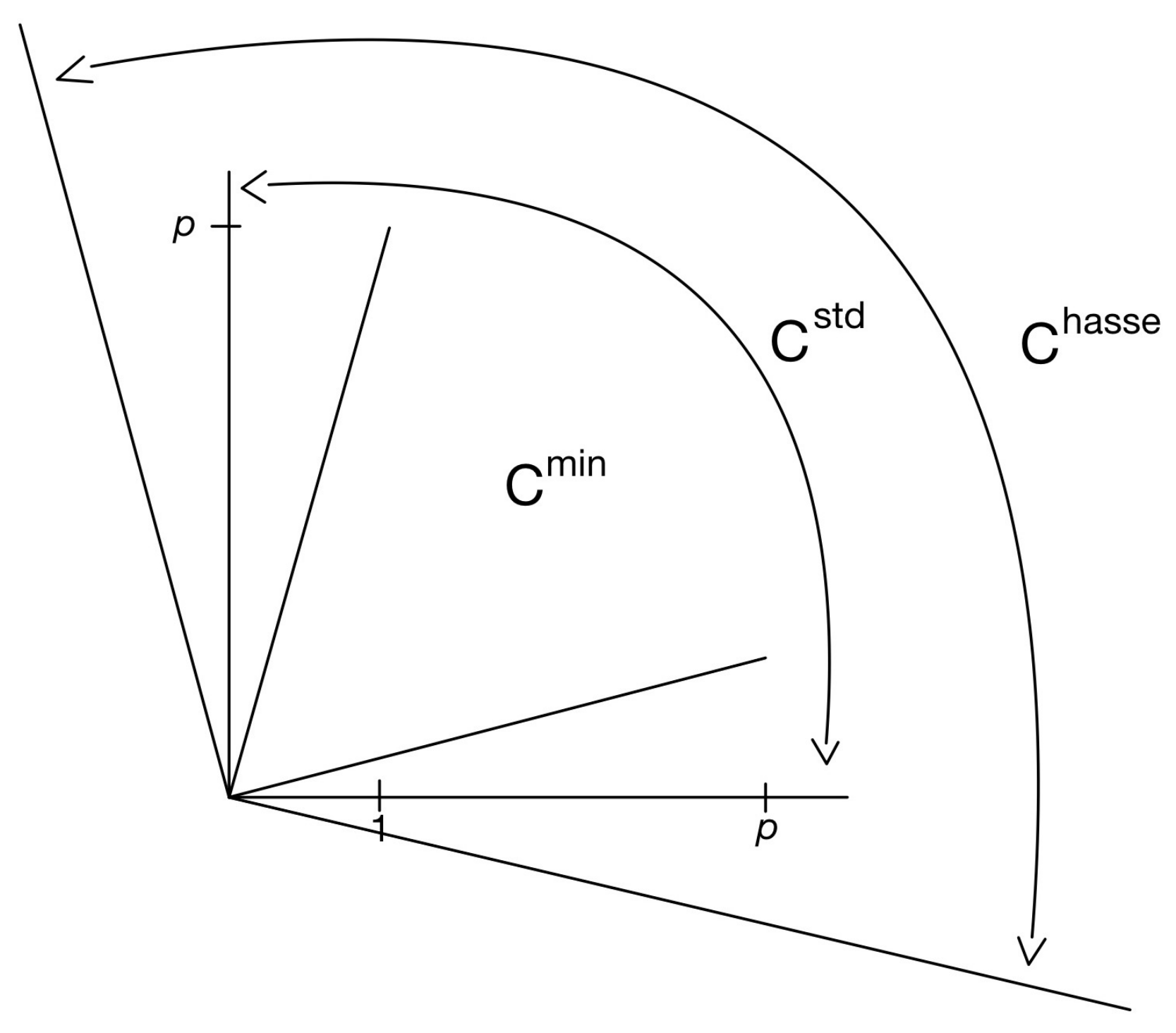}
\end{center}

\bigskip{}

Let $f\in M_{k}(N,\kappa).$ According to \cite{A-G} Proposition
8.9, any other $g\in M_{k'}(N,\kappa)$ having the same $q$-expansions
as $f$ is a product of $f$ and partial Hasse invariants raised to
integral (possibly negative) powers. Furthermore, there exists a
weight $\Phi(f),$ called the \emph{filtration} of $f$, and $g\in M_{\Phi(f)}(N,\kappa)$
having the same $q$-expansions as $f,$ such that any other $g_{1}\in M_{k'}(N,\kappa)$
with the same $q$-expansions is of the form
\[
g_{1}=g\prod_{\sigma\in\mathbb{B}}h_{\sigma}^{n_{\sigma}}
\]
for some integers $n_{\sigma}\ge0.$ These results generalize well-known
results of Serre for elliptic modular forms in characteristic $p$.
\begin{thm}
\label{thm:Diamond Kassaei on filtration of HMF}(\cite{DK}, Corollary
1.2) The filtration $\Phi(f)\in C^\text{\rm min}.$
\end{thm}

\subsubsection{Hilbert modular varieties with Iwahori level structure at $p$}

The last piece of input needed for the proof of Theorem \ref{Main Theorem HMV}
concerns $\Gamma_{0}(\mathfrak{p})$-moduli problems, for the primes
$\mathfrak{p}$ of $L$ dividing $p.$ References for the results
quoted below are \cite{G-K,Pa,St}.

Fix a prime $\mathfrak{p}$ of $L$ dividing $p$. Let $\mathscr{M}_{0}(\mathfrak{p})$ be the moduli problem over $W(\kappa)$
classifying isomorphism types of tuples $(\underline{A},H)$ where
$\underline{A}\in\mathscr{M}$ and $H$ is a finite flat $\mathcal{O}_{L}$-invariant
isotropic subgroup scheme of $A${[}$\mathfrak{p}]$ of rank $p^{f}$,
where $f=f(\mathfrak{p}/p).$ The meaning of ``isotropic'' is the
following.

Since we assumed that $\mathfrak{c}$ is relatively prime to $p$
there is a canonical isomorphism of group schemes $\mathfrak{c}\otimes_{\mathcal{O}_{L}}A[p]\simeq A[p]$
(sending $\alpha\otimes u$ to $\iota(\alpha)u$). The canonical $e_{p}$-pairing
$A[p]\otimes A^{t}[p]\to\mu_{p}$ therefore induces, via $\lambda,$
a perfect Weil pairing
\[
\left\langle .,.\right\rangle _{\lambda}:A[p]\otimes A[p]\to\mu_{p}.
\]
It restricts to a perfect pairing on $A[\mathfrak{p}]$ (since the
Rosati involution on $\mathcal{O}_{L}$ is the identity), and $H$
is a maximal isotropic subgroup scheme of $A[\mathfrak{p}].$

If $\mathfrak{c}$ is not relatively prime to $p,$ we may change
it, in the definition of $\mathscr{M},$ to $\mathfrak{c}'=\gamma\mathfrak{c}$,
where $\gamma\gg0$, so that $\mathfrak{c'}$ is now relatively prime
to $p$, and get an isomorphic moduli problem (the isomorphism depending
on $\gamma$). By ``isotropic'' we then mean that $H$ is isotropic
with respect to the Weil pairing induced by the isomorphism $\mathfrak{c}'\otimes_{\mathcal{O}_{L}}A[p]\simeq A[p]$
as above.

As before, the moduli problem $\mathscr{M}_{0}(\mathfrak{p})$ is
represented by a scheme, flat over $W(\kappa)$, which we denote by
the same letter, and the forgetful morphism is a proper morphism $\mathscr{M}_{0}(\mathfrak{p})\to\mathscr{M}.$
The scheme $\mathscr{M}_{0}(\mathfrak{p})$ is normal and Cohen-Macaulay.
We let 
\[
M_{0}(\mathfrak{p})\to M
\]
be the characteristic $p$ fiber of this morphism, over the field
$\kappa$. This morphism has been studied in detail in \cite{G-K}.
Away from the ordinary locus, it is neither finite nor flat.

Let $M^\text{\rm ord}$ be the ordinary locus of $M$, the open dense subset
where none of the partial Hasse invariants $h_{\sigma}$ vanishes.
If $k$ is an algebraically closed field containing~$\kappa$ and
$x:\text{\rm Spec}(k)\to M$ a $k$-valued point of $M,$ then $x$ lies on $M^\text{\rm ord}$
if and only if the corresponding abelian variety $A_{x}=x^{*}(A^\text{\rm univ})$
is ordinary.

Let $M_{0}(\mathfrak{p})^\text{\rm ord}$ be the open subset of $M_{0}(\mathfrak{p})$
which lies over $M^\text{\rm ord}.$ Then $M_{0}(\mathfrak{p})^\text{\rm ord}$ is the
disjoint union of two smooth varieties. The component $M_{0}(\mathfrak{p})^\text{\rm ord,m}$
(the \emph{multiplicative} component) classifies tuples $(\underline{A},H)$
where $H$ is multiplicative (its Cartier dual is étale). The forgetful
morphism is an isomorphism $M_{0}(\mathfrak{p})^\text{\rm ord,m}\simeq M^\text{\rm ord},$
its inverse given by the section $\underline{A}\mapsto(\underline{A},A[\text{\rm Fr}]\cap A[\mathfrak{p}]).$
Here $\text{\rm Fr}:A\to A^{(p)}$ is the relative Frobenius morphism (everything
over a base scheme $S$ lying over $\kappa$).

The second component $M_{0}(\mathfrak{p})^\text{\rm ord,et}$ (the \emph{étale}
component) classifies pairs $(\underline{A},H)$ where $H$ is étale.
The forgetful morphism $M_{0}(\mathfrak{p})^\text{\rm ord,et}\to M^\text{\rm ord}$
is finite flat, purely inseparable of height 1 and degree $p^{f}.$

To discuss the Atkin-Lehner map $w$ we have to bring the polarization
module~$\mathfrak{c}$ back into the picture, because $w$ will in
general change it. We therefore recall that all our constructions
depended on an auxiliary ideal~$\mathfrak{c}$ (at least on its narrow
ideal class in $Cl^{+}(L)$), and use the notation $M^{\mathfrak{c}},$
$M_{0}^{\mathfrak{c}}(\mathfrak{p})$ etc. to emphasize this dependency.
We now define
\[
w:M_{0}^{\mathfrak{c}}(\mathfrak{p})\to M_{0}^{\mathfrak{cp}}(\mathfrak{p}),\,\,\,w(\underline{A},H)=(\underline{A}/H,A[\mathfrak{p}]/H).
\]
Here $\underline{A}/H$ is the tuple $(A/H,\iota',\lambda',\eta')$
where $\iota'$ and $\eta'$ are induced by $\iota$ and $\eta$.
The $\mathfrak{cp}$-polarization $\lambda'$ of $A/H$ is obtained
as follows. Giving $\lambda$ is equivalent to giving a homomorphism
\[
\psi_{1}:\mathfrak{c}\to\mathrm{Hom}_{\mathcal{O}_{L}}(A,A^{t})_{sym}=\mathcal{P}_{A}
\]
such that $\psi_{1}(\alpha)$ is in the cone of $\mathcal{O}_{L}$-polarizations
$\mathcal{P}_{A}^{+}$ if and only if $\alpha\gg0,$ and such that
$\psi_{1}$ induces an isomorphism
\[
\mathfrak{c}\otimes_{\mathcal{O}_{L}}A\simeq A^{t}.
\]
Denoting $A/H$ by $B$ and letting $f:A\to B$ be the canonical homomorphism,
the polarization $\lambda'$ is determined by a similar homomorphism
\[
\psi_{2}:\mathfrak{cp}\to\mathrm{Hom}_{\mathcal{O}_{L}}(B,B^{t})_{sym}=\mathcal{P}_{B},
\]
defined by the relation
\[
f^{t}\circ\psi_{2}(\alpha)\circ f=\psi_{1}(\alpha)
\]
for all $\alpha\in\mathfrak{pc}\subset\mathfrak{c}.$ The \emph{existence
}of $\psi_{2}(\alpha)$ stems from the fact that $H$ is isotropic
for the pairing $A[p]\times A[p]\to\mu_{p}$ induced by $\psi_{1}(\alpha).$
Its \emph{uniqueness} is obvious, and because of this uniqueness $\psi_{2}(\cdot)$
is a homomorphism. See \cite{Pa}, 2.2. The subgroup scheme $A[\mathfrak{p}]/H\subset(A/H)[\mathfrak{p}]=\mathfrak{p}^{-1}H/H$
is finite and flat, $\mathcal{O}_{L}$-invariant, and isotropic with
respect to $\lambda'$, of rank $p^{f}.$ The tuple $(\underline{A}/H,A[\mathfrak{p}]/H)$
therefore lies on $M_{0}^{\mathfrak{cp}}(\mathfrak{p})$, and $w$
is well defined.

In the definition of $w$ we may replace $\mathfrak{cp}$ by a polarization
module $\mathfrak{c}'$ which is relatively prime to $p,$ as was
assumed for $\mathfrak{c},$ by multiplying by an appropriate $\gamma\gg0$,
keeping the strict ideal class of $\mathfrak{cp}$ unchanged.

The Atkin-Lehner map is not an involution if the class of $\mathfrak{p}$
in $Cl^{+}(L)$ is not trivial. Indeed,
\[
w^{2}(\underline{A},H)=(\underline{A}/A[\mathfrak{p}],\mathfrak{p}^{-1}H/A[\mathfrak{p}])\in M_{0}^{\mathfrak{cp}^{2}}(\mathfrak{p}).
\]
Nevertheless, as in the case of modular curves, it preserves the ordinary
locus and exchanges the ordinary étale and ordinary multiplicative
components: 
\[
w:M_{0}^{\mathfrak{c}}(\mathfrak{p})^\text{\rm ord,m}\simeq M_{0}^{\mathfrak{cp}}(\mathfrak{p})^\text{\rm ord,et},\,\,\,w:M_{0}^{\mathfrak{c}}(\mathfrak{p})^\text{\rm ord,et}\simeq M_{0}^{\mathfrak{cp}}(\mathfrak{p})^\text{\rm ord,m}.
\]

We now define $M_{0}^{\mathfrak{c}}(\mathfrak{p})^{\rm m}$ to be the
\emph{Zariski closure} of $M_{0}^{\mathfrak{c}}(\mathfrak{p})^\text{\rm ord,m}$
in $M_{0}^{\mathfrak{c}}(\mathfrak{p})$, and $M_{0}^{\mathfrak{c}}(\mathfrak{p})^\text{\rm et}$
the Zariski closure of $M_{0}^{\mathfrak{c}}(\mathfrak{p})^\text{\rm ord,et}.$
The following proposition summarizes the situation.
\begin{prop}
\label{prop:etale and multiplicative components}Both $M_{0}^{\mathfrak{c}}(\mathfrak{p})^{\rm m}$
and $M_{0}^{\mathfrak{c}}(\mathfrak{p})^\text{\rm et}$ are finite flat over
$M^{\mathfrak{c}}.$ The forgetful morphism is an isomorphism $M_{0}^{\mathfrak{c}}(\mathfrak{p})^{\rm m}\simeq M^{\mathfrak{c}},$
and $M_{0}^{\mathfrak{c}}(\mathfrak{p})^\text{\rm et}$ is purely inseparable
of height 1 and degree $p^{f}$ over $M^{\mathfrak{c}}.$ The map
$w$ induces isomorphisms
\[
w:M_{0}^{\mathfrak{c}}(\mathfrak{p})^{\rm m}\simeq M_{0}^{\mathfrak{cp}}(\mathfrak{p})^\text{\rm et},\,\,\,w:M_{0}^{\mathfrak{c}}(\mathfrak{p})^\text{\rm et}\simeq M_{0}^{\mathfrak{cp}}(\mathfrak{p})^{\rm m}.
\]
Both $M_{0}^{\mathfrak{c}}(\mathfrak{p})^{\rm m}$ and $M_{0}^{\mathfrak{c}}(\mathfrak{p})^\text{\rm et}$
are therefore smooth over $\kappa$.
\end{prop}

\begin{proof}
Taking Zariski closures, it is clear that $w$ maps $M_{0}^{\mathfrak{c}}(\mathfrak{p})^{\rm m}$
to $M_{0}^{\mathfrak{cp}}(\mathfrak{p})^\text{\rm et}$ and $M_{0}^{\mathfrak{c}}(\mathfrak{p})^\text{\rm et}$
to $M_{0}^{\mathfrak{cp}}(\mathfrak{p})^{\rm m}$. Let $h$ be the order
of $[\mathfrak{p}]$ in $Cl^{+}(L)$ and let $\gamma$ be a totally
positive element of $L$ such that $(\gamma)=\mathfrak{p}^{h}.$ Applying
$w$ successively $2h$ times we get the map
\[
M_{0}^{\mathfrak{c}}(\mathfrak{p})\to M_{0}^{\gamma^{2}\mathfrak{c}}(\mathfrak{p}),\,\,\,(\underline{A},H)\mapsto(\underline{A}/A[\gamma],\gamma^{-1}H/A[\gamma]).
\]
Denote by $\gamma_{*}:M_{0}^{\gamma^{2}\mathfrak{c}}(\mathfrak{p})\to M_{0}^{\mathfrak{c}}(\mathfrak{p})$
the isomorphism
\[
(A',\iota',\lambda',\eta';H')\mapsto(A',\iota',\lambda'\circ(\gamma^{2}\otimes1),\eta'\circ(\gamma\otimes1);H').
\]
Identifying $A/A[\gamma]$ with $A$ under $\gamma$ (carrying $\gamma^{-1}H/A[\gamma]$
to $H$) it is easy to check that $\gamma_{*}$ maps the tuple $(\underline{A}/A[\gamma],\gamma^{-1}H/A[\gamma])\in M_{0}^{\gamma^{2}\mathfrak{c}}(\mathfrak{p})$
back to $(\underline{A},H)\in M_{0}^{\mathfrak{c}}(\mathfrak{p}).$
We see that
\[
\gamma_{*}\circ w^{2h}
\]
is the identity, hence $w:M_{0}^{\mathfrak{c}}(\mathfrak{p})\to M_{0}^{\mathfrak{cp}}(\mathfrak{p})$
is an isomorphism. 

The multiplicative component maps isomorphically to $M^{\mathfrak{c}}$
since the section $\underline{A}\mapsto(\underline{A},A[\text{\rm Fr}]\cap A[\mathfrak{p}])$
extends from the ordinary part to all of $M^{\mathfrak{c}}$ and must
map it to $M_{0}^{\mathfrak{c}}(\mathfrak{p})^{\rm m}$ by continuity.
It follows that $M_{0}^{\mathfrak{c}}(\mathfrak{p})^{\rm m}$ is smooth.
Since $M_{0}^{\mathfrak{c}}(\mathfrak{p})^\text{\rm et}$ is isomorphic to
$M_{0}^{\mathfrak{cp}}(\mathfrak{p})^{\rm m}$, it is also smooth. 

It remains to prove that the morphism $\pi:M_{0}^{\mathfrak{c}}(\mathfrak{p})^\text{\rm et}\to M^{\mathfrak{c}}$
is finite flat, purely inseparable of height 1 and degree $p^{f}$.
Consider $(M_{0}^{\mathfrak{c}}(\mathfrak{p})^\text{\rm et})^{(p)}=\kappa\times_{\phi,\kappa}M_{0}^{\mathfrak{c}}(\mathfrak{p})^\text{\rm et}$,
which we canonically identify with $M_{0}^{\mathfrak{c}}(\mathfrak{p})^\text{\rm et}$,
since it is actually defined over $\mathbb{F}_{p},$ and the relative
Frobenius morphism (we use $Y$ as a shorthand for $M_{0}^{\mathfrak{c}}(\mathfrak{p})^\text{\rm et}$)
\[
\text{\rm Fr}_{Y/\kappa}:M_{0}^{\mathfrak{c}}(\mathfrak{p})^\text{\rm et}\to M_{0}^{\mathfrak{c}}(\mathfrak{p})^\text{\rm et}.
\]
As $M_{0}^{\mathfrak{c}}(\mathfrak{p})^\text{\rm et}$ is smooth, $\text{\rm Fr}_{Y/\kappa}$
is finite and flat of degree $p^{g}$. We claim that there is a morphism
$\theta:M^{\mathfrak{c}}\to M_{0}^{\mathfrak{c}}(\mathfrak{p})^\text{\rm et}$
such that $\text{\rm Fr}_{Y/\kappa}=\theta\circ\pi$. This will force both $\theta$
and $\pi$ to be finite flat and purely inseparable of height 1, as
finite morphisms between regular schemes are flat (``miracle flatness'',
\cite{Stacks}, Lemma 10.128.1). As the degree of $\pi$ over the
ordinary locus is $p^{f},$ this is its degree everywhere, and as
a by-product we get that the degree of $\theta$ is $p^{g-f}.$

We define $\theta$ in the language of moduli problems. Let $\sigma:M^{\mathfrak{c}}\simeq M_{0}^{\mathfrak{c}}(\mathfrak{p})^{\rm m}$
be the section $\underline{A}\mapsto(\underline{A},A[\text{\rm Fr}]\cap A[\mathfrak{p}])$
described before. Then $w\circ\sigma:M^{\mathfrak{c}}\simeq M_{0}^{\mathfrak{cp}}(\mathfrak{p})^\text{\rm et}$
is an isomorphism. Let $\mathfrak{p}'=\prod_{\mathfrak{q}\ne\mathfrak{p}}\mathfrak{q}$
be the product of the primes of $L$ dividing $p$ that are different
from $\mathfrak{p}.$ Let
\[
\theta':M_{0}^{\mathfrak{cp}}(\mathfrak{p})\to M_{0}^{\mathfrak{cpp'}}(\mathfrak{p})=M_{0}^{\mathfrak{c}p}(\mathfrak{p})
\]
be the map 
\[
\theta':(\underline{A},H)\mapsto(\underline{A}/A[\text{\rm Fr}]\cap A[\mathfrak{p}'],H\mod A[\text{\rm Fr}]\cap A[\mathfrak{p}']).
\]
 As $H\subset A[\mathfrak{p}]$, and $\mathfrak{p}$ and $\mathfrak{p}'$
are relatively prime, this map is well-defined and in fact sends $M_{0}^{\mathfrak{cp}}(\mathfrak{p})^\text{\rm et}$
to $M_{0}^{\mathfrak{c}p}(\mathfrak{p})^\text{\rm et}$ (it is enough to check
this on the ordinary locus). We let
\[
\theta=\theta'\circ w\circ\sigma:M^{\mathfrak{c}}\to M_{0}^{\mathfrak{c}p}(\mathfrak{p})^\text{\rm et}\simeq M_{0}^{\mathfrak{c}}(\mathfrak{p})^\text{\rm et}.
\]
In the last step, we have used the fact that $\mathfrak{pp'}=(p)$
is principal to identify $M_{0}^{\mathfrak{c}p}(\mathfrak{p})^\text{\rm et}\simeq M_{0}^{\mathfrak{c}}(\mathfrak{p})^\text{\rm et}$
sending $(A,\iota,\lambda,\eta;H)$ to $(A,\iota,\tilde{\lambda},\eta;H)$
where if $\lambda:\mathfrak{c}p\otimes_{\mathcal{O}_{L}}A\simeq A^{t}$
is a $\mathfrak{c}p$-polarization, $\tilde{\lambda}:\mathfrak{c}\otimes_{\mathcal{O}_{L}}A\simeq A^{t}$
is the $\mathfrak{c}$-polarization given by $\tilde{\lambda}=\lambda\circ(p\otimes1).$
It follows that
\begin{equation}\label{equation for theta}
\theta(\underline{A},H)=(\underline{A}/A[\text{\rm Fr}],A[\mathfrak{p}]\mod A[\text{\rm Fr}]).
\end{equation}

To conclude the proof we have to show that for $(\underline{A},H)\in M_{0}^{\mathfrak{c}}(\mathfrak{p})^\text{\rm et}$
\[
(\underline{A}/A[\text{\rm Fr}],A[\mathfrak{p}]\mod A[\text{\rm Fr}])=\text{\rm Fr}_{Y/\kappa}(\underline{A},H).
\]
It is enough to verify this for $(\underline{A},H)\in M_{0}^{\mathfrak{c}}(\mathfrak{p})^\text{\rm ord,et},$
as two morphisms that coincide on a dense open set, are equal. Assume
that $(\underline{A},H)\in M_{0}^{\mathfrak{c}}(\mathfrak{p})^\text{\rm ord,et}(k)$
is a $k$-valued point, for a $\kappa$-algebra $k.$ Then $\text{\rm Fr}=\text{\rm Fr}_{A/k}$
is the relative Frobenius of $A$ over $k$ (not to be confused with
$\text{\rm Fr}_{Y/\kappa}$ on the right hand side!). But $\underline{A}/A[\text{\rm Fr}]\simeq\underline{A}^{(p)/k}$
(base change with respect to the absolute Frobenius of $k$), and
since $H$ is the unique étale subgroup of $A[\mathfrak{p}]$ of order
$p^{f},$ and $A[\mathfrak{p}]\mod A[\text{\rm Fr}]$ is the unique étale subgroup
of $A^{(p)/k}[\mathfrak{p}]$ of order $p^{f}$, we must also have
$A[\mathfrak{p}]\mod A[\text{\rm Fr}]\simeq H^{(p)/k}.$ 

Finally, let us explain the equality $(\underline{A}^{(p)/k},H^{(p)/k})=\text{\rm Fr}_{Y/\kappa}(\underline{A},H)$.
Intuitively, ``the moduli of the object obtained by Frobenius base
change is the Frobenius base change of the original moduli''. However,
as there are two \emph{different }relative Frobenii involved, care
must be taken. Let $(\underline{A},H)\in M_{0}^{\mathfrak{c}}(\mathfrak{p})^\text{\rm ord,et}(k)$
correspond to the point $x:\text{\rm Spec}(k)\to M_{0}^{\mathfrak{c}}(\mathfrak{p})^\text{\rm ord,et}$
(over $\kappa$). By the functoriality of Frobenius, we have the following
commutative diagram, where we have substituted $Y=M_{0}^{\mathfrak{c}}(\mathfrak{p})^\text{\rm ord,et}$:
\[\xymatrix@C=1.5cm{\text{\rm Spec}(k)\ar[r]^x\ar[d]_{\text{\rm Fr}_{k\!/\! \kappa}} & Y\ar[d]^{\text{\rm Fr}_{Y\!/\! \kappa}}\\
\text{\rm Spec}(k^{(p)\!/\! \kappa}) \ar[r]^{x^{(p)\!/\! \kappa}}\ar[d] & Y^{(p)\!/\! \kappa}\ar[d]\\
\text{\rm Spec}(k) \ar[r]^x & Y
}\]
Here the vertical unmarked arrows are the base change maps with respect
to the absolute Frobenius of $\kappa.$ The composition of the two
vertical arrows on the right is the absolute Frobenius $\Phi_{Y}$
of $Y$, and the composition of the two arrows on the left is the
absolute Frobenius $\Phi_{k}$ of $\text{\rm Spec}(k)$. We therefore have
\[
(\underline{A}^{(p)/k},H^{(p)/k})=\Phi_{k}^{*}(x^{*}(\underline{A}^\text{\rm univ},H^\text{\rm univ}))=(\Phi_{Y}\circ x)^{*}(\underline{A}^\text{\rm univ},H^\text{\rm univ})=\text{\rm Fr}_{Y/\kappa}(\underline{A},H).
\]
(In the last step, we identified $Y^{(p)/\kappa}$ with $Y$, as it
is defined over $\mathbb{F}_{p},$ hence we may identify $\text{\rm Fr}_{Y/\kappa}(\underline{A},H)$
with $\Phi_{Y}(\underline{A},H).$)
\end{proof}
\begin{cor}
The smooth $\kappa$-variety $M_{0}^{\mathfrak{c}}(\mathfrak{p})^\text{\rm et}$
is the quotient of $M^{\mathfrak{c}}$ by a smooth $p$-foliation
of rank $g-f$.
\end{cor}

\begin{proof}
The finite flat morphism $\theta$ defined in the proof of the previous
proposition is purely inseparable of height 1 and degree $p^{g-f}.$
The corollary now follows easily from Theorem \ref{thm:Quotient by foliation}.
\end{proof}

\subsection{Proof of Theorem \ref{Main Theorem HMV}}

\subsubsection{Proof of parts (i) and (ii)}

Let $\Sigma\subset\mathbb{B}.$ We have seen that the foliation $\mathscr{F}_{\Sigma}$
is involutive. The obstruction $\kappa_{\mathscr{F}_{\Sigma}}$ to
it being $p$-closed (\S\ref{subsec:obstruction p closed}) lies in
\[
\mathrm{Hom}(\Phi_{M}^{*}\mathscr{F}_{\Sigma},\mathcal{T}_{M/\kappa}/\mathscr{F}_{\Sigma})\simeq\bigoplus_{\sigma\in\Sigma}\bigoplus_{\tau\notin\Sigma}\mathrm{Hom}(\mathscr{L}_{\sigma}^{-2p},\mathscr{L}_{\tau}^{-2})
\]
\[
=\bigoplus_{\sigma\in\Sigma}\bigoplus_{\tau\notin\Sigma}\Gamma(M,\mathscr{L}_{\sigma}^{2p}\otimes\mathscr{L}_{\tau}^{-2}).
\]
Here we use the fact that under the absolute Frobenius
$\Phi_{M}\colon M\to M$ the pullback of a \emph{line} bundle $\mathscr{L}$ is isomorphic
to $\mathscr{L}^{p}.$ Indeed, in characteristic $p$ the map $f\otimes s\mapsto fs^{\otimes p}$
is an isomorphism $\mathcal{O}_{M}\otimes_{\phi,\mathcal{O}_{M}}\mathscr{L}\simeq\mathscr{L}^{p}.$
\begin{lem}
\label{lem:Nonexistence of HMF of certain weights}Let $\sigma\ne\tau$.
We have
\[
\Gamma(M,\mathscr{L}_{\sigma}^{2p}\otimes\mathscr{L}_{\tau}^{-2})=\begin{cases}
\begin{array}{c}
0\\
\kappa\cdot h_{\tau}^{2}
\end{array} & \begin{array}{c}
\tau\ne\phi\circ\sigma;\\
\tau=\phi\circ\sigma.
\end{array}\end{cases}
\]
\end{lem}

\begin{proof}
Let $h$ be a non-zero Hilbert modular form on $M$ of weight $2p[\sigma]-2[\tau].$
According to Theorem \ref{thm:Diamond Kassaei on filtration of HMF}
there exist integers $a_{\beta}\ge0$ such that
\[
\Phi(h)=2p[\sigma]-2[\tau]-\sum_{\beta\in\mathbb{B}}a_{\beta}(p[\beta]-[\phi\circ\beta])\in C^\text{\rm min}\subset C^\text{\rm std}
\]
(we end up using only the weaker result that the left hand side lies
in $C^\text{\rm std}$). If $\tau$ is not in the $\phi$-orbit of $\sigma$
and $\tau\in\mathbb{B}_{\mathfrak{p}}$ we get, upon summing the coefficients
of $\beta\in\mathbb{B}_{\mathfrak{p}}$ in $\Phi(h)$ that
\[
-2-(p-1)\sum_{\beta\in\mathbb{B}_{\mathfrak{p}}}a_{\beta}\ge0,
\]
a clear contradiction. It follows that $\tau=\phi^{i}\sigma$ for
some $1\le i\le f-1$ where $f\ge2$ is the length of the $\phi$-orbit
of $\sigma$ (if $\sigma\in\mathbb{B}_{\mathfrak{p}},$ then $f=f(\mathfrak{p}/p)$).

Labelling the $\beta\in\mathbb{B}_{\mathfrak{p}}$ by $0,1,\dots,f-1$
so that $\phi\circ[i]=[i+1\mod f]$, and assuming, without loss of
generality, that $[\sigma]=[0]$ and $[\tau]=[i]$ for some $1\le i\le f-1$
we have
\[
2p[0]-2[i]-\sum_{j=0}^{f-1}a_{j}(p[j]-[j+1\mod f])=\sum_{j=0}^{f-1}k_{j}[j],
\]
where $a_{j}\ge0$ and $k_{j}\ge0.$

Summing over the coefficients we get $2p-2-(p-1)\sum a_{j}=\sum k_{j}\ge0,$
hence $\sum a_{j}=0,1$ or 2. We can not have $a_{j}=0$ for all $j,$
since $\Phi(h)\in C^\text{\rm std}$. If $a_{m}=1$ and all the other $a_{j}=0$
we again reach a contradiction since the coefficient of $[i]$ comes
out negative, no matter what $m$ is. There remains the case where
$\sum a_{j}=2.$ In this case all the $k_{j}=0$. Looking at the coefficient
of $[0]$ we must have $2p-pa_{0}+a_{f-1}=0.$ This forces $a_{0}=2$
hence all the other $a_{j}=0$ and $i=1$.

This means that $\tau=\phi\circ\sigma$ and $hh_{\tau}^{-2}$ has
weight 0, i.e. is a constant in $\kappa,$ proving the lemma.
\end{proof}
The lemma implies that if $\Sigma$ is invariant under $\phi$ (i.e.
is a union of certain $\mathbb{B}_{\mathfrak{p}}$'s for the primes
$\mathfrak{p}$ above $p$) then $\kappa_{\mathscr{F}_{\Sigma}}$
vanishes, and $\mathscr{F}_{\Sigma}$ is therefore $p$-closed. To
prove the converse, completing the proof of part (i) of Theorem \ref{Main Theorem HMV},
it is enough to show that when $\Sigma=\{\sigma\}$ the obstruction
$\kappa_{\mathscr{F}_{\sigma}}$ is a \emph{non-zero} multiple of
$h_{\phi\circ\sigma}^{2}.$ This will establish, at the same time,
claim (ii) of the Theorem.

To this end we use $q$-expansions. Let $\mathcal{R}$ be one of the
rings associated with the cusps as in \ref{subsec:Tate-objects},
and $R=\mathcal{R}/p\mathcal{R}.$ The pull back of the tangent bundle
of $M$ to $\text{\rm Spec}(R)$ is identified with the Lie algebra $\text{\rm Der}(R/\kappa)$
and the Kodaira-Spencer isomorphism yields the isomorphism (\ref{eq:KS-q})
\[
\text{\rm KS}:\mathcal{T}_{\text{\rm Spec}(R)/\kappa}=\text{\rm Der}(R/\kappa)\simeq\mathcal{O}_{L}\otimes R=\bigoplus_{\sigma\in\mathbb{B}}R_{\sigma}.
\]
Let $\{e_{\sigma}\}$ be the idempotents of $\mathcal{O}_{L}\otimes R$
corresponding to this decomposition. Then $\text{\rm KS}:\mathscr{L}_{\sigma}^{-2}|_{\text{\rm Spec}(R)}\simeq Re_{\sigma}=R_{\sigma}$.
The ring $\mathcal{O}_{L}\otimes R$ has an endomorphism
\[
\varphi(\alpha\otimes r)=\alpha\otimes r^{p}
\]
and $\varphi(e_{\sigma})=e_{\phi\circ\sigma}.$ Indeed,
\[
\alpha\otimes1\cdot\varphi(e_{\sigma})=\varphi(\alpha\otimes1\cdot e_{\sigma})=\varphi(1\otimes\sigma(\alpha)\cdot e_{\sigma})=1\otimes\sigma(\alpha)^{p}\cdot\varphi(e_{\sigma})=1\otimes(\phi\circ\sigma)(\alpha)\cdot\varphi(e_{\sigma})
\]
so $\varphi(e_{\sigma}),$ being an idempotent in $R_{\phi\circ\sigma},$
must equal $e_{\phi\circ\sigma}$. 

Let $\xi_{\sigma}\in\mathscr{L}_{\sigma}^{-2}$ be the derivation
mapping to $e_{\sigma}$ under $\text{\rm KS}$. If $e_{\sigma}=\sum_{j}\gamma_{j}\otimes r_{j}$
($\gamma_{j}\in\mathcal{O}_{L},$~$r_{j}\in R$) then by Katz' formula
(Lemma \ref{lem:Katz' formula})
\[
\xi_{\sigma}(\sum_{\alpha}a_{\alpha}q^{\alpha})=\sum_{\alpha}a_{\alpha}(\sum_{j}r_{j}\text{\rm Tr}_{L/\mathbb{Q}}(\alpha\gamma_{j}))q^{\alpha}.
\]
It follows that when we iterate $\xi_{\sigma}$ $p$ times we get
the derivation
\[
\xi_{\sigma}^{p}(\sum_{\alpha}a_{\alpha}q^{\alpha})=\sum_{\alpha}a_{\alpha}(\sum_{j}r_{j}^{p}\text{\rm Tr}_{L/\mathbb{Q}}(\alpha\gamma_{j}))q^{\alpha},
\]
i.e. the derivation corresponding to $\varphi(e_{\sigma})=e_{\phi\circ\sigma}.$
We conclude that $\xi_{\sigma}^{p}=\xi_{\phi\circ\sigma}\ne0$, hence
$\kappa_{\mathscr{F}_{\sigma}}$ is a non-zero section of $\mathscr{L}_{\phi\circ\sigma}^{-2}\otimes\mathscr{L}_{\sigma}^{2p}.$
\begin{rem*}
The reader might have noticed that the same $q$-expansion computation
can be used to give an alternative proof of all of (i) and (ii) in
Theorem \ref{Main Theorem HMV}. However, we found Lemma \ref{lem:Nonexistence of HMF of certain weights}
and the relation to the main result of \cite{DK} of independent interest,
especially considering the extension of our results to other PEL Shimura
varieties.
\end{rem*}

\subsubsection{Proof of (iii)}

We now turn to the last part of the theorem, identifying the quotient
of $M$ by the smooth $p$-foliation
\[
\mathscr{G}=\bigoplus_{\mathfrak{q}\ne\mathfrak{p}}\mathscr{F}_{\mathfrak{q}}
\]
with the purely inseparable, finite flat map $\theta:M\to M_{0}(\mathfrak{p})^\text{\rm et}$
constructed in Proposition \ref{prop:etale and multiplicative components}.
Note that the rank of $\mathscr{G}$ is $g-f,$ while $\text{\rm deg}(\theta)=p^{g-f}.$
Since $M_{0}(\mathfrak{p})^\text{\rm et}$ is the quotient of $M$ by the smooth
$p$-foliation $\ker(d\theta)$, and the latter is also of rank $g-f,$
it is enough to prove that $d\theta$ annihilates $\mathscr{G}$ to
conclude that
\[
\ker(d\theta)=\mathscr{G},
\]
thus proving part (iii) of Theorem \ref{Main Theorem HMV}. 

To simplify the notation write $N=M_{0}(\mathfrak{p})^\text{\rm et},$ let
$k$ be an algebraically closed field containing $\kappa,$ $x\in M(k)$
and $y=\theta(x)\in N(k).$ Let $\underline{A}$ be the tuple parametrized
by~$x.$ Then, by (\ref{equation for theta}), $y$ parametrizes the tuple $(\underline{A}^{(p)},\text{\rm Fr}(A[\mathfrak{p}]))$
where $\text{\rm Fr}=\text{\rm Fr}_{A/k}.$ 

Write $TM$ for the tangent bundle $\mathcal{T}_{M/\kappa}$ and $T_{x}M$
for its fiber at $x,$ the tangent space to $M$ at $x$. Similar
meanings are attached to the symbols $TN$ and $T_{y}N.$ Let $k[\epsilon]$
be the ring of dual numbers, $\epsilon^{2}=0.$ In terms of the moduli
problem,
\[
T_{x}M=\{\underline{\widetilde{A}}\in M(k[\epsilon])|\,\underline{\widetilde{A}}\mod\epsilon=\underline{A}\}/\simeq
\]
and its origin is the ``constant'' tuple $\text{\rm Spec}(k[\epsilon])\times_{\text{\rm Spec}(k)}\underline{A}.$
Similarly,
\[
T_{y}N=\{(\underline{\widetilde{B}},\widetilde{H})\in N(k[\epsilon])|\,(\underline{\widetilde{B}},\widetilde{H})\mod\epsilon=(\underline{A}^{(p)},\text{\rm Fr}(A[\mathfrak{p}]))\}/\simeq
\]
and its origin is the ``constant'' tuple $\text{\rm Spec}(k[\epsilon])\times_{\text{\rm Spec}(k)}(\underline{A}^{(p)},\text{\rm Fr}(A[\mathfrak{p}])).$

Let $\widetilde{x}=\underline{\widetilde{A}}\in T_{x}M$ be a tangent
vector at $x$. In terms of moduli problems
\[
d\theta(\widetilde{x})=\theta(\underline{\widetilde{A}})=(\underline{\widetilde{A}}^{(p)},\text{\rm Fr}(\widetilde{A}[\mathfrak{p}]))\in N(k[\epsilon])
\]
where now $\widetilde{A}^{(p)}=\widetilde{A}^{(p)/k[\epsilon]}$ is
the base change of $\widetilde{A}$ with respect to the raising-to-power
$p$ homomorphism $\phi_{k[\epsilon]}\colon k[\epsilon]\to k[\epsilon],$
the PEL structure $\iota,\lambda,\eta$ accompanying $\widetilde{A}$
in the definition of $\underline{\widetilde{A}}$ undergoes the same
base change, and $\text{\rm Fr}=\text{\rm Fr}_{\widetilde{A}/k[\epsilon]}:\widetilde{A}\to\widetilde{A}^{(p)}$
is the relative Frobenius of $\widetilde{A}$ over $k[\epsilon]$. 

We have to show that if $\widetilde{x}\in\mathscr{G}_{x}\subset T_{x}M$
then $d\theta(\widetilde{x})=0$, namely that the tuple $(\underline{\widetilde{A}}^{(p)},\text{\rm Fr}(\widetilde{A}[\mathfrak{p}]))$
is \emph{constant }along $\text{\rm Spec}(k[\epsilon]).$ That $\underline{\widetilde{A}}^{(p)}$
is constant along $\text{\rm Spec}(k[\epsilon])$ is always true, regardless
of whether $\widetilde{x}\in\mathscr{G}_{x}$ or not, simply becasue
$\phi_{k[\epsilon]}$ factors as the projection modulo $\epsilon$,
$k[\epsilon]\twoheadrightarrow k,$ followed by $\phi_{k}$, and then
by the inclusion $k\hookrightarrow k[\epsilon]:$
\[
\phi_{k[\epsilon]}:k[\epsilon]\twoheadrightarrow k\overset{\phi_{k}}{\to}k\hookrightarrow k[\epsilon].
\]
It all boils down to the identity $(a+b\epsilon)^{p}=a^{p}.$

Suppose $\widetilde{x}\in\mathscr{G}_{x}$ and let us show that $\text{\rm Fr}(\widetilde{A}[\mathfrak{p}])$
is also constant. Recall the local models $M^\text{\rm loc}$, $N^\text{\rm loc}$ of $M$ and $N$ constructed
in \cite{Pa}, 3.3. Let $R$ be any $\kappa$-algebra. Let $W=(\mathcal{O}_{L}\otimes R)^{2}$
with the induced $\mathcal{O}_{L}$ action and the perfect alternate
pairing ($e_{1},e_{2}$ is the standard basis)
\[
\left\langle a\otimes t\cdot e_{1},b\otimes s\cdot e_{2}\right\rangle =-\left\langle b\otimes s\cdot e_{2},a\otimes t\cdot e_{1}\right\rangle =\text{\rm Tr}_{L/\mathbb{Q}}(ab)ts\in R,
\]
\[
\left\langle a\otimes t\cdot e_{i},b\otimes s\cdot e_{i}\right\rangle =0\,\quad(i=1,2).
\]
Then $M^\text{\rm loc}(R)$ is the set of rank-1 $\mathcal{O}_{L}\otimes R$
local direct summands $\omega\subset W$ which are totally isotropic
(equal to their own annihilator) under $\left\langle ,\right\rangle $.

Similarly, the $R$-points of $N^\text{\rm loc}$ are given by the following
data. Fix $a\in\mathcal{O}_{L}$ such that $\mathfrak{p}=(a,p)$ but
$a\equiv1\mod\mathfrak{q}$ for every prime $\mathfrak{q}\ne\mathfrak{p}$
above $p$. Let $u:W\to W$ be the $\mathcal{O}_{L}\otimes R$-linear
map sending $e_{1}$ to $e_{1}$ and $e_{2}$ to $a\otimes1\cdot e_{2}$.
Equivalently, if we decompose $W=\oplus_{\mathfrak{q}|p}W[\mathfrak{q}],$
$w=\sum_{\mathfrak{q}|p}w(\mathfrak{q})$, then $u$ sends $e_{2}(\mathfrak{p})$
to $0$, but every $e_{2}(\mathfrak{q})$ for $\mathfrak{q}\ne\mathfrak{p}$
to itself. Then $N^\text{\rm loc}(R)$ is the set of pairs $(\omega,\omega')$
of rank-1 $\mathcal{O}_{L}\otimes R$ local direct summands $\omega,\omega'\subset W$
which are totally isotropic such that $u(\omega)\subset\omega'$.

The scheme $N^\text{\rm loc}$ is a closed subscheme of a product of two Grassmannians,
and its projection to the first factor is $M^\text{\rm loc}$.

We shall use the $R$-points of the local models with $R=k$ or $k[\epsilon]$
to study the map $d\theta$ between $T_{x}M$ and $T_{y}N.$ Let $W=(\mathcal{O}_{L}\otimes k)^{2}$,
$\widetilde{W}=(\mathcal{O}_{L}\otimes k[\epsilon])^{2},$ and suppose
$x=\underline{A}$ corresponds to $\xi=(\omega\subset W)\in M^\text{\rm loc}(k).$
We may fix an identification $W=H_{dR}^{1}(A/k)$ so that $\omega=\omega_{A}=H^{0}(A,\Omega_{A/k}^{1}).$

Let $\widetilde{x}=\underline{\widetilde{A}}\in M(k[\epsilon])$ map
to $x$ modulo $\epsilon$. Let $\widetilde{\xi}=(\widetilde{\omega}\subset\widetilde{W})\in M^\text{\rm loc}(k[\epsilon])$
correspond to $\widetilde{x}$ under the isomorphism $T_{x}M\simeq T_{\xi}M^\text{\rm loc},$
i.e. $\widetilde{\omega}\mod\epsilon=\omega.$ Fix $\alpha\in\omega$
and suppose that $\widetilde{\alpha}=\alpha+\epsilon\beta$ ($\beta\in W)$
is an element of $\widetilde{\omega}$ mapping to $\alpha$ modulo
$\epsilon$. Since $\widetilde{\alpha}$ is uniquely determined modulo
$\epsilon\widetilde{\omega}=\epsilon\omega$, the image $\overline{\beta}$
of $\beta$ in
\[
W/\omega=H^{1}(A,\mathcal{O})=\text{\rm Lie}(A^{t})=\omega_{A^{t}}^{\vee}
\]
is well-defined. We have thus associated to $\widetilde{x}$ a map
from $\omega=\omega_{A}$ to $\omega_{A^{t}}^{\vee}.$ Using the polarization
$\lambda$ we view this as a map $\text{\rm KS}(\widetilde{x})$ from $\omega_{A}$
to $\omega_{A}^{\vee}.$ It is straightforward to prove that $\widetilde{\omega}$
is totally isotropic if and only if $\text{\rm KS}(\widetilde{x})$ is symmetric,
i.e. $\text{\rm KS}(\widetilde{x})^{\vee}=\text{\rm KS}(\widetilde{x}).$ The following
is a reformulation of the Kodaira-Spencer isomorphism:
\begin{prop*}
The map $\text{\rm KS}$ is an isomorphism from $T_{x}M$ onto the space of symmetric
$\mathcal{O}_{L}\otimes k$-homomorphisms from $\omega_{A}$ to $\omega_{A}^{\vee}.$
The fiber $\mathscr{G}_{x}$ of the foliation $\mathscr{G}$ consists
of those tangent vectors $\widetilde{x}$ for which $\text{\rm KS}(\widetilde{x})$
annihilates $\omega_{A}[\mathfrak{p}].$
\end{prop*}
In the decomposition $\mathcal{O}_{L}\otimes k=\bigoplus_{\sigma\in\mathbb{B}}k_{\sigma}$,
$\mathfrak{p}\otimes k$ is the ideal of all $(\alpha_{\sigma})$
for which $\alpha_{\sigma}=0$ whenever $\sigma\in\mathbb{B}_{\mathfrak{p}}.$
We therefore have
\[
\omega[\mathfrak{p}]=\bigoplus_{\sigma\in\mathbb{B}_{\mathfrak{p}}}\omega[\sigma]
\]
($\omega[\mathfrak{p}]$ denotes the kernel of $\mathfrak{p,}$ $\omega[\sigma]$
denotes the $\sigma$-isotypical component of $\omega$), and recall
that each $\omega[\sigma]$ is one-dimensional over $k$. We conclude
that if $\widetilde{x}\in\mathscr{G}_{x}$ then whenever $\alpha\in\omega$
and $\widetilde{\alpha}=\alpha+\epsilon\beta\in\widetilde{\omega}$
maps to it modulo $\epsilon,$ then for any $\sigma\in\mathbb{B}_{\mathfrak{p}}$
the $\sigma$-component of $\beta$ is proportional to the $\sigma$-component
of $\alpha.$ In other words, $\widetilde{\omega}[\mathfrak{p}]=k[\epsilon]\otimes_{k}\omega[\mathfrak{p}].$

By the crystalline deformation theory of Grothendieck (equivalently,
by the local model), the group scheme $\widetilde{A}[p]$ over $k[\epsilon]$
is completely determined by the lifting~$\widetilde{\omega}$ of $\omega$
to $k[\epsilon]$. In the notation used above $\widetilde{\omega}$
determines the point $\widetilde{\xi}=(\widetilde{\omega}\subset\widetilde{W})\in M^\text{\rm loc}(k[\epsilon]),$
hence the deformation $\widetilde{x}=\underline{\widetilde{A}}\in M(k[\epsilon]),$
and therefore $\widetilde{A}[p].$ Furthermore, the subgroup scheme
$\widetilde{A}[\mathfrak{p}]$ is constant along $k[\epsilon]$ if
and only if $\widetilde{\omega}[\mathfrak{p}]=k[\epsilon]\otimes_{k}\omega[\mathfrak{p}].$
We have therefore proved that if $\widetilde{x}\in\mathscr{G}_{x},$
$\widetilde{A}[\mathfrak{p}]$ is constant along $k[\epsilon].$ With
it, so are $\widetilde{A}[\mathfrak{p}][\text{\rm Fr}]$ and the quotient $\text{\rm Fr}(\widetilde{A}[\mathfrak{p}])$.

\subsection{Integral varieties} \label{subsec integral varieties}

Our goal in this section is to prove the following theorem.
\begin{thm}
\label{thm:HBMV integral varieties}(i) Assume that $\emptyset\subsetneqq\Sigma\subsetneqq\mathbb{B}$.
Then the foliation $\mathscr{F}_{\Sigma}$ does not have any (algebraic)
integral variety in the generic fiber $\mathscr{M}_{\mathbb{Q}}$
of the Hilbert modular variety.

(ii) Assume that $\Sigma$ is invariant under $\phi$. Let $\Sigma^{c}$
be its complement. Then the Goren-Oort stratum $M_{\Sigma}$ (see
definition (\ref{eq:GO stratum}) below) is an integral variety of
$\mathscr{F}_{\Sigma^{c}}$ in the characteristic $p$ fiber $M$
of the Hilbert modular variety.
\end{thm}

\begin{proof}
The proof of part (i) is transcendental. Had there been an integral
variety to $\mathscr{F}_{\Sigma}$ in characteristic 0, it would provide
an algebraic integral variety over $\mathbb{C}$. But over the universal
covering $\mathfrak{H}^{g}$ the \emph{analytic} leaves
of the foliation are easily determined. If $z_{0}=(z_{0,\sigma})_{\sigma\in\mathbb{B}}$
is a point of $\mathfrak{H}^{g}$, the leaf through it
is the coordinate ``plane'' 
\[
H_{\Sigma}(z_{0})=\{z\in\mathfrak{H}^{g}|\,z_{\tau}=z_{0,\tau}\,\,\forall\tau\notin\Sigma\}.
\]
Unless $\Sigma$ is empty or the whole of $\mathbb{B},$ these coordinate
``planes'' do not descend to algebraic varieties in $\Gamma\setminus\mathfrak{H}^{g}$
because the map $H_{\Sigma}(z_{0})\to\Gamma\setminus\mathfrak{H}^{g}$ has a dense image. In fact, \cite{RT} Proposition 3.4 shows that the analytic leaves of these foliations do not even contain any algebraic curves. 


(ii) Let $\mathcal{H}=H_{dR}^{1}(A^\text{\rm univ}/M)$ be the relative de
Rham cohomology of the universal abelian variety over $M$ and $\nabla:\mathcal{H}\to\mathcal{H}\otimes_{\mathcal{O}_{M}}\Omega_{M/\kappa}^{1}$
the Gauss-Manin connection. For a vector field $\xi\in\mathcal{T}_{M/\kappa}(U)$
over a Zariski open set $U\subset M$ we denote by
\[
\nabla_{\xi}:\mathcal{H}(U)\to\mathcal{H}(U)
\]
the $\xi$-derivation obtained by contracting $\Omega_{M/\kappa}^{1}$
with $\xi$. It satisfies ($a\in\mathcal{O}_{M}(U)$)
\[
\nabla_{\xi}(at)=\xi(a)t+a\nabla_{\xi}(t).
\]

If $\sigma,\tau\in\mathbb{B}$ and $\tau$ is not in the $\phi$-orbit
of $\sigma,$ and if $\xi\in\mathscr{F}_{\tau}(U)\subset\mathcal{T}_{M/\kappa}(U)$
then, since $\text{\rm KS}(\xi)$ annihilates 
\[
\mathscr{L}_{\sigma}=\underline{\omega}[\sigma]\subset\underline{\omega}\subset\mathcal{H},
\]
$\nabla_{\xi}$ induces an $\mathcal{O}_{M}$-derivation of the line
bundle $\mathscr{L}_{\sigma}$ over $U$. The same holds with $\mathscr{L}_{\phi\circ\sigma}$
since $\tau\ne\phi\circ\sigma$. By the usual rules of derivations,
we obtain a derivation $\nabla_{\xi}$ of the line bundle $\mathrm{Hom}(\mathscr{L}_{\phi\circ\sigma},\mathscr{L}_{\sigma}^{p}),$
of which the partial Hasse invariant $h_{\phi\circ\sigma}$ is a global
section.

Let us elaborate on the last statement. If $t$ is a section of $\mathscr{L}_{\sigma}$
then the induced derivation of 
\[
\mathscr{L}_{\sigma}^{p}=\mathscr{L}_{\sigma}^{\otimes p}\simeq\mathscr{L_{\sigma}}^{(p)}
\]
(the first is the $p$-th tensor product, the second the base-change
by the absolute Frobenius of $M$) is given by
\[
\nabla_{\xi}^{(p)}(at^{p})=\xi(a)t^{p}
\]
(equivalently, on Frobenius base-change, $\nabla_{\xi}^{(p)}(a\otimes t)=\xi(a)\otimes t$;
note that this is a \emph{canonical} derivation, independent of the
original $\nabla_{\xi}$). If $\mathcal{H}$ is the relative de Rham
cohomology of $A^\text{\rm univ},$ and $\mathcal{H}^{(p)}$ the relative de
Rham cohomology of $A^{{\rm univ}(p)},$ then the same formula applied to
the Gauss-Manin connection of $A^\text{\rm univ}$ gives the Gauss-Manin connection
of $A^{{\rm univ}(p)}.$ Once again, the latter is the \emph{canonical}
connection which exists on the Frobenius base change of any vector
bundle. \emph{Any} section of the form $1\otimes t$ is flat, just
as any function which is a $p$-th power is annihilated by all the
derivations.

Finally, if $h:\mathscr{L}_{\phi\circ\sigma}\to\mathscr{L}_{\sigma}^{p}$
is a homomorphism of line bundles, then$\nabla_{\xi}h$ is the homomorphism
\[
(\nabla_{\xi}h)(t)=\nabla_{\xi}^{(p)}(h(t))-h(\nabla_{\xi}(t)).
\]
\begin{lem} Let $\sigma, \tau \in \mathbb{B}$, not in the same $\phi$-orbit.
The partial Hasse invariant $h_{\phi\circ\sigma}$ is horizontal for
$\xi\in\mathscr{F}_{\tau},$ i.e. $\text{\rm \ensuremath{\nabla_{\xi}}}(h_{\phi\circ\sigma})=0.$
\end{lem}

\begin{proof}
Identify $\mathscr{L}_{\sigma}^{p}$ with $\mathscr{L}_{\sigma}^{(p)}=\Phi_{M}^{*}\mathscr{L}_{\sigma}=\Phi_{M}^{*}(\underline{\omega}[\sigma])=(\Phi_{M}^{*}\underline{\omega})[\phi\circ\sigma]$.
By definition, $h_{\phi\circ\sigma}$ is the $\phi\circ\sigma$ component
of the $\mathcal{O}_{L}$-homomorphism
\[
V:\underline{\omega}\to\underline{\omega}^{(p)}=\Phi_{M}^{*}\underline{\omega},
\]
induced by the Verschiebung isogeny $\text{\rm Ver}_{A^\text{\rm univ}/M}:A^{{\rm univ}(p)}\to A^\text{\rm univ}.$
It is horizontal since the Gauss-Manin connection commutes, in general,
with any map on $\mathcal{H}=H_{dR}^{1}(A^\text{\rm univ}/M)$ induced by an
isogeny, and in particular
\[
\nabla_{\xi}(h_{\phi\circ\sigma})=\nabla_{\xi}^{(p)}\circ h_{\phi\circ\sigma}-h_{\phi\circ\sigma}\circ\nabla_{\xi}=0.
\]
\end{proof}
Let $H_{\phi\circ\sigma}$ be the hypersurface defined by the vanishing
of $h_{\phi\circ\sigma}$ in $M.$ By the results of \cite{G-O},
it is smooth, and $h_{\phi\circ\sigma}$ vanishes on it to first order.
Furthermore, for different $\sigma$'s these hypersurfaces intersect
transversally. Let $x\in H_{\phi\circ\sigma}$ and choose a Zariski
open neighborhood $U$ of $x$ on which $\mathrm{Hom}(\mathscr{L}_{\phi\circ\sigma},\mathscr{L}_{\sigma}^{p})$
is a trivial invertible sheaf. Let~$e$ be a basis of $\mathrm{Hom}(\mathscr{L}_{\phi\circ\sigma},\mathscr{L}_{\sigma}^{p})$
over $U$ and write $h_{\phi\circ\sigma}=he$ for some $h\in\mathcal{O}_{M}(U)$.
Then $H_{\phi\circ\sigma}\cap U$ is given by the equation $h=0$
and $h$ vanishes on it to first order. Furthermore, if $\xi\in\mathscr{F}_{\tau}(U),$
by the Lemma we have
\[
0=\nabla_{\xi}(h_{\phi\circ\sigma})=\xi(h)\cdot e+h\nabla_{\xi}(e),
\]
so along $H_{\phi\circ\sigma}=\{h=0\}$ we also have $\xi(h)=0.$
This proves that $\xi$ is parallel to $H_{\phi\circ\sigma},$ i.e.
$\xi_{x}\in T_{x}H_{\phi\circ\sigma}\subset T_{x}M.$ Let $\Sigma$
be a $\phi$-invariant subset of $\mathbb{B}$. Since the same analysis
holds for every $\sigma\in\Sigma$ and every $\tau\notin\Sigma$ we
get that at every point $x$ of
\begin{equation}
M_{\Sigma}:=\{x|\,h_{\sigma}(x)=0\,\,\forall\sigma\in\Sigma\}\label{eq:GO stratum}
\end{equation}
the $p$-foliation
\[
\mathscr{F}_{\Sigma^{c}}=\oplus_{\tau\notin\Sigma}\mathscr{F}_{\tau}
\]
is contained in $T_{x}M_{\Sigma}\subset T_{x}M.$ As both $\mathscr{F}_{\Sigma^{c}}$
and $TM_{\Sigma}$ are vector bundles of the same rank $g-\#(\Sigma)$,
and both are local direct summands of $TM$, we have shown that the
Goren-Oort stratum $M_{\Sigma}$ is an integral variety of $\mathscr{F}_{\Sigma^{c}}.$
\end{proof}

\section{$V$-foliations on unitary Shimura varieties}

\subsection{Notation and preliminary results on unitary Shimura varieties}

\subsubsection{The moduli scheme}

We now turn to the second type of foliations considered in this paper,
on unitary Shimura varieties in characteristic $p$. Let $K$ be a
CM field, $[K:\mathbb{Q}]=2g$ and $L=K^{+}$ its totally real subfield.
Let $\rho\in Gal(K/L)$ denote complex conjugation. Let $E\subset\mathbb{C}$,
\emph{the field of definition,} be a number field containing all the
conjugates\footnote{We do not insist on the field of definition being the minimal possible
one, i.e. the reflex field of the CM type.} of $K$. For 
\[
\tau\in\mathscr{I}:=\mathrm{Hom}(K,E)=\mathrm{Hom}(K,\mathbb{C})
\]
we write $\bar{\tau}=\tau\circ\rho$. We let $\mathscr{I}^{+} =\mathrm{Hom}(L,E)=\mathscr{I}/\left\langle \rho\right\rangle$ be the set of orbits of $\mathscr{I}$ under the action of $\rho$,
and write its elements as unordered pairs $\{\tau,\bar{\tau}\}.$

Let $d\ge1$ and fix a \emph{PEL-type $\mathcal{O}_{K}$-lattice $(\Lambda,\left\langle ,\right\rangle ,h)$
}of rank $d$ over $\mathcal{O}_{K}$ (\cite{Lan}, 1.2.1.3). Thus
$\Lambda$ is a projective $\mathcal{O}_{K}$-module of rank $d$
(regarded, if we forget the $\mathcal{O}_{K}$-action, as a lattice
of rank $2gd$), $\left\langle ,\right\rangle $ is a non-degenerate
alternating bilinear form $\Lambda\times\Lambda\to2\pi i\mathbb{Z},$
satisfying $\left\langle ax,y\right\rangle =\left\langle x,\bar{a}y\right\rangle $
for $a\in\mathcal{O}_{K}$, and
\[
h:\mathbb{C}\to {\rm End}_{\mathcal{O}_{K}}(\Lambda\otimes\mathbb{R})
\]
is an $\mathbb{R}$-linear ring homomorphism satisfying (i) $\left\langle h(z)x,y\right\rangle =\left\langle x,h(\bar{z})y\right\rangle $
(ii) $(x,y)=(2\pi i)^{-1}\left\langle x,h(i)y\right\rangle $ is an
inner product (symmetric and positive definite) on the real vector
space $\Lambda\otimes\mathbb{R}$. 

The $2gd$-dimensional complex vector space $V=\Lambda\otimes\mathbb{C}$
breaks up as a direct sum
\[
V=V_{0}\oplus V_{0}^{c}
\]
of two $\left\langle ,\right\rangle $-isotropic subspaces, where
$V_{0}=\{v\in V|\,h(z)v=1\otimes z\cdot v\}$ and $V_{0}^{c}=\{v\in V|\,h(z)v=1\otimes\bar{z}\cdot v\}$.
The inclusion $\Lambda\otimes\mathbb{R}\subset\Lambda\otimes\mathbb{C}=V$
allows us to identify $V_{0}$ with the real vector space $\Lambda\otimes\mathbb{R}$,
and then its complex structure is given by $J=h(i).$ 

As representations of $\mathcal{O}_{K}$
\[
V_{0}\simeq\sum_{\tau\in\mathscr{I}}r_{\tau}\tau,\,\,\,V_{0}^{c}\simeq\sum_{\tau\in\mathscr{I}}r_{\tau}\bar{\tau},
\]
where the $r_{\tau}$ are non-negative integers satisfying $r_{\tau}+r_{\bar{\tau}}=d$
for each $\tau$. We call the collection $\{r_{\tau}\}$ (or the formal
sum $\sum_{\tau\in\mathscr{I}}r_{\tau}\tau$) the \emph{signature}
of $(\Lambda,\left\langle ,\right\rangle ,h)$ (\cite{Lan} 1.2.5.2),
or the \emph{CM} \emph{type}.

Let $N\ge3$ (the \emph{tame level}) be an integer which is relatively
prime to the discriminant of the lattice $(\Lambda,\left\langle ,\right\rangle ).$
Let $S$ be the set of \emph{bad} \emph{primes}, defined to be the
rational primes that ramify in $K$, divide $N$, or divide the discriminant
of $\Lambda$. The primes $p\notin S$ are called \emph{good}, and
we fix once and for all such a prime $p$.

Consider the following moduli problem $\mathscr{M}$ over $\mathcal{O}_{E}[1/S]$.
For an $\mathcal{O}_{E}[1/S]$-algebra $R$, the set $\mathscr{M}(R)$
is the set of isomorphism classes of tuples $\underline{A}=(A,\iota,\lambda,\eta)$
where: 
\begin{itemize}
\item $A$ is an abelian scheme of relative dimension $gd$ over $R$.
\item $\iota:\mathcal{O}_{K}\hookrightarrow {\rm End}(A/R)$ is an embedding of
rings, rendering $\text{\rm Lie}(A/R)$ an $\mathcal{O}_{K}$-module of type
$\sum_{\tau\in\mathscr{I}}r_{\tau}\tau.$
\item $\lambda:A\to A^{t}$ is a $\mathbb{Z}_{(p)}^{\times}$-polarization
whose Rosati involution preserves $\iota(\mathcal{O}_{K})$ and induces
on it complex conjugation.
\item $\eta$ is a full level-$N$ structure compatible via $\lambda$ with
$(\Lambda\otimes\widehat{\mathbb{Z}}^{(p)},\left\langle ,\right\rangle ).$
\end{itemize}
See \cite{Lan}, 1.4.1.2 for more details, in particular pertaining
to the level-$N$ structure.

The moduli problem $\mathscr{M}$ is representable by a smooth scheme
over $\mathcal{O}_{E}[1/S]$, which we denote by the same letter.
Its complex points form a finite disjoint union of Shimura varieties
associated with the unitary group of signature $\{r_{\tau}\}.$ Denote
by
\[
\underline{A}^\text{\rm univ}=(A^\text{\rm univ},\iota^\text{\rm univ},\lambda^\text{\rm univ},\eta^\text{\rm univ})
\]
the universal tuple over $\mathscr{M}.$

\bigskip{}

We let $\kappa$ be a finite field, large enough to contain all the
residue fields of the primes of $E$ above $p.$ Fix, once and for
all, an embedding $E\hookrightarrow W(\kappa)[1/p]$, and consider
\[
M=\kappa\times_{\mathcal{O}_{E}[1/S]}\mathscr{M},
\]
the special fiber at the chosen prime of $\mathcal{O}_{E}[1/S]$,
base-changed to $\kappa.$ It is a smooth variety over $\kappa$ of
dimension $\sum_{\{\tau,\bar{\tau}\}\in\mathscr{I}^{+}}r_{\tau}r_{\bar{\tau}}.$
We let $\mathcal{T}$ denote its tangent bundle.

Via the fixed embedding of $\mathcal{O}_{E}$ in $W(\kappa)$ we regard
$\mathscr{I}$ also as the set of homomorphisms of $\mathcal{O}_{K}$
to $\kappa.$ For a prime $\mathfrak{P}$ of $\mathcal{O}_{K}$ above
$p$ we let $\mathscr{I}_{\mathfrak{P}}$ be those homomorphisms that
factor through $\kappa(\mathfrak{P})=\mathcal{O}_{K}/\mathfrak{P},$
\[
\mathscr{I}=\coprod_{\mathfrak{P}|p}\mathscr{I}_{\mathfrak{P}},\,\,\,\mathscr{I}_{\mathfrak{P}}=\mathrm{Hom}(\kappa(\mathfrak{P}),\kappa)=\mathrm{Hom}(\mathcal{O}_{K,\mathfrak{P}},W(\kappa)).
\]
The Frobenius $\phi(x)=x^{p}$ acts on $\mathscr{I}$ on the left
via $\tau\mapsto\phi\circ\tau$ and the $\mathscr{I}_{\mathfrak{P}}$
are its orbits, each of them permuted cyclically by $\phi$. Following
Moonen's convention \cite{Mo}, when we use $\mathscr{I}_{\mathfrak{P}}$
as an indexing set, we shall also write $i$ for $\tau$ and $i+1$
for $\phi\circ\tau$. This will be done without further notice to
avoid the heavy notation $\tau_{i+1}=\phi\circ\tau_{i}.$

\subsubsection{The Kodaira-Spencer isomorphism}

Let $\pi:A^\text{\rm univ}\to\mathscr{M}$ be the structure morphism of the
universal abelian variety, and 
\[
\underline{\omega}=\pi_{*}(\Omega_{A^\text{\rm univ}/\mathscr{M}}^{1})\subset\mathcal{H}=\mathbb{R}^{1}\pi_{*}(\Omega_{A^\text{\rm univ}/\mathscr{M}}^{\bullet})
\]
its relative de-Rham cohomology $\mathcal{H}$ and Hodge bundle $\underline{\omega}$.
These are vector bundles on $\mathscr{M}$ of ranks $2gd$ and $gd$
respectively. The Hodge bundle $\underline{\omega}$ is the dual bundle
to the relative \text{\rm Lie} algebra $\underline{\text{\rm Lie}}=\text{\rm Lie}(A^\text{\rm univ}/\mathscr{M})$.
Since $E$ contains all the conjugates of~$K,$ $S$ contains all
the ramified primes in $K$, and $\mathcal{O}_{\mathscr{M}}$ is an
$\mathcal{O}_{E}[1/S]$-algebra, these vector bundles decompose under
the action of $\mathcal{O}_{K}$ into isotypical parts
\[
\underline{\omega}=\bigoplus_{\tau\in\mathscr{I}}\underline{\omega}[\tau]\subset\bigoplus_{\tau\in\mathscr{I}}\mathcal{H}[\tau]=\mathcal{H}.
\]
For each $\tau$ the rank of $\underline{\omega}[\tau]$ is $r_{\tau}$,
and we have a short exact sequence (the Hodge filtration)
\[
0\to\underline{\omega}[\tau]\to\mathcal{H}[\tau]\to R^{1}\pi_{*}(\mathcal{O}_{A^\text{\rm univ}})[\tau]\to0.
\]
Since $\lambda^\text{\rm univ}$ is a prime-to-$p$ quasi-isogeny,\emph{ }and
 the Rosati involution on $\iota^\text{\rm univ}(\mathcal{O}_{K})$ is
complex conjugation, \emph{after} we base-change from $\mathcal{O}_{E}[1/S]$
to $W(\kappa)$ the polarization induces an isomorphism 
\[
R^{1}\pi_{*}(\mathcal{O}_{A^\text{\rm univ}})[\tau]\simeq R^{1}\pi_{*}(\mathcal{O}_{(A^\text{\rm univ})^{t}})[\bar{\tau}]=\underline{\text{\rm Lie}}[\bar{\tau}]=\underline{\omega}[\bar{\tau}]^{\vee}.
\]
Since $\mathrm{rk}(\underline{\omega}[\tau])+\mathrm{rk}(\underline{\omega}[\bar{\tau}])=d,$
each $\mathcal{H}[\tau]$ is of rank $d$.

We introduce the shorthand notation
\[
\mathcal{P}_{\tau}=\underline{\omega}[\tau].
\]
The Hodge filtration exact sequence can be written therefore, over
$W(\kappa)$, as
\begin{equation}
0\to\mathcal{P}_{\tau}\to\mathcal{H}[\tau]\to\mathcal{P}_{\bar{\tau}}^{\vee}\to0.\label{eq:Hodge Filtration Exact Sequence}
\end{equation}

For any abelian scheme $A/R$, there is a canonical perfect pairing
\[
\{,\}_{dR}:H_{dR}^{1}(A/R)\times H_{dR}^{1}(A^{t}/R)\to R.
\]
In our case, using the prime-to-$p$ quasi-isogeny $\lambda$ to identify
$H_{dR}^{1}(A/R)[\bar{\tau}]$ with $H_{dR}^{1}(A^{t}/R)[\tau]$,
we get a pairing
\[
\{,\}_{dR}:\mathcal{H}[\tau]\times\mathcal{H}[\bar{\tau}]\to\mathcal{O}_{\mathscr{M}}.
\]
Under this pairing $\mathcal{P}_{\tau}$ and $\mathcal{P}_{\bar{\tau}}$
are exact annihilators of each other, and the induced pairing between
$\mathcal{P}_{\tau}$ and $\mathcal{P}_{\tau}^{\vee}$ is the natural
one.

The Gauss-Manin connection is a flat connection
\[
\nabla:\mathcal{H}\to\mathcal{H}\otimes\Omega_{\mathscr{M}}^{1}.
\]
If $\xi\in\mathcal{T}$ is a vector field (on an open set in $\mathscr{M}$,
omitted from the notation), as in \S\ref{subsec integral varieties}, we denote by $\nabla_{\xi}:\mathcal{H}\to\mathcal{H}$
the $\xi$-derivation of $\mathcal{H}$ obtained by contracting $\Omega_{\mathscr{M}}^{1}$
with~$\xi$. Since $\nabla$ commutes with the endomorphisms, $\nabla_{\xi}$
preserves the $\tau$-isotypical parts $\mathcal{H}[\tau]$ for every
$\tau\in\mathscr{I}.$ When $\nabla_{\xi}$ is applied to $\mathcal{P}_{\tau}$
and the result is projected to $\mathcal{P}_{\bar{\tau}}^{\vee}$,
we get an $\mathcal{O}_{\mathscr{M}}$\emph{-linear} homomorphism
\[
\text{\rm KS}^{\vee}(\xi)_{\tau}\in {\rm Hom}(\mathcal{P}_{\tau},\mathcal{P}_{\bar{\tau}}^{\vee})\simeq\mathcal{P}_{\tau}^{\vee}\otimes\mathcal{P}_{\bar{\tau}}^{\vee}.
\]
Using the formalism of the Gauss-Manin connection and the symmetry
of the polarization, it is easy to check that when we identify $\mathcal{P}_{\tau}^{\vee}\otimes\mathcal{P}_{\bar{\tau}}^{\vee}$
with $\mathcal{P}_{\bar{\tau}}^{\vee}\otimes\mathcal{P}_{\tau}^{\vee}$,
\[
\text{\rm KS}^{\vee}(\xi)_{\tau}=\text{\rm KS}^{\vee}(\xi)_{\bar{\tau}}.
\]
Thus $\text{\rm KS}^{\vee}(\xi)_{\tau}$ depends only on the pair $\{\tau,\bar{\tau}\}\in\mathscr{I}^{+}$,
i.e. on $\tau|_{L}.$ When we combine these maps, we get an $\mathcal{O}_{\mathscr{M}}$\emph{-linear}
homomorphism 
\[
\text{\rm KS}^{\vee}(\xi)\in\bigoplus_{\{\tau,\bar{\tau}\}\in\mathscr{I}^{+}}\mathcal{P}_{\tau}^{\vee}\otimes\mathcal{P}_{\bar{\tau}}^{\vee}.
\]

\begin{prop}[The Kodaira-Spencer isomorphism]
 The map 
\[\xymatrix@C=0.6cm{
\text{\rm KS}^{\vee}\colon \mathcal{T}\ar[r]^>>>>>\sim&\underset{{\{\tau,\bar{\tau}\}\in\mathscr{I}^{+}}}{\bigoplus}\mathcal{P}_{\tau}^{\vee}\otimes\mathcal{P}_{\bar{\tau}}^{\vee}
}\]
 sending $\xi$ to $\text{\rm KS}^{\vee}(\xi)$ is an isomorphism.
\end{prop}

We let $\text{\rm KS}$ be the isomorphism dual to $\text{\rm KS}^{\vee},$ namely
\begin{equation}\label{eq:KS-1}\xymatrix@C=0.6cm{
\text{\rm KS}\colon\underset{{\{\tau,\bar{\tau}\}\in\mathscr{I}^{+}}}{\bigoplus}\mathcal{P}_{\tau}\otimes\mathcal{P}_{\bar{\tau}}\ar[r]^>>>>>\sim&\Omega_{\mathscr{M}}^{1}.}
\end{equation}

\subsubsection{\label{subsec:The-mu-ordinary-locus}The $\mu$-ordinary locus of
$M$}

For a general signature, the abelian varieties parametrized by the
mod $p$ fiber $M$ of our moduli space are never ordinary. There
is, however, a dense open set $M^\text{\rm ord}\subset M$ at whose geometric
points $A^\text{\rm univ}$ is ``as ordinary as possible''. To make this
precise we introduce, following \cite{Mo} 1.2.3, certain standard
Dieudonné modules and their associated $p$-divisible groups.

Let $k$ be an algebraically closed field containing $\kappa.$ If
$A$ is an abelian variety over the field $k$, with endomorphisms
by $\mathcal{O}_{K}$, its $p$-divisible group $A[p^{\infty}]$ breaks
up as a product
\[
A[p^{\infty}]=\prod_{\mathfrak{P}|p}A[\mathfrak{P}^{\infty}]
\]
over the primes of $\mathcal{O}_{K}$ above $p$, and $A[\mathfrak{P}^{\infty}]$
becomes a $p$-divisible group with $\mathcal{O}_{K,\mathfrak{P}}$-action.

Fix a prime $\mathfrak{P}|p$ and write $\mathcal{O}=\mathcal{O}_{K,\mathfrak{P}}$.
Let 
\[
\mathfrak{f}:\mathscr{I}_{\mathfrak{P}}=\mathrm{Hom}(\mathcal{O},W(k))\to[0,d]
\]
 be an integer-valued function. Let $M(d,\mathfrak{f})$ be the following
Dieudonné module with $\mathcal{O}$-action over $W(k).$ First,
\[
M(d,\mathfrak{f})=\oplus_{i\in\mathscr{I}_{\mathfrak{P}}}M_{i}
\]
where $M_{i}=\oplus_{j=1}^{d}W(k)e_{i,j}$ is a free $W(k)$-module
of rank $d$. We let $\mathcal{O}$ act on $M_{i}$ via the homomorphism
$i:\mathcal{O}\to W(k).$ We let $F$ (resp. $V)$ be the $\phi$-semilinear
(resp. $\phi^{-1}$-semilinear) endomorphism of $M(d,\mathfrak{f})$
satisfying (recall the convention that if $i$ refers to the embedding
$\tau$ then $i+1$ refers to $\phi\circ\tau$)
\[
F(e_{i,j})=\begin{cases}
\begin{array}{c}
e_{i+1,j}\\
pe_{i+1,j}
\end{array} & \begin{array}{c}
1\le j\le d-\mathfrak{f}(i)\\
d-\mathfrak{f}(i)<j\le d
\end{array}\end{cases}
\]
and
\[
V(e_{i+1,j})=\begin{cases}
\begin{array}{c}
pe_{i,j}\\
e_{i,j}
\end{array} & \begin{array}{c}
1\le j\le d-\mathfrak{f}(i)\\
d-\mathfrak{f}(i)<j\le d
\end{array}\end{cases}.
\]
Then $M(d,\mathfrak{f})$ is a Dieudonné module with $\mathcal{O}$-action,
of rank $[\mathcal{O}:\mathbb{Z}_{p}]d$ over $W(k).$ We let $X(d,\mathfrak{f})$
be the unique $p$-divisible group with $\mathcal{O}$-action over
$k$ whose contravariant Dieudonné module is $M(d,\mathfrak{f}).$

Let $N(d,\mathfrak{f})=M(d,\mathfrak{f})/pM(d,\mathfrak{f}).$ This
is the Dieudonné module of the finite group scheme $Y(d,\mathfrak{f})=X(d,\mathfrak{f})[p].$
The cotangent space of $X(d,\mathfrak{f})$ is canonically isomorphic
to $N(d,\mathfrak{f})[F]=\bigoplus_{i\in\mathscr{I}_{\mathfrak{P}}}\bigoplus_{j=d-\mathfrak{f}(i)+1}^{d}ke_{i,j}.$
It inherits a $\kappa(\mathfrak{P})$-action and its $i$-isotypic
subspace $\bigoplus_{j=d-\mathfrak{f}(i)+1}^{d}ke_{i,j}$ is $\mathfrak{f}(i)$-dimensional.

Let $X$ be a $p$-divisible group with $\mathcal{O}$-action over
$k.$ In \cite{Mo}, Theorem 1.3.7, it is proved that if either $X$
is isogenous to $X(d,\mathfrak{f})$, or $X[p]$ is isomorphic to
$Y(d,\mathfrak{f})$, then $X$ is already isomorphic to $X(d,\mathfrak{f})$.
\begin{defn}
Let $A$ be an abelian variety over $k$ with $\mathcal{O}_{K}$-action,
of dimension $2gd.$ Then $A$ is called $\mu$\emph{-ordinary} if
every $A[\mathfrak{P}^{\infty}]$ with its $\mathcal{O}_{K,\mathfrak{P}}$-action
is isomorphic to some $X(d,\mathfrak{f}).$
\end{defn}

Let $A$ be a $\mu$-ordinary abelian variety over $k$ with $\mathcal{O}_{K}$-action
and CM type $\{r_{\tau}\}$ as before. The cotangent space of $A$
may be identified with that of $A[p^{\infty}]$. From the relation
between the cotangent space of $A[p^{\infty}]$ and its Dieudonné
module, it follows that if $\tau=i\in\mathscr{I}_{\mathfrak{P}}$,
$\mathfrak{f}(i)=r_{\tau}.$

It follows that the function $\mathfrak{f}$ is determined by the
signature, hence all the $\mu$-ordinary $A/k$ parametrized by geometric
points of $M$ have isomorphic $p$-divisible groups. Wedhorn proved
the following fundamental theorem.
\begin{thm*}
\cite{Wed,Mo} There is a dense open set $M^\text{\rm ord}\subset M$ such
that for any geometric point $x\in M(k)$ the abelian variety $A_{x}^\text{\rm univ}$
is $\mu$-ordinary if and only if $x\in M^\text{\rm ord}(k).$
\end{thm*}
Using the slope decomposition explained below, it is possible to attach
a Newton polygon to a $p$-divisible group with $\mathcal{O}_{K}$-action,
and the points of $M^\text{\rm ord}$ are characterized also as those whose
Newton polygon lies \emph{below }every other Newton polygon (Newton
polygons \emph{go up under specialization}).

\subsubsection{\label{subsec:Slope-decomposition}Slope decomposition over $M^\text{\rm ord}$}

We next review the slope decomposition of the $\mathfrak{P}$-divisible
group of a $\mu$-ordinary abelian variety $A$ with $\mathcal{O}_{K}$-action
over an algebraically closed field $k$ containing $\kappa.$ Recall
that $A[\mathfrak{P}^{\infty}]\simeq X(d,\mathfrak{f}).$ For each
$1\le j\le d$ the submodule
\[
M(d,\mathfrak{f})^{j}=\bigoplus_{i\in\mathscr{I}_{\mathfrak{P}}}W(k)e_{i,j}
\]
of $M(d,\mathfrak{f})$ is a sub-Dieudonné module of $\mathcal{O}$-height
1. We define its \emph{slope} as the rational number
\[
\frac{|\{i|\,j>d-\mathfrak{f}(i)\}|}{|\mathscr{I}_{\mathfrak{P}}|}.
\]
(This is the same as the slope in the classification of $p$-divisible groups and Dieudonn\'e modules.)
Note that the slope of $M(d,\mathfrak{f})^{j+1}$ is greater or equal
than the slope of $M(d,\mathfrak{f})^{j}$, and if they are equal,
$M(d,\mathfrak{f})^{j+1}\simeq M(d,\mathfrak{f})^{j}$.

Let $0\le\lambda_{1}<\lambda_{2}<\cdots<\lambda_{r}\le1$ be the distinct
slopes obtained in this way, and for $1\le\nu\le r$ let $d^{\nu}$
be the number of $j$'s with $\mathrm{slope}(M(d,\mathfrak{f})^{j})=\lambda_{\nu}.$
Define a function $\mathfrak{f}^{\nu}:\mathscr{I}_{\mathfrak{P}}\to\{0,d^{\nu}\}$
by
\[
\mathfrak{f}^{\nu}(i)=\begin{cases}
\begin{array}{c}
0\\
d^{\nu}
\end{array} & \begin{array}{c}
\mathrm{if\,\,\,}\sum_{\ell=1}^{\nu-1}d^{\ell}<d-\mathfrak{f}(i)\\
\mathrm{if\,\,\,}\sum_{\ell=1}^{\nu-1}d^{\ell}\ge d-\mathfrak{f}(i).
\end{array}\end{cases}
\]
Then grouping together the $M(d,\mathfrak{f})^{j}$ of slope $\lambda_{\nu}$
we get an isoclinic Dieudonné module isomorphic to $M(d^{\nu},\mathfrak{f}^{\nu}).$
We arrive at the \emph{slope decomposition}
\[
M(d,\mathfrak{f})=\bigoplus_{\nu=1}^{r}M(d^{\nu},\mathfrak{f}^{\nu}),
\]
and similarly for the $p$-divisible group
\[
X(d,\mathfrak{f})=\prod_{\nu=1}^{r}X(d^{\nu},\mathfrak{f}^{\nu}).
\]

This description is valid for $\mu$-ordinary $p$-divisible groups
(with $\mathcal{O}$-action) over algebraically closed fields only.
Its significance stems from the fact that when we study deformations,
the isoclinic $p$-divisible groups with $\mathcal{O}$-action deform
uniquely (are rigid), and the deformations arise only from non-trivial
extensions of one isoclinic subquotient by another one, of a higher
slope. Over an artinian ring with residue field $k$ the slope decomposition
is replaced by a \emph{slope filtration}. The study of the universal
deformation space via these extensions lead Moonen to introduce his
\emph{cascade} structures, which are the main topic of \cite{Mo}.
Finally, we remind the reader that by the Serre-Tate theorem, deformations
of a tuple $\underline{A}\in M(k)$ correspond to deformations of
$A[p^{\infty}]$ with its $\mathcal{O}_{K}$-structure and polarization.
Moonen's theory of cascades supplies therefore ``coordinates'' at
a $\mu$-ordinary point $x\in M^\text{\rm ord}(k),$ reminiscent of the Serre-Tate
coordinates at an ordinary point of the usual modular curve.

\subsubsection{Duality}

Finally, let us examine duality. Quite generally, the Cartier
dual of $A[p^{\infty}]$ is $A^{t}[p^{\infty}].$ A $\mathbb{Z}_{(p)}^{\times}$-polarization
$\lambda$ therefore makes $A[p^{\infty}]$ self-dual. In the presence
of an $\mathcal{O}_{K}$-action as above, the duality induced by $\lambda$
sets $A[\mathfrak{P}^{\infty}]$ in duality with $A[\bar{\mathfrak{P}}^{\infty}].$
Let $\mathfrak{p}$ be the prime of $L$ underlying $\mathfrak{P}.$
We distinguish two cases.

(a) If $\mathfrak{p}\mathcal{O}_{K}=\mathfrak{P\bar{P}}$ is split,
there are no further restrictions on $A[\mathfrak{P}^{\infty}],$
but $A[\bar{\mathfrak{P}}^{\infty}]$ is completely determined by
$A[\mathfrak{P}^{\infty}]$, being its dual group.

(b) If $\mathfrak{p}\mathcal{O}_{K}=\mathfrak{P}$ is inert, $A[\mathfrak{P}^{\infty}]$
is self-dual. In this case let $m=[\kappa(\mathfrak{p}):\mathbb{F}_{p}]$
be the inertia degree of $\mathfrak{p},$ so that $[\kappa(\mathfrak{P}):\mathbb{F}_{p}]=2m$.
Complex conjugation $\rho\in Gal(K/L)$ fixes $\mathfrak{P},$ so
induces an automorphism of $\kappa(\mathfrak{P})$, and for $\tau\in\mathscr{I}_{\mathfrak{P}}$
we have $\tau\circ\rho=\phi^{m}\circ\tau.$ Recall that we denoted
$\tau$ by $i$ and $\phi^{m}\circ\tau$ by $i+m$. If $A$ is $\mu$-ordinary,
the self-duality of $A[\mathfrak{P}^{\infty}]$ is manifested (\cite{Mo}
3.1.1, Moonen's $\varepsilon=+1$ in our case) in a \emph{perfect
symmetric} $W(k)$-linear pairing 
\[
\varphi:M(d,\mathfrak{f})\times M(d,\mathfrak{f})\to W(k)
\]
 such that
\[
\varphi(Fx,y)=\varphi(x,Vy)^{\phi}
\]
\[
\varphi(ax,y)=\varphi(x,\bar{a}y)
\]
($a\in\mathcal{O}_{K})$. This (or the relation $r_{\tau}+r_{\bar{\tau}}=d$)
implies that $\mathfrak{f}(i)+\mathfrak{f}(i+m)=d$. In fact, $M_{i}$
is orthogonal to $M_{i'}$ unless $i'=i+m$ and we can choose the
basis $\{e_{i,j}\}$ in such a way that
\[
\varphi(e_{i,j},e_{i+m,j'})=c_{i,j}\delta_{j',d+1-j}
\]
with some $c_{i,j}\in W(k)^{\times}$ ($\delta_{a,b}$ is Kronecker's
delta). This means that the Dieudonné modules $M(d,\mathfrak{f})^{j}$
and $M(d,\mathfrak{f})^{d+1-j}$ are dual under this pairing. See
\cite{Mo} 3.2.3, case (AU).

\subsection{The $V$-foliations on the $\mu$-ordinary locus}

\subsubsection{Construction}

Consider the universal abelian variety $A^\text{\rm univ}$ over $M,$ which
we now denote for brevity $A,$ and its Verschiebung isogeny
\[
\text{\rm Ver}:A^{(p)}=M\times_{\Phi_{M},M}A\to A.
\]
The relative de Rham cohomology of $A^{(p)},$ denoted $\mathcal{H}^{(p)}$,
may be identified with $\Phi_{M}^{*}\mathcal{H}$, and its Hodge bundle
$\underline{\omega}^{(p)}$ with $\Phi_{M}^{*}\underline{\omega}.$
Letting $a\in\mathcal{O}_{K}$ act on $A^{(p)}$ as $\iota^{(p)}(a)=1\times\iota(a)$
we get an induced action $\iota^{(p)}$ of $\mathcal{O}_{K}$ on $\mathcal{H}^{(p)}$
and on $\underline{\omega}^{(p)}.$ However, for $\tau\in\mathscr{I}$
\[
\mathcal{H}[\tau]^{(p)}:=\Phi_{M}^{*}(\mathcal{H}[\tau])=\mathcal{H}^{(p)}[\phi\circ\tau],
\]
because if $x\in\mathcal{H}[\tau]$ and $1\otimes x\in\mathcal{O}_{M}\otimes_{\phi,\mathcal{O}_{M}}\mathcal{H}[\tau]=\Phi_{M}^{*}(\mathcal{H}[\tau]),$
then 
\[
\iota^{(p)}(a)(1\otimes x)=1\otimes\tau(a)x=\tau(a)^{p}\otimes x=\phi\circ\tau(a)\cdot(1\otimes x).
\]
The isogeny $\text{\rm Ver}$ commutes with the endomorphisms,
\[
\text{\rm Ver}\circ\iota^{(p)}(a)=\iota(a)\circ \text{\rm Ver},
\]
and therefore induces a homomorphism of vector bundles 
\[
V:\mathcal{H}[\tau]\to\text{\rm \ensuremath{\mathcal{H}^{(p)}[\tau]=\mathcal{H}[\phi^{-1}\circ\tau]^{(p)},}}
\]
and similarly on $\underline{\omega}[\tau]=\mathcal{P}_{\tau}$
\[
V:\mathcal{P}_{\tau}\to(\mathcal{P}^{(p)})_{\tau}=(\mathcal{P}_{\phi^{-1}\circ\tau})^{(p)}.
\]
We shall use the notation $\mathcal{P}_{\phi^{-1}\circ\tau}^{(p)}$
for the right hand side.

At a $\mu$-ordinary geometric point $x\in M^\text{\rm ord}(k)$ we may identify
the fiber $\mathcal{H}_{x}=H_{dR}^{1}(A_{x}/k)$ with the (contravariant)
Dieudonné module of $A_{x}[p],$ and the linear map $V:\mathcal{H}_{x}\to\mathcal{H}_{x}^{(p)}$
with the $\phi^{-1}$-semilinear endomorphism $V$ of the Dieudonné
module. Let $\tau=i\in\mathscr{I}_{\mathfrak{P}}.$ Recalling the
shape of the Dieudonné module $M(d,\mathfrak{f})$ of $A_{x}[\mathfrak{P}^{\infty}]$
we conclude that $\text{\rm \ensuremath{\mathcal{P}_{\tau,x}[V]=}}\underline{\omega}_{x}[\tau][V]=0$
if $\mathfrak{f}(i-1)\ge\mathfrak{f}(i)$, and
\[
\text{\rm \ensuremath{\mathcal{P}_{\tau,x}[V]=}}\sum_{j=d-\mathfrak{f}(i)+1}^{d-\mathfrak{f}(i-1)}ke_{i,j}
\]
if $\mathfrak{f}(i-1)<\mathfrak{f}(i)$. We recall that $\mathfrak{f}(i)=r_{\tau}$
and $\mathfrak{f}(i-1)=r_{\phi^{-1}\circ\tau}.$

The following is the main definition of the second part of our paper.
\begin{defn}
Let $\Sigma\subset\mathscr{I}^{+}.$ Let\footnote{The reader might have noticed a twist in our notation. While the foliations
denoted $\mathscr{F}_{\Sigma}$ on Hilbert modular varieties, in the
first part of our paper, grow with $\Sigma$, our current $\mathscr{F}_{\Sigma}$
become smaller when $\Sigma$ grows. This could be solved by labelling
our $\mathscr{F}_{\Sigma}$ by the \emph{complement} of $\Sigma$,
but as the two types of foliations are distinct and of a different
nature, we did not find it necessary to reconcile the two conventions.} $\mathscr{F}_{\Sigma}\subset\mathcal{T}$ be the subsheaf on $M^\text{\rm ord}$
which is the \emph{annihilator}, under the pairing between $\mathcal{T}$
and $\Omega_{M}^{1},$ of
\[
\text{\rm KS}\left(\sum_{\{\tau,\bar{\tau}\}\in\Sigma}(\mathcal{P}_{\tau}\otimes\mathcal{P}_{\bar{\tau}})[V\otimes V]\right).
\]
\end{defn}

Our first goal is to prove that $\mathscr{F}_{\Sigma}$ is a $p$-foliation.
Note that at every $x\in M(k)$
\begin{equation}
\mathcal{P}_{\tau,x}\otimes\mathcal{P}_{\bar{\tau},x}[V\otimes V]=\mathcal{P}_{\tau,x}[V]\otimes\mathcal{P}_{\bar{\tau},x}+\mathcal{P}_{\tau,x}\otimes\mathcal{P}_{\bar{\tau},x}[V].\label{eq:ker_V_times_V}
\end{equation}
By the discussion above, if $x\in M^\text{\rm ord}(k),$ the first term is
a subspace whose dimension is 
\[
\max\{0,r_{\tau}-r_{\phi^{-1}\circ\tau}\}\cdot(d-r_{\tau})
\]
and the second is of dimension $r_{\tau}\cdot\max\{0,r_{\phi^{-1}\circ\tau}-r_{\tau}\}.$
Here we used the relations $\overline{\phi^{-1}\circ\tau}=\phi^{-1}\circ\tau\circ\rho=\phi^{-1}\circ\bar{\tau}$
and $r_{\bar{\tau}}=d-r_{\tau}.$ At most one of the terms is non-zero,
and both are zero if and only if either $r_{\tau}=r_{\phi^{-1}\circ\tau},$$r_{\tau}=0$
or $r_{\tau}=d.$ In particular,
\[
\dim_{k}\mathcal{P}_{\tau,x}\otimes\mathcal{P}_{\bar{\tau},x}[V\otimes V]
\]
is the same for all $x\in M^\text{\rm ord}(k).$
\begin{lem}
\label{lem:Foliation is smooth}For $\{\tau,\bar{\tau}\}\in\mathscr{I}^{+}$
let
\[
r_{V}^\text{\rm ord}\{\tau,\bar{\tau}\}=\max\{0,r_{\tau}-r_{\phi^{-1}\circ\tau}\}\cdot(d-r_{\tau})+r_{\tau}\cdot\max\{0,r_{\phi^{-1}\circ\tau}-r_{\tau}\}
\]
(this quantity is symmetric in $\tau$ and $\bar{\tau}).$ Over $M^\text{\rm ord}$,
the subsheaf $\mathscr{F}_{\Sigma}$ is a vector sub-bundle of $\mathcal{T}$
of corank $r_{V}(\Sigma)=\sum_{\{\tau,\bar{\tau}\}\in\Sigma}r_{V}^\text{\rm ord}\{\tau,\bar{\tau}\}.$
Its rank is given by the formula

\[
\mathrm{rk}(\mathscr{F}_{\Sigma})=\sum_{\{\tau,\bar{\tau}\}\notin\Sigma}r_{\tau}r_{\bar{\tau}}+\sum_{\{\tau,\bar{\tau}\}\in\Sigma}\min\{r_{\tau},r_{\phi^{-1}\circ\tau}\}\cdot\min\{r_{\bar{\tau}},r_{\phi^{-1}\circ\bar{\tau}}\}.
\]
\end{lem}

\begin{proof}
The sheaf $\mathscr{G}_{\tau}=(\mathcal{P}_{\tau}\otimes\mathcal{P}_{\bar{\tau}})[V\otimes V]$
is a subsheaf of $\mathcal{P}_{\tau}\otimes\mathcal{P}_{\bar{\tau}}$,
and $\mathscr{F}_{\Sigma},$ the annihilator of $\text{\rm KS}(\sum_{\{\tau,\bar{\tau}\}\in\Sigma}\mathscr{G}_{\tau})$
in $\mathcal{T}$, is \emph{saturated}. This is because $\mathcal{O}_{M}$
has no zero divisors: if $f\in\mathcal{O}_{M}$, $\xi\in\mathcal{T}$
and $f\xi\in\mathscr{F}_{\Sigma}$, then for every $\omega\in \text{\rm KS}(\sum_{\{\tau,\bar{\tau}\}\in\Sigma}\mathscr{G}_{\tau})$
we have $f\left\langle \omega,\xi\right\rangle =\left\langle \omega,f\xi\right\rangle =0$,
so $\left\langle \omega,\xi\right\rangle =0.$

Quite generally, if $M$ is a variety over a field $k$, $\mathcal{F}$
and $\mathcal{G}$ are locally free sheaves of rank $n$, and $T\in\mathrm{Hom}_{\mathcal{O}_{M}}(\mathcal{F},\mathcal{G})$
is such that the fiber rank $\mathrm{rk}(T_{x})=m$ is constant on
$M$, then $\ker(T)$ is locally a direct summand (i.e. a vector sub-bundle)
of rank $n-m$. Compare \cite{Mu}, II.5 Lemma 1, p.51. It is in fact
enough to verify the constancy of $\mathrm{rk}(T_{x})$ at closed
points, because every other point $y\in M$ contains closed points
in its closure, and the fiber rank can only go up under specialization.
Applying this with $\mathcal{F}=\mathcal{P}_{\tau}\otimes\mathcal{P}_{\bar{\tau}}$,
$\mathcal{G}=\mathcal{F}^{(p)}$ and $T=V\otimes V$ we deduce that
$\mathscr{G}_{\tau}$, hence also $\mathscr{F}_{\Sigma}$, are vector
sub-bundles.

The rank of $\mathscr{F}_{\Sigma}$ is
\[
\sum_{\{\tau,\bar{\tau}\}\in\mathscr{I}^{+}}r_{\tau}r_{\bar{\tau}}-\sum_{\{\tau,\bar{\tau}\}\in\Sigma}r_{V}^\text{\rm ord}\{\tau,\bar{\tau}\}.
\]
For $\{\tau,\bar{\tau}\}\in\Sigma$ we first assume that $r_{\phi^{-1}\circ\tau}\le r_{\tau}.$
We then have
\[
r_{\tau}r_{\bar{\tau}}-r_{V}^\text{\rm ord}\{\tau,\bar{\tau}\}=r_{\tau}r_{\bar{\tau}}-(r_{\tau}-r_{\phi^{-1}\circ\tau})r_{\bar{\tau}}=r_{\phi^{-1}\circ\tau}r_{\bar{\tau}}.
\]
But by our assumption $r_{\phi^{-1}\circ\tau}=\min\{r_{\tau},r_{\phi^{-1}\circ\tau}\}$
and $r_{\bar{\tau}}=\min\{r_{\bar{\tau}},r_{\phi^{-1}\circ\bar{\tau}}\}.$
The case $r_{\phi^{-1}\circ\tau}\ge r_{\tau}$ is treated similarly.
\end{proof}

\subsubsection{Closure under Lie brackets and $p$-power}
\begin{lem}
The vector bundle $\mathscr{F}_{\Sigma}$ is involutive: if $\xi,\eta$
are sections of $\mathscr{F}_{\Sigma},$ so is $[\xi,\eta].$
\end{lem}

\begin{proof}
The proof is essentially the same as the proof of Proposition 3 in
\cite{G-dS1}. For $\alpha\in\mathcal{P}_{\tau}$ and $\beta\in\mathcal{P}_{\bar{\tau}}$
we have the formula
\[
\left\langle \text{\rm KS}(\alpha\otimes\beta),\xi\right\rangle =\{\nabla_{\xi}(\alpha),\beta\}_{dR}\in\mathcal{O}_{M}
\]
(loc. cit. Lemma 4). Thus, $\xi\in\mathscr{F}_{\Sigma}$ if and only
if for every $\tau$ such that $\{\tau,\bar{\tau}\}\in\Sigma$ we
have
\[
\nabla_{\xi}(\mathcal{P}_{\tau}[V])\perp\mathcal{P}_{\bar{\tau}}
\]
under the pairing $\{,\}_{dR}:\mathcal{H}[\tau]\times\mathcal{H}[\bar{\tau}]\to\mathcal{O}_{M}.$
But the left annihilator of $\mathcal{P}_{\bar{\tau}}$ is $\mathcal{P}_{\tau}.$
The Gauss-Manin connection commutes with isogenies, so in particular
carries $\mathcal{H}[V]$ to itself. It follows that $\xi\in\mathscr{F}_{\Sigma}$
if and only if
\begin{equation}
\nabla_{\xi}(\mathcal{P}_{\tau}[V])\subset\mathcal{P}_{\tau}[V]\label{eq:criterion}
\end{equation}
for any $\tau$ such that $\{\tau,\bar{\tau}\}\in\Sigma$. The Gauss-Manin
connection is well-known to be flat, i.e.
\[
\nabla_{[\xi,\eta]}=\nabla_{\xi}\circ\nabla_{\eta}-\nabla_{\eta}\circ\nabla_{\xi},
\]
so if both $\xi$ and $\eta$ lie in $\mathscr{F}_{\Sigma},$ $(\ref{eq:criterion})$
implies that so does $[\xi,\eta].$
\end{proof}
\begin{lem}
The vector bundle $\mathscr{F}_{\Sigma}$ is $p$-closed: if $\xi$
is a section of $\mathscr{F}_{\Sigma}$, so is $\xi^{p}$.
\end{lem}

\begin{proof}
Again we follow the proof of Proposition 3 in \cite{G-dS1}. The $p$-curvature
\[
\psi(\xi)=\nabla_{\xi}^{p}-\nabla_{\xi^{p}}
\]
does not vanish identically, but is only a nilpotent endomorphism
of $\mathcal{H}$ (\cite{Ka-Tur}, §5). However, denoting by $\mathcal{H}^{(p)}=\Phi_{M}^{*}\mathcal{H}=H_{dR}^{1}(A^{(p)}/M)$
$(A=A^\text{\rm univ})$ the relative de Rham cohomology of $A^{(p)},$ and
by $F$ the map induced by the relative Frobenius $\text{\rm Fr}:A\to A^{(p)}$
on cohomology, we have
\[
\mathcal{H}[V]=F^{*}\mathcal{H}^{(p)}.
\]
Furthermore, since the Gauss-Manin connection commutes with any isogeny,
the following diagram commutes
\begin{equation}
\label{diagram1}
 \xymatrix@=1cm{ \mathcal{H}^{(p)} \ar[r]^{F^\ast} \ar[d]_{\nabla^{(p)} = \nabla_{\text{\rm can}}} &\mathcal{H}\ar[d]^{\nabla} \\
\mathcal{H}^{(p)} \otimes \Omega_M \ar[r]^{F^{*}\otimes1} & \mathcal{H} \otimes \Omega_M.
}
\end{equation}
Here $\nabla^{(p)},$ the Gauss-Manin connection of $A^{(p)}$, turns
out to be the canonical connection $\nabla_{can}$ that exists on
\emph{any} vector bundle of the form $\Phi_{M}^{*}\mathcal{H},$ namely
if $f\otimes\alpha\in\mathcal{O}_{M}\otimes_{\phi,\mathcal{O}_{M}}\mathcal{H},$
\[
\nabla_{can}(f\otimes\alpha)=df\otimes\alpha
\]
(that this is well-defined follows from the rule $d(g^{p}f)=g^{p}df$).
The commutativity of $(\ref{diagram1})$ implies that the restriction
of $\nabla$, the Gauss-Manin connection of $A$, to $F^{*}\mathcal{H}^{(p)}=\mathcal{H}[V]$
is the connection denoted $\nabla_{can}$ in \cite{Ka-Tur}. It follows
from Cartier's theorem (loc. cit. Theorem 5.1) that $\psi(\xi)$ vanishes
when restricted to $\mathcal{H}[V]$. We conclude the proof as in
the previous lemma, using the criterion $(\ref{eq:criterion})$ for
$\xi$ to lie in $\mathscr{F}_{\Sigma}$.
\end{proof}
Altogether we proved the following.
\begin{thm}
The sheaf $\mathscr{F}_{\Sigma}$ is a smooth $p$-foliation on $M^\text{\rm ord}$.
\end{thm}

We stress the difference between the tautological foliations on Hilbert
modular varieties, which were $p$-closed only if $\Sigma$ was invariant
under the action of $\phi,$ and the $V$-foliations on unitary Shimura
varieties that are always $p$-closed. The reason lies in the last
Lemma, in the delicate relation between the $p$-curvature of the
Gauss-Manin connection and the kernel of Verschiebung.

\subsection{Relation with Moonen's cascade structure}

\subsubsection{Moonen's cascade structure}

In \cite{Mo} Moonen generalized the notion of Serre-Tate coordinates
to any Shimura variety of PEL type in characteristic $p$. Here we
recall his results in the case of our unitary Shimura variety $M$.

Fix a $\mu$-ordinary geometric point $x:\text{\rm Spec}(k)\to M^\text{\rm ord}$. Let
$W=W(k)$ be the ring of Witt vectors over $k$, and consider the
category $\mathbf{C}_{W}$ of local artinian $W$-algebras with residue
field $k$, morphisms being local homomorphisms inducing the identity
on $k$. Let $\mathbf{FS}_{W}$ be the category of affine formal schemes
$\mathfrak{X}$ over $W$ with the property that $\Gamma(\mathfrak{X},\mathcal{O}_{\mathfrak{X}})$
is a profinite $W$-algebra (the last condition regarded as a ``smallness''
condition). By a theorem of Grothendieck, associating to $\mathfrak{X}\in\mathbf{FS}_{W}$
the functor of points $R\mapsto\mathfrak{X}(R)$ ($R\in\mathbf{C}_{W})$
identifies $\mathbf{FS}_{W}$ with the category of left-exact functors
from $\mathbf{C}_{W}$ to sets.

Equip $\mathbf{FS}_{W}$ with the flat topology and let $\mathfrak{T}=\mathbf{\widehat{FS}}_{W}$
be the topos of sheaves of sets on it. Since the flat topology is
coarser than the canonical topology, $\mathbf{FS}_{W}$ embeds in
$\mathbf{\widehat{FS}}_{W}$ (Yoneda's lemma). In particular, we may
consider ${\rm Spf}(\widehat{\mathcal{O}}_{\mathscr{M},x})$ as a sheaf
of sets on $\mathbf{FS}_{W}$.

Let $\mathbb{D}=\mathrm{Def}(\underline{X}')$ be the universal deformation
space of the pair
\[
\underline{X}'=(A_{x}^\text{\rm univ}[p^{\infty}],\iota_{x}^\text{\rm univ}).
\]
For every $\mathfrak{X}\in\mathbf{FS}_{W}$ the deformations of $\underline{X}'$
over $\mathfrak{X}$ make up the set $\mathbb{D}(\mathfrak{X})$,
and $\mathbb{D}\in\mathbf{\widehat{FS}}_{W}.$ In fact, it is representable
by a formal scheme in $\mathbf{FS}_{W}$.

If $\underline{X}'^{,D}$ is the Cartier dual of $X$, with the $\mathcal{O}_{K}$-action
\[
\iota^{D}(a):=\iota(\bar{a})^{D},
\]
then $\mathrm{Def}(\underline{X}'^{,D})=\mathbb{D}$ as well, since
any deformation of $\underline{X}'$ yields a deformation of $\underline{X}'^{,D}$
and vice versa. (The cascade structures defined below will be dual,
though.) The polarization $\lambda_{x}^\text{\rm univ}:X\simeq X^{D}$ (intertwining
the actions $\iota$ and $\iota^{D}$) induces an automorphism
\[
\gamma:\mathbb{D}=\mathrm{Def}(\underline{X}')\simeq\mathrm{Def}(\underline{X}'^{,D})=\mathbb{D},
\]
and the subsheaf $\mathbb{D}^{\lambda}=\{x\in\mathbb{D}|\,\gamma(x)=x\}$
of symmetric elements in $\mathbb{D}$ is the universal deformation
space of 
\[
\underline{X}=(A_{x}^\text{\rm univ}[p^{\infty}],\iota_{x}^\text{\rm univ},\lambda_{x}^\text{\rm univ}).
\]
See \cite{Mo}, 3.3.1. Like $\mathbb{D},$ its sub-functor $\mathbb{D}^{\lambda}\in\mathbf{\widehat{FS}}_{W},$
and is representable. Indeed, $\mathbb{D}^{\lambda}$ is represented
by the formal scheme ${\rm Spf}(\widehat{\mathcal{O}}_{\mathscr{M},x})$,
and its characteristic $p$ fiber by ${\rm Spf}(\widehat{\mathcal{O}}_{M,x})$.

\bigskip{}

We refer to \cite{Mo}, 2.2.1, for the precise definition of an $r$\emph{-cascade}
in a topos $\mathfrak{T}$. A $1$-cascade is a point, a $2$-cascade
is a sheaf of commutative groups, and a $3$-cascade is a biextension.
The general structure of an $r$-cascade in $\mathfrak{T}$ generalizes
these notions (see below).

We have a decomposition $\underline{X} = \bigtimes_{\mathfrak{p}} \underline{X}[\mathfrak{p}^\infty]$, where $\mathfrak{p}$ runs over primes of $L$ dividing~$p$, and a corresponding decomposition $\mathbb{D} = \bigtimes_\mathfrak{p} \mathbb{D}_\mathfrak{p}$. The products are understood as fibered products over ${\rm Spec}(W)$. In \cite{Mo}, 2.3.6, Moonen defines a structure of an $r$-cascade
on each $\mathbb{D}_\mathfrak{p}$, where $r = r(\mathfrak{p})$ is the number of slopes of $X[\mathfrak{p}^\infty]$ (see
§\ref{subsec:Slope-decomposition}).  One might refer to $\mathbb{D}$ as a multi-cascade. Thus ${\rm Spf}(\widehat{\mathcal{O}}_{\mathscr{M},x})$
(or rather, the sheaf that it represents) is endowed with the structure
of symmetric elements in the self-dual multi-cascade $\mathbb{D}.$

\subsubsection{Previous results}

When $K$ is quadratic imaginary and $p$ is inert we showed in \cite{G-dS1},
Theorem 13, that the foliation we have constructed (denoted there
$\mathcal{T}S^{+}$, and here $\mathscr{F}_{\Sigma}$) is compatible
with Moonen's cascade structure on $M^\text{\rm ord}$ in the following sense.
Let $(n,m)$ be the signature, $n> m$, so that $\mathrm{rk}(\mathscr{F}_{\Sigma})=m^{2}.$
Let $x\in M^\text{\rm ord}(k).$ The number of slopes of $A_{x}^\text{\rm univ}[p^{\infty}]$
turns out to be three, and ${\rm Spf}(\widehat{\mathcal{O}}_{M,x})$ acquires
from the cascade structure a structure of a $\widehat{\mathbb{G}}_{m}^{m^{2}}$-torsor.
The formal torus $\widehat{\mathbb{G}}_{m}^{m^{2}}$ gives rise to
an $m^{2}$-dimensional subspace of the $mn$-dimensional tangent space at~$x$, and
in loc. cit. we proved that this subspace coincided with the foliation. 

We find that while the cascade structure ``lives'' only on the formal
neighborhood of $x$, and does not globalize, its ``trace'' on the
tangent space globalizes to the foliation that we constructed in the
tangent bundle of $M^\text{\rm ord}$.

\subsubsection{The general case: subspaces of the tangent bundle defined by the
cascade}

We shall now describe how this result generalizes to the setting of
our paper. For simplicity let us assume that~$p$ is inert in $L,$
$p\mathcal{O}_{L}=\mathfrak{p}$, $\mathbb{D} = \mathbb{D}_\mathfrak{p}$. Both the foliation and the cascade
structure break up as products of corresponding structures over the
primes~$\mathfrak{p}$ of $L$ dividing $p,$ so with a little more
book-keeping the general case follows the same pattern as the inert
case. 

We shall also assume that $\mathfrak{p}$ splits in $K,$ and write
as before $\mathfrak{p}\mathcal{O}_{K}=\mathfrak{P}\bar{\mathfrak{P}}.$
The other case ($\mathfrak{p}$ inert in $K$) is a little more complicated,
but can be handled in a similar way. The assumption that $\mathfrak{p}$
splits in $K$ allows us to concentrate on the deformation space of
$\underline{X}=A_{x}^\text{\rm univ}[\mathfrak{P}^{\infty}],$ as a $p$-divisible
group with $\mathcal{O}_{K}$-action, and ignore the polarization.
The resulting deformation space is just an $r$-cascade $\mathbb{D}_{\mathfrak{P}}.$ Then, $\mathbb{D}=\mathbb{D}_{\mathfrak{P}}\times\mathbb{D}_{\mathfrak{\mathfrak{P}}}^{\vee},$ and
the polarization induces an isomorphism $\mathbb{D}_{\mathfrak{P}}\simeq\mathbb{D}_{\mathfrak{P}}^{\vee}=\mathbb{D}_{\bar{\mathfrak{P}}}$
that we denote $x\mapsto\gamma(x)$; $\mathbb{D}^{\lambda}$ is the
set of pairs $(x,\gamma(x))\in\mathbb{D}$, and is therefore isomorphic
to $\mathbb{D}_{\mathfrak{P}}.$

Because of the relation $\mathscr{F}_{\Sigma_{1}\cup\Sigma_{2}}=\mathscr{F}_{\Sigma_{1}}\cap\mathscr{F}_{\Sigma_{2}}$
it is enough to determine the relation of $\mathscr{F}_{\{\tau,\bar{\tau}\}}$
to the cascade structure, and the general case will follow from it.
As we have seen in $(\ref{eq:KS-1}),$ the Kodaira-Spencer map yields
an isomorphism
\[\xymatrix@C=0.6cm{
\text{\rm KS}^{\vee}\colon \mathcal{T}\ar[r]^>>>>>\sim&\underset{\{\sigma,\bar{\sigma}\}\in\mathscr{I}^{+}}{\bigoplus}\mathcal{\mathcal{P}_{\sigma}^{\vee}\otimes}\mathcal{P}_{\bar{\sigma}}^{\vee},}
\]
and we may write
\[
\mathcal{T}_{\{\sigma,\bar{\sigma}\}}=(\text{\rm KS}^{\vee})^{-1}(\mathcal{\mathcal{P}_{\sigma}^{\vee}\otimes}\mathcal{P}_{\bar{\sigma}}^{\vee}).
\]
We then have
\begin{equation}
\mathscr{F}_{\{\tau,\bar{\tau}\}}=\bigoplus_{\{\sigma,\bar{\sigma}\}\ne\{\tau,\bar{\tau}\}}\mathcal{T}_{\{\sigma,\bar{\sigma}\}}\oplus\mathscr{E}_{\{\tau,\bar{\tau}\}}\label{eq:F_=00007B=00005Ctau,=00005Ctaubar=00007D}
\end{equation}
where $\mathscr{E}_{\{\tau,\bar{\tau}\}}\subset\mathcal{T}_{\{\tau,\bar{\tau}\}}$
is the annihilator of $\text{\rm KS}(\mathcal{P}_{\tau}\otimes\mathcal{P}_{\bar{\tau}}[V\otimes V])$.
If we choose the labelling so that $r_{\phi^{-1}\circ\tau}\le r_{\tau}$
then over $M^\text{\rm ord}$ we have $\mathcal{P}_{\tau}\otimes\mathcal{P}_{\bar{\tau}}[V\otimes V]=\mathcal{P}_{\tau}[V]\otimes\mathcal{P}_{\bar{\tau}}$
and
\[
\mathrm{rk}\;\mathscr{E}_{\{\tau,\bar{\tau}\}}=r_{\phi^{-1}\tau}\cdot(d-r_{\tau}).
\]
Although our $p$-foliation is $\mathscr{F}_{\{\tau,\bar{\tau}\}}$
and not $\mathscr{E}_{\{\tau,\bar{\tau}\}}$ (the latter is an involutive
sub-bundle but need not be $p$-closed!), for our purpose it will
be enough to relate $\mathscr{E}_{\{\tau,\bar{\tau}\}}$ to the cascade
structure.

Choose the notation $\mathfrak{p}\mathcal{O}_{K}=\mathfrak{P}\bar{\mathfrak{P}}$
so that $\tau\in\mathscr{I}_{\mathfrak{P}}$ and $\bar{\tau}\in\mathscr{I}_{\bar{\mathfrak{P}}}$.
As in §\ref{subsec:The-mu-ordinary-locus} let $X(d,\mathfrak{f})$
(now for $\mathfrak{f}:\mathscr{I}_{\mathfrak{P}}\to[0,d]$) be the
standard $p$-divisible group with $\mathcal{O}_{K}$-structure over
$k$ whose Dieudonné module is $M(d,\mathfrak{f})$. As before, $\mathfrak{f}(i)=r_{\sigma}$
if $i=\sigma\in\mathscr{I}_{\mathfrak{P}}.$ Let $r$ be the number
of distinct slopes of $X(d,\mathfrak{f})$ (not to be confused with
the CM type $\{r_{\sigma}\}$).

Fix $x\in M^\text{\rm ord}(k)$. We recall some notation related to the cascade
structure at $x.$ There is a canonical $r$-cascade structure 
\[
\mathscr{C}=\{\Gamma^{(a,b)},G^{(a,b)}|\,1\le a<b\le r\}
\]
on the deformation space
\[
{\rm Spf}(\widehat{\mathcal{O}}_{M,x})=\mathbb{D}_{\mathfrak{P}}.
\]
The $\Gamma^{(a,b)}$ are formal schemes supported at $x,$ and $\Gamma^{(1,r)}=\mathbb{D}_{\mathfrak{P}}.$
They are equipped with (left and right) morphisms
\[
\lambda:\Gamma^{(a,b)}\to\Gamma^{(a,b-1)},\,\,\,\rho:\Gamma^{(a,b)}\to\Gamma^{(a+1,b)},
\]
satisfying $\lambda\circ\rho=\rho\circ\lambda$ (where applicable).
Each $\Gamma^{(a,b)}$ is endowed with a structure of a relative bi-extension
of
\[
\Gamma^{(a,b-1)}\times_{\Gamma^{(a+1,b-1)}}\Gamma^{(a+1,b)}
\]
by a formal group that we denote by $G^{(a,b)}$ (in the category
of formal schemes \emph{over} $\Gamma^{(a+1,b-1)}$). See the following
diagram.

\[ \xymatrix@C=1.7cm{& \Gamma^{(a, b)}\ar[ddl]_\lambda\ar[ddr]^\rho\ar@{..>}[d]_{G^{(a, b)}} & \\ 
&  { \Gamma^{(a, b-1)}\!\!\!\!\! \!\!\!\underset{{\Gamma^{(a+1, b-1)}}}{\times}\!\!\!\! \!\!\!\!\Gamma^{(a+1, b)}}\ar@{..>}[dr]\ar@{..>}[dl]& \\
\Gamma^{(a, b-1)}\ar[dr]^\rho & &\Gamma^{(a+1, b)}\ar[dl]_\lambda\\
& \Gamma^{(a+1, b-1)}& 
 }\]

In fact, $G^{(a,b)}={\rm Ext}(X^{(a)},X^{(b)})$ where $X^{(\nu)}$ is the
$\nu$-th isoclinic component of $A_{x}^\text{\rm univ}[\mathfrak{P}^{\infty}]\simeq X(d,\mathfrak{f})$.
This identification should be interpreted as an identity between fppf
sheaves of $\mathcal{O}_{K}$-modules; each $X^{(\nu)},$ with its
$\mathcal{O}_{K}$-structure, is rigid, so admits a unique canonical
lifting $X_{R}^{(\nu)}$ to any local artinian ring $R$ with residue
field $k$, and
\[
G^{(a,b)}(R)={\rm Ext}_{R}(X_{R}^{(a)},X_{R}^{(b)}).
\]

Let
\[
\mathcal{T}_{x}^{(a,b)}=\ker d(\lambda,\rho)
\]
where $d(\lambda,\rho)$ is the differential of the map
\[
(\lambda,\rho):\Gamma^{(a,b)}\to\Gamma^{(a,b-1)}\times_{\Gamma^{(a+1,b-1)}}\Gamma^{(a+1,b)}.
\]
This is the subspace of the tangent space of $\Gamma^{(a,b)}$ ``in
the direction'' of $G^{(a,b)}.$ Thus,
\[
\mathcal{T}_{x}^{(a,b)}={\rm Ext}_{k[\varepsilon]}(X_{k[\varepsilon]}^{(a)},X_{k[\varepsilon]}^{(b)})
\]
is the $k$-vector space of all the extensions, over the ring of dual
numbers $k[\varepsilon]$, of the formal $\mathcal{O}_{K}$-module
$X^{(a)}$ by the formal $\mathcal{O}_{K}$-module $X^{(b)}.$ 

Using these spaces we define subspaces
\[
\mathcal{T}_{x}^{[a,b]}\subset\mathcal{T}_{x}
\]
by a descending induction on $b-a$ for $1\le a<b\le r$. First, $\mathcal{T}_{x}^{[1,r]}=\mathcal{T}_{x}^{(1,r)}.$
Suppose $\mathcal{T}_{x}^{[a,b]}$ have been defined when $b-a>s$,
let $b-a=s$ and consider $U=\mathcal{T}_{x}/\sum_{[a,b]\subsetneq[a',b']}\mathcal{T}_{x}^{[a',b']}.$
Then $\mathcal{T}_{x}^{(a,b)}$ is a subspace of $U$ and we let $\mathcal{T}_{x}^{[a,b]}$
be its preimage in $\mathcal{T}_{x}$, so that
\[
\mathcal{T}_{x}^{(a,b)}=\mathcal{T}_{x}^{[a,b]}/\sum_{[a,b]\subsetneq[a',b']}\mathcal{T}_{x}^{[a',b']}.
\]
For example, if $r=3$ then ${\rm Spf}(\widehat{\mathcal{O}}_{M,x})=\Gamma^{(1,3)}$
has the structure of a bi-extension of $\Gamma^{(1,2)}\times\Gamma^{(2,3)}$
by $G^{(1,3)},$ $\mathcal{T}_{x}^{[1,3]}$ is the tangent space to
$G^{(1,3)},$ $\mathcal{T}_{x}^{[1,2]}/\mathcal{T}_{x}^{[1,3]}=\mathcal{T}_{x}^{(1,2)}$
is the tangent space to $\Gamma^{(1,2)}$ and, likewise, $\mathcal{T}_{x}^{[2,3]}/\mathcal{T}_{x}^{[1,3]}=\mathcal{T}_{x}^{(2,3)}$
is the tangent space to $\Gamma^{(2,3)}$ (all tangent spaces are
at the origin). In general, we have defined a collection of subspaces
of $\mathcal{T}_{x}$ indexed by closed intervals $[a,b]$ with $1\le a<b\le r,$ so
that $\mathcal{T}_{x}^{I}\supset\mathcal{T}_{x}^{J}$ whenever $I\subset J.$
The filtration of $\mathcal{T}_{x}$ by the $\mathcal{T}_{x}^{[a,b]}$
is not linearly ordered, but its graded pieces are the $\mathcal{T}_{x}^{(a,b)},$
the tangent spaces to the $G^{(a,b)}.$

\subsubsection{The relation between the cascade and $\mathscr{E}_{\{\tau,\bar{\tau}\}}$}

Depending on our $\tau$, we define two integers $0\le p_{\tau}\le q_{\tau}\le r.$
Recall that the $\nu$-th isoclinic group $X^{(\nu)}$ is of the form
\[
X^{(\nu)}=X(d^{\nu},\mathfrak{f}^{(\nu)})=X(1,\mathfrak{g}^{(\nu)})^{d^{\nu}}
\]
with $\mathfrak{g}^{(\nu)}:\mathscr{I}_{\mathfrak{P}}\to\{0,1\},$
$\mathfrak{f}^{(\nu)}=d^{\nu}\mathfrak{g}^{(\nu)},$ and $\mathfrak{g}^{(\nu+1)}\ge\mathfrak{g}^{(\nu)}.$
If $\tau=i\in\mathscr{I}_{\mathfrak{P}}$ we let
\[
p_{\tau}=\#\{\nu|\,\mathfrak{g}^{(\nu)}(i)=0\},\,\,\,q_{\tau}=p_{\phi^{-1}\circ\tau}=\#\{\nu|\,\mathfrak{g}^{(\nu)}(i-1)=0\},
\]
so that $\mathfrak{f}(i)=\sum_{\nu=p_{\tau}+1}^{r}d^{\nu}$ and similarly
$\mathfrak{f}(i-1)=\sum_{\nu=q_{\tau}+1}^{r}d^{\nu}$.
\begin{prop}
The fiber of $\mathscr{E}_{\{\tau,\bar{\tau}\}}$ at $x$ coincides
with the subspace $\mathcal{T}_{x}^{[p_{\tau},q_{\tau}+1]}\cap\mathcal{T}_{\{\tau,\bar{\tau}\},x}$.
If $p_{\tau}=0$ or $q_{\tau}=r$ this is 0.
\end{prop}

We can summarize the situation as follows. The Kodaira-Spencer isomorphism
induces a ``vertical'' decomposition of the tangent bundle into
the direct sum of the $\mathcal{T}_{\{\sigma,\bar{\sigma}\}}$. The
foliation $\mathscr{F}_{\{\tau,\bar{\tau}\}}$ is ``vertical'' in
the sense that it combines the subspace $\mathscr{E}_{\{\tau,\bar{\tau}\}}\subset\mathcal{T}_{\{\tau,\bar{\tau}\}}$
with the full $\mathcal{T}_{\{\sigma,\bar{\sigma}\}}$ for all $\{\sigma,\bar{\sigma}\}\ne\{\tau,\bar{\tau}\}.$
The cascade structure, on the other hand, results from the slope decomposition
of $X(d,\mathfrak{f}),$ and induces a ``horizontal'' filtration
by the $\mathcal{T}^{[a,b]}$ on the tangent bundle. The proposition
describes the way these vertical decomposition and horizontal filtration
interact.
\begin{proof}
We first compute the dimension of $\mathcal{T}_{x}^{[p_{\tau},q_{\tau}+1]}\cap\mathcal{T}_{\{\tau,\bar{\tau}\},x}$.
The graded pieces of $\mathcal{T}_{x}^{[p_{\tau},q_{\tau}+1]}$ are
the $\mathcal{T}_{x}^{(a,b)}$ for $1\le a\le p_{\tau}$ and $q_{\tau}+1\le b\le r.$
The dimension of $\mathcal{T}_{x}^{(a,b)}$ is the dimension of the
formal group 
\[
G^{a,b}={\rm Ext}(X^{(a)},X^{(b)})={\rm Ext}(X(1,\mathfrak{g}^{(a)}),X(1,\mathfrak{g}^{(b)}))^{d^{a}d^{b}},
\]
 which is $d^{a}d^{b}\sum_{i\in\mathscr{I}_{\mathfrak{P}}}(\mathfrak{g}^{(b)}(i)-\mathfrak{g}^{(a)}(i)).$
The dimension of $\mathcal{T}_{x}^{[p_{\tau},q_{\tau}+1]}$ is therefore
\[
\sum_{1\le a\le p_{\tau}}\sum_{q_{\tau}+1\le b\le r}d^{a}d^{b}\sum_{i\in\mathscr{I}_{\mathfrak{P}}}(\mathfrak{g}^{(b)}(i)-\mathfrak{g}^{(a)}(i)),
\]
and the dimension of $\mathcal{T}_{x}^{[p_{\tau},q_{\tau}+1]}\cap\mathcal{T}_{\{\tau,\bar{\tau}\},x}$
is the contribution of $i=\tau$ to this sum. But from the way $p_{\tau}$
and $q_{\tau}$ were defined it follows that for the particular index~$i$ corresponding to $\tau$ we have, for all $a$ and $b$ in the
above range, $\mathfrak{g}^{(b)}(i)=1$ and $\mathfrak{g}^{(a)}(i)=0.$
It follows that
\[
\dim\mathcal{T}_{x}^{[p_{\tau},q_{\tau}+1]}\cap\mathcal{T}_{\{\tau,\bar{\tau}\},x}=\sum_{1\le a\le p_{\tau}}\sum_{q_{\tau}+1\le b\le r}d^{a}d^{b}=r_{\phi^{-1}\circ\tau}\cdot(d-r_{\tau})=\mathrm{rk}\;\mathscr{E}_{\{\tau,\bar{\tau}\}}.
\]
To conclude the proof of the Proposition it is therefore enough to
show the inclusion 
\[
W_{x}:=\mathcal{T}_{x}^{[p_{\tau},q_{\tau}+1]}\cap\mathcal{T}_{\{\tau,\bar{\tau}\},x}\subset\mathscr{E}_{\{\tau,\bar{\tau}\}},
\]
i.e. that $W_{x}$ annihilates $\text{\rm KS}(\mathcal{P}_{\tau}[V]\otimes\mathcal{P}_{\bar{\tau}}).$
Let
\[
i:\mathfrak{S}\hookrightarrow {\rm Spf}(\widehat{\mathcal{O}}_{M,x})
\]
be the infinitesimal neighborhood of $x$ in the direction of $W_{x}.$
More precisely,
\[
\mathfrak{S}={\rm Spf}((\mathcal{O}_{M,x}/\mathfrak{m}_{x}^{2})/((\mathfrak{m}_{x}/\mathfrak{m}_{x}^{2})[W_{x}]))
\]
where $(\mathfrak{m}_{x}/\mathfrak{m}_{x}^{2})[W_{x}]$ is the subspace
of the cotangent space at $x$ annihilated by $W_{x}\subset\mathcal{T}_{x}.$
To conform with \cite{G-dS1} we introduce the following notation.
Let $\mathcal{A}=i^{*}A^\text{\rm univ}$, 
\[
\mathcal{H}=H_{dR}^{1}(\mathcal{A}/\mathfrak{S})=i^{*}H_{dR}^{1}(A^\text{\rm univ}/\widehat{\mathcal{O}}_{M,x}),\,\,\,\mathcal{P}=i^{*}\mathcal{P}_{\tau},\,\,\,\mathcal{P}_{0}=\mathcal{P}[V],\,\,\,\mathcal{Q}=i^{*}\mathcal{P}_{\bar{\tau}}.
\]
Over $M^\text{\rm ord}$, the $p$-divisible group $A^\text{\rm univ}[\mathfrak{P}^{\infty}]$
has an $\mathcal{O}_{K}$-stable slope filtration by $p$-divisible
groups
\[
\mathscr{X}^{(r)}\subset\mathscr{X}^{(r-1,r)}\subset\cdots\subset\mathscr{X}^{(1,r)}=A^\text{\rm univ}[\mathfrak{P}^{\infty}],
\]
characterized by the fact that at each geometric point $x\in M^\text{\rm ord}$
\[
\mathscr{X}_{x}^{(a,r)}=X^{(a)}\times\cdots\times X^{(r)},
\]
where $X^{(\nu)}$ is the $\nu$-th isoclinic factor of $A_{x}^\text{\rm univ}[\mathfrak{P}^{\infty}]$
(the slopes increasing with $\nu$). 

Let
\[
0\subset {\rm Fil}^{2}=i^{*}(\mathscr{X}^{(q_{\tau}+1,r)})\subset {\rm Fil}^{1}=i^{*}(\mathscr{X}^{(p_{\tau}+1,r)})\subset {\rm Fil}^{0}=\mathcal{A}[\mathfrak{P}^{\infty}].
\]
It follows from the construction of the cascade $\mathscr{C}=\{\Gamma^{(a,b)}\}$
that while the full $\mathcal{A}[\mathfrak{P}^{\infty}]$ does deform
over $\mathfrak{S}$, its subquotient $p$-divisible groups ${\rm Fil}^{1}$
and ${\rm Fil}^{0}/{\rm Fil}^{2}$ are constant there: they are obtained (with
their $\mathcal{O}_{K}$-structure) by base change from the fiber
at $x$. Indeed, as the only non-zero graded pieces of $\mathcal{T}_{x}^{[p_{\tau},q_{\tau}+1]}$
are the $\mathcal{T}_{x}^{(a,b)}$ for $a\le p_{\tau}$ and $q_{\tau}+1\le b$,
only extensions of $X^{(a)}$ by $X^{(b)}$ for $a$ an $b$ in these
ranges contribute to the deformation of $\mathcal{A}[\mathfrak{P}^{\infty}]$
over $\mathfrak{S}$. But such an $X^{(a)}$ does not participate
in ${\rm Fil}^{1},$ and $X^{(b)}$ does not participate in ${\rm Fil}^{0}/{\rm Fil}^{2}$.

It also follows from our choice of $p_{\tau}$ and $q_{\tau}$ that
$\mathcal{P}_{0}$ pairs trivially with the tangent space to ${\rm Fil}^{2}$,
so can be considered a subspace of the cotangent space of the $p$-divisible
group $G={\rm Fil}^{0}/{\rm Fil}^{2}.$ Let $D(G)=\mathbb{D}(G)_{\mathfrak{S}}$
be the evaluation of the Dieudonné crystal associated to $G$ on $\mathfrak{S}.$
This is an $\mathcal{O}_{\mathfrak{S}}$-module with $\mathcal{O}_{K}$-action,
and it is equipped with an integrable connection $\nabla$. We have
$D(G)\subset D(\mathcal{A}[p^{\infty}])=H_{dR}^{1}(\mathcal{A}/\mathfrak{S})=\mathcal{H}$
and 
\[
\nabla:D(G)\to D(G)\otimes\Omega_{\mathfrak{S}}^{1}
\]
is induced by the Gauss-Manin connection on $\mathcal{H}.$ Therefore,
the diagram
\[ \xymatrix@C=1cm@M=0.3cm{\mathcal{P} = \omega_{\mathcal{A}/\mathfrak{S}}[\tau] \ar@{^{(}->}[r]\ar[d]_{{\rm KS}_{\mathcal{A}/\mathfrak{S}}} & \mathcal{H}[\tau]\ar[d]^\nabla\\ 
\mathcal{Q}^\vee\otimes \Omega^1_{\mathfrak{S}} \underset{\lambda}{\stackrel{\sim}{\longrightarrow}} \omega^\vee_{\mathcal{A}^t/\mathfrak{S}}[\tau]\otimes \Omega^1_\mathfrak{S}& \mathcal{H}[\tau] \otimes  \Omega^1_\mathfrak{S} \ar@{->>}[l]
}\]
used to compute $\text{\rm KS}_{\mathcal{A}/\mathfrak{\mathfrak{S}}},$ can be
replaced, when we restrict to $\mathcal{P}_{0}=\mathcal{P}[V]$, by
the diagram
\[ \xymatrix@C=1cm@M=0.3cm{\mathcal{P}_0 \ar@{^{(}->}[r]\ar[d]_{{\rm KS}_{\mathcal{A}/\mathfrak{S}}} & D(G)[\tau]\ar[d]^\nabla\\ \mathcal{Q}^\vee\otimes \Omega^1_{\mathfrak{S}} & D(G)[\tau] \otimes  \Omega^1_\mathfrak{S}. \ar[l]
}\]
However, by the constancy of $G={\rm Fil}^{0}/{\rm Fil}^{2}$ over $\mathfrak{S}$,
the left arrow vanishes. We must therefore have $\text{\rm KS}_{\mathcal{A}/\mathfrak{S}}(\mathcal{P}_{0}\otimes\mathcal{Q})=0.$
By the functoriality of Kodaira-Spencer homomorphisms with respect
to base change we conclude that $W_{x},$ the tangent space to $\mathfrak{S},$
annihilates $\text{\rm KS}_{A^\text{\rm univ}/M}(\mathcal{P}_{\tau}[V]\otimes\mathcal{P}_{\bar{\tau}})\subset\Omega_{M/k}^{1},$
i.e. is contained in $\mathscr{E}_{\{\tau,\bar{\tau}\}},$ and this
completes the proof.
\end{proof}

\subsection{The extension of the foliations to inner Ekedahl-Oort strata}

\subsubsection{The problem}

This section is largely combinatorial. The variety $M$ has a stratification
by locally closed subsets $M_{w}$, of which $M^\text{\rm ord}$ is the largest,
called the Ekedahl-Oort (EO) strata of $M$. They are labelled by
certain Weyl group cosets $w$ (or by their distinguished representatives
of shortest length), recalled below. They are characterized by the
fact that the isomorphism class of $A_{x}^\text{\rm univ}[p]$, with its $\mathcal{O}_{K}$-structure
and polarization, is the same for all geometric points $x$ in a given
stratum, and the strata are maximal with respect to this property.
See \cite{Mo2,Wed2}.

Consider a (not necessarily closed) point $x\in M.$ Let $k(x)$ be
its residue field. By $(\ref{eq:ker_V_times_V})$ and the fact that
\[
\mathcal{P}_{\tau,x}[V]\otimes\mathcal{P}_{\bar{\tau},x}\cap\mathcal{P}_{\tau,x}\otimes\mathcal{P}_{\bar{\tau},x}[V]=\mathcal{P}_{\tau,x}[V]\otimes\mathcal{P}_{\bar{\tau},x}[V],
\]
the dimension of $\ker(V\otimes V:(\mathcal{P_{\tau}}\otimes\mathcal{P}_{\bar{\tau}})_{x}\to(\mathcal{P}_{\phi^{-1}\circ\tau}^{(p)}\otimes\mathcal{P}_{\phi^{-1}\circ\bar{\tau}}^{(p)})_{x}$)
is given by
\begin{equation}
r_{V}\{\tau,\bar{\tau}\}(x)=\dim\mathcal{P}_{\tau,x}[V]\cdot r_{\bar{\tau}}+r_{\tau}\cdot\dim\mathcal{P}_{\bar{\tau},x}[V]-\dim\mathcal{P}_{\tau,x}[V]\cdot\dim\mathcal{P}_{\bar{\tau},x}[V].\label{eq:r_v(x) recall}
\end{equation}

This quantity depends only on $A_{x}^\text{\rm univ}[p]$, and is therefore
\emph{constant along each EO strata}. At $x\in M^\text{\rm ord}$ it was expressed
in terms of the CM type by the formula
\[
r_{V}^\text{\rm ord}\{\tau,\bar{\tau}\}=\max\{0,r_{\tau}-r_{\phi^{-1}\circ\tau}\}\cdot(d-r_{\tau})+r_{\tau}\cdot\max\{0,r_{\phi^{-1}\circ\tau}-r_{\tau}\}.
\]
(We did it at geometric points, but the same works at every schematic
point of $M^\text{\rm ord},$ not necessarily closed.) Since $r_{V}\{\tau,\bar{\tau}\}(x)$
can only increase under specialization, $r_{V}^\text{\rm ord}\{\tau,\bar{\tau}\}$
is the minimal value attained by it, and
\begin{equation}
M_{\Sigma}=\{x\in M|\,r_{V}\{\tau,\bar{\tau}\}(x)=r_{V}^\text{\rm ord}\{\tau,\bar{\tau}\}\,\,\forall\{\tau,\bar{\tau}\}\in\Sigma\}\label{eq:M_=00005CSigma}
\end{equation}
is a Zariski open set, which is a union of EO strata. Clearly,
\[
M_{\Sigma_{1}\cup\Sigma_{2}}=M_{\Sigma_{1}}\cap M_{\Sigma_{2}}.
\]
In our earlier work \cite{G-dS1}, where there was only one $\Sigma$
to consider, this set was denoted by $M_{\sharp}.$ 

The proof of Lemma \ref{lem:Foliation is smooth} shows that over
$M_{\Sigma}$ the sheaf
\begin{equation}
\mathscr{F}_{\Sigma}=\text{\rm KS}\left(\sum_{\{\tau,\bar{\tau}\}\in\Sigma}(\mathcal{P}_{\tau}\otimes\mathcal{P}_{\bar{\tau}})[V\otimes V]\right)^{\perp}\label{eq:the foliation}
\end{equation}
remains a vector sub-bundle of $\mathcal{T}$ of corank $r_{V}(\Sigma)=\sum_{\{\tau,\bar{\tau}\}\in\Sigma}r_{V}^\text{\rm ord}\{\tau,\bar{\tau}\}.$
The properties of being involutive and $p$-closed extend by continuity
from the open dense $\mu$-ordinary stratum. We conclude:
\begin{thm}
The vector bundle $\mathscr{F}_{\Sigma}$ extends as a \emph{smooth}
$p$-foliation to the Zariski open set $M_{\Sigma}$.
\end{thm}

Our task is therefore to calculate $r_{V}\{\tau,\bar{\tau}\}(x)$
at $x\in M_{w}$ for the different EO strata $M_{w}$, and see for
which of them it equals $r_{V}^\text{\rm ord}\{\tau,\bar{\tau}\}.$ This will
give us an explicit description of $M_{\Sigma}$.

\subsubsection{Previous results}

When $L=\mathbb{Q}$ and $K$ is quadratic imaginary, and when $p\mathcal{O}_{K}=\mathfrak{P}$
is inert (the split case being trivial for quadratic imaginary $K$),
we showed in \cite{G-dS1} that there exists a smallest EO stratum
$M^\text{\rm fol}$ in $M_{\sharp}$. Any other stratum lies in $M_{\sharp}$
if and only if it contains $M^\text{\rm fol}$ in its closure. In this case
$\Sigma=\{\tau,\bar{\tau}\}$, and writing $(r_{\tau},r_{\bar{\tau}})=(n,m)$
we may assume, without loss of generality, that $n\ge m.$ We find
that $r_{V}(\Sigma)=(n-m)m,$ hence 
\[
\mathrm{rk}\;\mathscr{F}_{\Sigma}=m^{2}.
\]
The dimension of $M^\text{\rm fol}$ is also $m^{2},$ while the dimension
of $M$ itself is $nm$. For example, for Shimura varieties attached
to the group $U(n,1),$ the foliation extends everywhere except to
the lowest, $0$-dimensional, EO stratum, consisting of the so-called
superspecial points.

In case $K$ is quadratic imaginary, the labelling of the EO strata
of $M$ is by $(n,m)$\emph{-shuffles}. We review it to motivate the
type of combinatorics that will show up in the general case. A permutation
$w$ of $\{1,\dots,n+m\}$ is called an $(n,m)$-shuffle if
\[
w^{-1}(1)<\cdots<w^{-1}(n),\,\,\,w^{-1}(n+1)<\cdots<w^{-1}(n+m),
\]
i.e. $w$ interlaces the intervals $[1,n]$ and $[n+1,n+m]$ but keeps
the order within each interval. If $w$ is an $(n,m)$-shuffle and
$M_{w}$ is the corresponding EO stratum, then its dimension is equal
to the length of $w$
\[
\dim(M_{w})=\ell(w)=\sum_{i=1}^{n}(w^{-1}(i)-i).
\]
The unique $w$ for which this gets the value $nm$ is
\begin{equation}
w^\text{\rm ord}=\left(\begin{array}{cccccc}
1 & \cdots & m & m+1 & \cdots & n+m\\
n+1 & \cdots & n+m & 1 & \cdots & n
\end{array}\right)\label{eq:ordinary w}
\end{equation}
and the corresponding $M_{w}$ is $M^\text{\rm ord}.$ The formula for $r_{V}\{\tau,\bar{\tau}\}(x)$
at $x\in M_{w}(k)$ reads
\[
r_{V}\{\tau,\bar{\tau}\}(x)=a(w)\cdot m=|\{1\le i\le n|\,1\le w^{-1}(i)\le n\}|\cdot m
\]
(compare \cite{G-dS1} (2.3)). For example, for $w^\text{\rm ord}$ this is
$(n-m)m.$ It is readily seen that there is a unique $(n,m)$-shuffle
$w^\text{\rm fol}$ for which $a(w)$ is still $n-m,$ but for which $\dim(M_{w})=m^{2}$
is the smallest possible. This is the permutation
\[
w^\text{\rm fol}=\left(\begin{array}{ccccccccc}
1 & \cdots & n-m & n-m+1 & \cdots & n & n+1 & \cdots & n+m\\
1 & \cdots & n-m & n+1 & \cdots & n+m & n-m+1 & \cdots & n
\end{array}\right),
\]
whose corresponding EO stratum is $M^\text{\rm fol}.$

\subsubsection{Weyl group cosets and EO strata}

We return to a general CM field $K$. Denote by $\Pi_{e,d-e}$ the
set of $(e,d-e)$-shuffles in the symmetric group $\mathfrak{S}_{d}.$
They serve as representatives (of minimal length) for
\[
\mathfrak{S}_{e}\times\mathfrak{S}_{d-e}\setminus\mathfrak{S}_{d}.
\]
For an $(e,d-e)$-shuffle $\pi$ let
\[
\check{\pi}=w_{0}\circ\pi\circ w_{0}
\]
where $w_{0}(\nu)=d+1-\nu$ is the element of maximal length in $\mathfrak{S}_{d}.$
This $\check{\pi}$ is a $(d-e,e)$-shuffle, and $\pi\mapsto\check{\pi}$
is a bijection between $\Pi_{e,d-e}$ and $\Pi_{d-e,e}.$ Explicitly,
\[
\check{\pi}^{-1}(d+1-\nu)=d+1-\pi^{-1}(\nu).
\]

Let $w=(w_{\tau})_{\tau\in\mathscr{I}}$ where $w_{\tau}\in\Pi_{r_{\tau},d-r_{\tau}}$
and $w_{\bar{\tau}}=\check{w}_{\tau}.$ Note that $w_{\tau}$ and
$w_{\bar{\tau}}$, being conjugate by $w_{0}$, have the same length.

Let $k$ be, as usual, an algebraically closed field containing $\kappa.$
Consider the following Dieudonné module with $\mathcal{O}_{K}$-structure
$N_{w}$ attached to $w$:
\begin{itemize}
\item $N_{w}=\bigoplus_{i\in\mathscr{I}}\bigoplus_{j=1}^{d}ke_{i,j},$ $\mathcal{O}_{K}$
acting on $\bigoplus_{j=1}^{d}ke_{i,j}$ via $i:\mathcal{O}_{K}\to k.$
\vspace{0.2cm}
\item $F(e_{i,j})=\begin{cases}
\begin{array}{c}
0\\
e_{i+1,m}
\end{array} & \begin{array}{c}
\mathrm{if}\,\,w_{i}(j)\le\mathfrak{f}(i)\\
\mathrm{if}\,\,w_{i}(j)=\mathfrak{f}(i)+m.
\end{array}\end{cases}$
\vspace{0.2cm}

\item $V(e_{i+1,j})=\begin{cases}
\begin{array}{c}
0\\
e_{i,n}
\end{array} & \begin{array}{c}
\mathrm{if}\,\,j\le d-\mathfrak{f}(i)\\
\mathrm{if}\,\,j=d-\mathfrak{f}(i)+w_{i}(n).
\end{array}\end{cases}$
\end{itemize}
This $N_{w}$ is endowed with a pairing, setting its $\tau$ and $\bar{\tau}$-components
in duality, but we suppress it from the notation, as it is irrelevant
to the computation that we have to make.
\begin{prop}
(\cite{Mo2}, Theorem 6.7) The EO strata $M_{w}$ of $M$ are in one-to-one
correspondence with the $w$'s as above. If $x\in M_{w}(k)$ is a
geometric point, the Dieudonné module of $A_{x}^\text{\rm univ}[p]$ is isomorphic,
with its $\mathcal{O}_{K}$-structure, to $N_{w}.$ The dimension
of $M_{w}$ is given by
\begin{equation}
\dim(M_{w})=\ell(w):=\sum_{\{\tau,\bar{\tau}\}\in\mathscr{I}^{+}}\ell(w_{\tau}).\label{eq:EO dimension formula}
\end{equation}
\end{prop}

As an example, the $\mu$-ordinary stratum $M^\text{\rm ord}$ corresponds
to $w=w^\text{\rm ord}=(w_{i}^\text{\rm ord})_{i\in\mathscr{I}}$ where $w_{i}^\text{\rm ord}$
is given by $(\ref{eq:ordinary w})$ with $(n,m)=(\mathfrak{f}(i),d-\mathfrak{f}(i)).$

We now consider $r_{V}\{\tau,\bar{\tau}\}(x)$, as in (\ref{eq:r_v(x) recall}),
for $x\in M_{w}(k).$ As usual we write $\tau=i$ and $r_{\tau}=\mathfrak{f}(i)$.
We have
\begin{equation}
\begin{split}
\dim\mathcal{P}_{\tau,x}[V]& =|\{j|\,j\le r_{\phi^{-1}\circ\bar{\tau}},\,w_{\tau}(j)\le r_{\tau}\}|\label{eq:P_tau_V}
\\
& =|\{j|\,j\le d-\mathfrak{f}(i-1),\,w_{i}(j)\le\mathfrak{f}(i)\}|
\\
&\ge\max\{0,\mathfrak{f}(i)-\mathfrak{f}(i-1)\},
\end{split}
\end{equation}
and
\begin{equation}
\begin{split}
\dim\mathcal{P}_{\bar{\tau},x}[V]&=|\{j|\,j\le r_{\phi^{-1}\circ\tau},\,w_{\bar{\tau}}(j)\le r_{\bar{\tau}}\}|
\\& =|\{\ell|\,\ell\le r_{\bar{\tau}},\,w_{\bar{\tau}}^{-1}(\ell)\le r_{\phi^{-1}\circ\tau}\}|
\\
&=|\{\ell|\,\ell\le r_{\bar{\tau}},\,\check{w}_{\tau}^{-1}(\ell)\le r_{\phi^{-1}\circ\tau}\}|
\\
&=|\{m|\,d+1-m\le r_{\bar{\tau}},\,d+1-w_{\tau}^{-1}(m)\le r_{\phi^{-1}\circ\tau}\}|
\\
&=|\{j|\,r_{\phi^{-1}\circ\bar{\tau}}+1\le j,\,r_{\tau}+1\le w_{\tau}(j)\}|
\\
&=|\{j|\,d-\mathfrak{f}(i-1)+1\le j,\,\mathfrak{f}(i)+1\le w_{i}(j)\}|
\\
&\ge\max\{0,\mathfrak{f}(i-1)-\mathfrak{f}(i)\}.
\end{split}
\end{equation}

\begin{lem}
We have $r_{V}\{\tau,\bar{\tau}\}(x)=r_{V}^\text{\rm ord}\{\tau,\bar{\tau}\}$
precisely when the following conditions are satisfied:
\end{lem}

\begin{itemize}
\item If $\mathfrak{f}(i-1)\le\mathfrak{f}(i)$, 
\begin{equation}
|\{j|\,j\le d-\mathfrak{f}(i-1),\,w_{i}(j)\le\mathfrak{f}(i)\}|=\mathfrak{f}(i)-\mathfrak{f}(i-1),\label{eq:Cond_1-1}
\end{equation}
\item If $\mathfrak{f}(i)\le\mathfrak{f}(i-1)$,
\begin{equation}
|\{j|\,d-\mathfrak{f}(i-1)+1\le j,\,\mathfrak{f}(i)+1\le w_{i}(j)\}|=\mathfrak{f}(i-1)-\mathfrak{f}(i).\label{eq:Cond_2-1}
\end{equation}
\end{itemize}
\begin{proof}
We begin by noting that when $\mathfrak{f}(i-1)=\mathfrak{f}(i)$
the two conditions agree with each other, since each of them is equivalent,
in this case, to $w_{i}=w_{i}^\text{\rm ord}$, so the Lemma is consistent.

Let us dispose first of the cases where $\mathfrak{f}(i)\in\{0,d\}.$
In these extreme cases $\mathcal{P}_{\tau}$ or $\mathcal{P}_{\bar{\tau}}$
are zero, so
\[
r_{V}\{\tau,\bar{\tau}\}(x)=r_{V}^\text{\rm ord}\{\tau,\bar{\tau}\}=0
\]
for all $x.$ The conditions of the Lemma also hold then, trivially,
everywhere. We therefore assume from now on that $0<\mathfrak{f}(i)<d.$ 

Suppose $\mathfrak{f}(i-1)\le\mathfrak{f}(i).$ If $(\ref{eq:Cond_1-1})$
holds then 
\[
\dim\mathcal{P}_{\tau,x}[V]\otimes\mathcal{P}_{\bar{\tau},x}=(\mathfrak{f}(i)-\mathfrak{f}(i-1))\cdot(d-\mathfrak{f}(i))=r_{V}^\text{\rm ord}\{\tau,\bar{\tau}\}.
\]
We claim that $|\{j|\,d-\mathfrak{f}(i-1)+1\le j,\,\mathfrak{f}(i)+1\le w_{i}(j)\}|=0$,
and therefore $\mathcal{P}_{\bar{\tau},x}[V]=0$ and $r_{V}\{\tau,\bar{\tau}\}(x)=r_{V}^\text{\rm ord}\{\tau,\bar{\tau}\}$
as desired. Suppose, to the contrary, that for some $j$ we had $d-\mathfrak{f}(i-1)<j$
and $\mathfrak{f}(i)<w_{i}(j)$. Then there would be \emph{fewer}
than $\mathfrak{f}(i-1)$ values of $j$ such that $d-\mathfrak{f}(i-1)<j$
and $w_{i}(j)\le\mathfrak{f}(i),$ hence \emph{more than} $\mathfrak{f}(i)-\mathfrak{f}(i-1)$
values of $j$ such that $j\le d-\mathfrak{f}(i-1)$ and $w_{i}(j)\le\mathfrak{f}(i).$
This would violate condition $(\ref{eq:Cond_1-1})$.

Conversely, still under $\mathfrak{f}(i-1)\le\mathfrak{f}(i)$, assume
that $r_{V}\{\tau,\bar{\tau}\}(x)=r_{V}^\text{\rm ord}\{\tau,\bar{\tau}\}=(\mathfrak{f}(i)-\mathfrak{f}(i-1))\cdot(d-\mathfrak{f}(i)).$
Since we assumed that $d-\mathfrak{f}(i)\ne0$ we can not have $\dim\mathcal{P}_{\tau,x}[V]>\mathfrak{f}(i)-\mathfrak{f}(i-1)$
and condition $(\ref{eq:Cond_1-1})$ must hold.

The arguments when $\mathfrak{f}(i)\le\mathfrak{f}(i-1)$ are analogous.
\end{proof}
It is easy to see that if $\mathfrak{f}(i-1)\le\mathfrak{f}(i)$ the
$(\mathfrak{f}(i),d-\mathfrak{f}(i))$-shuffle $w_{i}^\text{\rm fol}$ given
by

{\scriptsize{}
\[
\left(\begin{array}{ccccccccc}
1 & \cdots & \mathfrak{f}(i)-\mathfrak{f}(i\!-\!1) & \mathfrak{f}(i)-\mathfrak{f}(i\!-\!1)+1 & \cdots & d-\mathfrak{f}(i\!-\!1) & d-\mathfrak{f}(i\!-\!1)+1 & \cdots & d\\
1 & \cdots & \mathfrak{f}(i)-\mathfrak{f}(i\!-\!1) & \mathfrak{f}(i)+1 & \cdots & d & \mathfrak{f}(i)-\mathfrak{f}(i\!-\!1)+1 & \cdots & \mathfrak{f}(i)
\end{array}\right)
\]
}satisfies $(\ref{eq:Cond_1-1})$, and this is the $(\mathfrak{f}(i),d-\mathfrak{f}(i))$-shuffle
of smallest length satisfying it. Its length is then
\[
\ell(w_{i}^\text{\rm fol})=(d-\mathfrak{f}(i))\mathfrak{f}(i-1).
\]

Similarly if $\mathfrak{f}(i)\le\mathfrak{f}(i-1)$ letting $\mathfrak{g}(i)=d-\mathfrak{f}(i)$
the same holds with $w_{i}^\text{\rm fol}$ given by

{\scriptsize{}
\[
\left(\begin{array}{ccccccccc}
1 & \cdots & \mathfrak{g}(i\!-\!1) & \mathfrak{g}(i\!-\!1)+1 & \cdots & \mathfrak{g}(i\!-\!1)+\mathfrak{f}(i) & \mathfrak{g}(i\!-\!1)+\mathfrak{f}(i)+1 & \cdots & d\\
\mathfrak{f}(i)+1 & \cdots & \mathfrak{g}(i\!-\!1)+\mathfrak{f}(i) & 1 & \cdots & \mathfrak{f}(i) & \mathfrak{g}(i\!-\!1)+\mathfrak{f}(i)+1 & \cdots & d
\end{array}\right).
\]
}In this case the length is
\[
\ell(w_{i}^\text{\rm fol})=\mathfrak{f}(i)(d-\mathfrak{f}(i-1)).
\]

We arrive at the following result.
\begin{thm}
The EO stratum $M_{w}\subset M_{\Sigma}$ if and only if for every
$\{\tau,\bar{\tau}\}\in\Sigma$, writing $\tau=i,$ $\phi^{-1}\circ\tau=i-1$
as usual, $(\ref{eq:Cond_1-1})$ or $(\ref{eq:Cond_2-1})$ hold.

There exists a unique EO stratum $M_{w}\subset M_{\Sigma}$ of smallest
dimension. It is given by the following recipe: $w_{i}=w_{i}^\text{\rm fol}$
if $\{\tau,\bar{\tau}\}\in\Sigma$ and $w_{i}=id.$ otherwise. Denote
this $M_{w}$ by $M_{\Sigma}^\text{\rm fol}.$ Its dimension is given by
\[
\dim M_{\Sigma}^\text{\rm fol}=\sum_{\{\tau,\bar{\tau}\}\in\Sigma}\min(r_{\tau},r_{\phi^{-1}\circ\tau})\cdot\min(r_{\bar{\tau}},r_{\phi^{-1}\circ\bar{\tau}}).
\]

Any other EO stratum $M_{w}$ lies in $M_{\Sigma}$ if and only if
$M_{\Sigma}^\text{\rm fol}$ lies in its closure.
\end{thm}

\begin{proof}
By the computations above, the unique $M_{w}$ of smallest dimension
which is still contained in $M_{\Sigma}$, i.e. for which $r_{V}\{\tau,\bar{\tau}\}(x)=r_{V}^\text{\rm ord}\{\tau,\bar{\tau}\}$
for all $\{\tau,\bar{\tau}\}\in\Sigma$, is obtained when $w_{i}=w_{i}^\text{\rm fol}$
whenever $\{\tau,\bar{\tau}\}\in\Sigma$ (this condition is symmetric
in $\tau$ and $\bar{\tau}$), and $w_{i}=id.$ otherwise. Its dimension
follows from $(\ref{eq:EO dimension formula})$, and the computation
of the lengths of the $w_{i}^\text{\rm fol}.$

If $M_{\Sigma}^\text{\rm fol}$ lies in the closure of $M_{w}$ then for $x\in M_{\Sigma}^\text{\rm fol}$
and $y\in M_{w}$ and $\{\tau,\bar{\tau}\}\in\Sigma$ we have
\[
r_{V}^\text{\rm ord}\{\tau,\bar{\tau}\}\le r_{V}\{\tau,\bar{\tau}\}(y)\le r_{V}\{\tau,\bar{\tau}\}(x)=r_{V}^\text{\rm ord}\{\tau,\bar{\tau}\}
\]
so equality holds and $M_{w}\subset M_{\Sigma}$. Conversely, suppose
that $M_{w}\subset M_{\Sigma}$. Assume, without loss of generality,
that $r_{\phi^{-1}\circ\tau}\le r_{\tau}.$ Then condition $(\ref{eq:Cond_1-1})$
holds, so writing $\tau=i$, we must have that $w_{i}$ is the permutation

{\scriptsize{}
\[
\left(\begin{array}{ccccccccc}
1 & \cdots & \mathfrak{f}(i)-\mathfrak{f}(i\!-\!1) & \mathfrak{f}(i)-\mathfrak{f}(i\!-\!1)+1 & \cdots & d-\mathfrak{f}(i\!-\!1) & d-\mathfrak{f}(i\!-\!1)+1 & \cdots & d\\
* & \cdots & * & * & \cdots & * & \mathfrak{f}(i)-\mathfrak{f}(i\!-\!1)+1 & \cdots & \mathfrak{f}(i)
\end{array}\right).
\]
}Since the blocks $[1,\mathfrak{f}(i)-\mathfrak{f}(i-1)]$ and $[\mathfrak{f}(i)+1,d]$
must appear in the bottom row in increasing order (but interlaced),
it is easy to check that the permutation $w_{i}^\text{\rm fol}$ is smaller
than or equal to $w_{i}$ in the Bruhat order on the Weyl group of
$GL_{d}.$ This is enough (although in general, not equivalent) for
$M_{\Sigma}^\text{\rm fol}$ to lie in the closure of $M_{w}.$ (For the closure
relation between EO strata, see \cite{V-W}.)
\end{proof}

\subsubsection{Integral varieties} The results of \S \ref{sec:general integral} imply that integral varieties of the foliation $\mathscr{F}_\Sigma$ abound. Nonetheless, it's interesting to identify specific examples. 
By Lemma \ref{lem:Foliation is smooth}, when $\Sigma=\mathscr{I}^{+},$
$\dim M_{\Sigma}^\text{\rm fol}=\mathrm{rk}\;\mathscr{F}_{\Sigma}.$ We expect
$M_{\Sigma}^\text{\rm fol}$ to be an integral variety of $\mathscr{F}_{\Sigma}$
in this case. This has been proved when $K$ is quadratic imaginary
in \cite{G-dS1}, Theorem 25, and would be analogous to Theorem \ref{thm:HBMV integral varieties}
above. The proof of Theorem 25 in \cite{G-dS1} was not conceptual,
and involved tedious computations with Dieudonné modules. 

\subsection{The moduli problem $M^{\Sigma}$ and the extension of the foliation
to it}

\subsubsection{The moduli scheme $M^{\Sigma}$}

As we have seen in the previous section, $\mathscr{F}_{\Sigma}$,
defined by $(\ref{eq:the foliation}),$ is a $p$-foliation on all
of $M$, but is smooth only on the Zariski open $M_{\Sigma}.$ Our
goal is to define a ``successive blow-up'' $\beta:M^{\Sigma}\to M,$
which is an isomorphism over $M_{\Sigma},$ and a natural extension
of $\mathscr{F}_{\Sigma}$ to a \emph{smooth} $p$-foliation on $M^{\Sigma}$.
See \cite{G-dS1} §4.1 for $K$ quadratic imaginary, where $M^{\Sigma}$
was denoted $M^{\sharp}.$

Fix $\{\tau,\bar{\tau}\}\in\Sigma$, and assume that $r_{\phi^{-1}\circ\tau}\le r_{\tau}$
(otherwise switch notation between $\tau$ and $\bar{\tau}$). Consider
the moduli problem $M^{\tau,\bar{\tau}}$ on $\kappa$-algebras $R$
given by
\[
M^{\tau,\bar{\tau}}(R)=\{(\underline{A},\mathcal{N})|\,\underline{A}\in M(R),\,\mathcal{N}\subset\mathcal{P}_{\tau}[V]\,\mathrm{a\,subbundle\,of\,rank}\,r_{\tau}-r_{\phi^{-1}\circ\tau}\}/\simeq.
\]
(Here, as before, by a \emph{subbundle}, we mean a subsheaf that is locally a direct summand.)
The forgetful map $\beta:M^{\tau,\bar{\tau}}\to M$ is bijective above
$M_{\Sigma}$, since 
\[
\mathrm{rk}(\mathcal{P}_{\tau}[V])=r_{\tau}-r_{\phi^{-1}\circ\tau}
\]
holds there.

Let $(n,m)=(r_{\tau},r_{\phi^{-1}\circ\tau})$ and consider the relative
Grassmannian $Gr(n-m,\mathcal{P}_{\tau})$ over $M$, classifying sub-bundles
$\mathcal{N}$ of rank $n-m$ in the rank $n$ bundle $\mathcal{P}_{\tau}.$
This is a smooth scheme over $M,$ of relative dimension $(n-m)m.$
As the condition $V(\mathcal{N})=0$ is closed, the moduli problem
$M^{\tau,\bar{\tau}}$ is representable by a closed subscheme of $Gr(n-m,\mathcal{P}_{\tau}).$
The fiber $M_{x}^{\tau,\bar{\tau}}=\beta^{-1}(x)$ is the Grassmannian
of $(n-m)$-dimensional subspaces in $\mathcal{P}_{\tau,x}[V].$ 

Suppose $x\in M_{w}$ where $w=(w_{\sigma})_{\sigma\in\mathscr{I}},$
$w_{\sigma}\in\Pi_{r_{\sigma},d-r_{\sigma}}$ and $w_{\bar{\sigma}}=\check{w}_{\sigma}$.
As we have computed in $(\ref{eq:P_tau_V})$
\[
\dim\mathcal{P}_{\tau,x}[V]=a_{\tau}(w):=|\{j| j\le d-r_{\phi^{-1}\circ\tau}, w_{\tau}(j)\le r_{\tau}\}|\ge n-m,
\]
and consequently
\[
\dim M_{x}^{\tau,\bar{\tau}}=(n-m)(a_{\tau}(w)-n+m).
\]

Denote by $M_{w}^{\tau,\bar{\tau}}=\beta^{-1}(M_{w}).$ We have shown
that it is smooth of relative dimension $(n-m)(a_{\tau}(w)-n+m)$
over $M_{w}.$ When the EO strata undergo specialization, this dimension
jumps.

The main result concerning $M^{\tau,\bar{\tau}}$ is the following.
\begin{thm}
The scheme $M^{\tau,\bar{\tau}}$ is a non-singular variety, and $\beta$
induces a bijection between its irreducible ( = connected) components
and those of $M$. In particular, $M^{\tau,\bar{\tau},{\rm ord}}$ is dense
in $M^{\tau,\bar{\tau}}.$
\end{thm}

\begin{proof}
This is, mutatis mutandis, the proof of Theorem 15 in \cite{G-dS1}
§4.1.3, and we refer to our earlier paper for details.
\end{proof}
\begin{cor}
The moduli problem which is the fiber product of the $M^{\tau,\bar{\tau}}$
over~$M$, for all $\{\tau,\bar{\tau}\}\in\Sigma$, is represented
by a smooth $\kappa$-variety $M^{\Sigma}$. The map $\beta:M^{\Sigma}\to M$
induces a bijection on irreducible components, and is an isomorphism
over $M_{\Sigma}.$ In particular, $M^{\Sigma,{\rm ord}}$ is dense in $M^{\Sigma}.$
Over an EO stratum $M_{w}$ that is not contained in $M_{\Sigma}$
the map $\beta$ is no longer an isomorphism, but it is smooth of
relative dimension
\[
\sum_{\{\tau,\bar{\tau}\}\in\Sigma}(r_{\tau}-r_{\phi^{-1}\circ\tau})(a_{\tau}(w)-r_{\tau}+r_{\phi^{-1}\circ\tau}).
\]
Here for each pair $\{\tau,\bar{\tau}\}$ we choose $\tau$ so that
$r_{\tau}\ge r_{\phi^{-1}\circ\tau}$, and
\[
a_{\tau}(w)=|\{j| j\le d-r_{\phi^{-1}\circ\tau}, w_{\tau}(j)\le r_{\tau}\}|.
\]
\end{cor}

We can also draw the following corollary.
\begin{cor}
The open set $M_{\Sigma}\subset M$ is the largest open set in $M$
to which the foliation $\mathscr{F}_{\Sigma}$ extends as a smooth
foliation.
\end{cor}

\begin{proof}
We have already noted that $\mathscr{F}_{\Sigma}$ is a saturated
foliation everywhere on $M.$ Outside $M_{\Sigma}$ the dimension
of its fibers is strictly larger than their dimension over $M^\text{\rm ord},$
so $\mathscr{F}_{\Sigma}$ can not be a vector sub-bundle on any open
set larger than $M_{\Sigma}.$
\end{proof}
\begin{rem}
Since $\beta:M^{\Sigma}\to M$ is a birational projective morphism
between non-singular varieties, it follows from the general theory
(\cite{Stacks} 29.43 and \cite{H} Theorem 7.17 and exercise 7.11(c))
that $\beta$ is a blow up at an ideal sheaf supported on $M-M_{\Sigma}.$
\end{rem}

\subsubsection{Extending the foliation to $M^{\Sigma}$}

Above $M_{\Sigma}$, the map $\beta$ is an isomorphism
\[
M^{\Sigma}\supset\beta^{-1}(M_{\Sigma})\simeq M_{\Sigma},
\]
so the foliation $\mathscr{F}_{\Sigma}$ induces a foliation on $\beta^{-1}(M_{\Sigma}).$
We now explain how to extend it to a smooth foliation on all of $M^{\Sigma}.$

Let $k$ be, as usual, an algebraically closed field containing $\kappa$,
and $y\in M^{\Sigma}(k)$ a geometric point, with image $x=\beta(y)\in M(k).$
The point $y$ ``is'' a pair $(\underline{A}_{x}^\text{\rm univ},\mathcal{N}_{y})$
where $\mathcal{N}_{y}=(\mathcal{N}_{y,\tau})$. Here
\begin{itemize}
\item $\tau$ range over representatives of the pairs $\{\tau,\bar{\tau}\}\in\Sigma$,
chosen in such a way that $r_{\phi^{-1}\circ\tau}\le r_{\tau},$
\item $\mathcal{N}_{y,\tau}\subset\mathcal{P}_{x,\tau}[V]$ is an $r_{\tau}-r_{\phi^{-1}\circ\tau}$-dimensional
subspace.
\end{itemize}
Let $\mathcal{H}_{x}=H_{dR}^{1}(A_{x}^\text{\rm univ}/k)$ and $\underline{\omega}_{x}=H^{0}(A_{x}^\text{\rm univ},\Omega^{1}).$
The polarization $\lambda_{x}$ on $A_{x}^\text{\rm univ}$ induces a perfect
pairing
\[
\{,\}_{\lambda}:\underline{\omega}_{x}\times\mathcal{H}_{x}/\underline{\omega}_{x}\to k
\]
satisfying $\{\iota(a)u,v\}_{\lambda}=\{u,\iota(\bar{a})v\}_{\lambda}$
for $a\in\mathcal{O}_{K}$. Adapting the proof of Theorem 15 in \cite{G-dS1}
§4.1.3 to our situation, we see that the tangent space $\mathcal{T}M_{y}^{\Sigma}$
at the point~$y$ can be described as the set of pairs $(\varphi,\psi)$
where
\begin{itemize}
\item $\varphi\in\mathrm{Hom}_{\mathcal{O}_{K}}(\underline{\omega}_{x},\mathcal{H}_{x}/\underline{\omega}_{x})^{sym}$
is an $\mathcal{O}_{K}$-linear homomorphism from $\underline{\omega}_{x}$
to $\mathcal{H}_{x}/\underline{\omega}_{x}$ which is symmetric with
respect to $\{,\}_{\lambda},$ i.e. satisfies $\{u,\varphi(v)\}_{\lambda}=\{v,\varphi(u)\}_{\lambda}$
for $u,v\in\underline{\omega}_{x}$. By Kodaira-Spencer,
such a $\varphi$ represents a tangent vector to $M$ at $x,$ and
the projection $(\varphi,\psi)\mapsto\varphi$ corresponds to the
map $d\beta:\mathcal{T}M_{y}^{\Sigma}\to\mathcal{T}M_{x}.$ We can
write $\varphi$ as a tuple $(\varphi_{\tau})$ where $\tau\in\mathscr{I}$
and $\varphi_{\tau}\in\mathrm{Hom}(\mathcal{P}_{x,\tau},\mathcal{H}_{x,\tau}/\mathcal{P}_{x,\tau}$).
The symmetry condition and the fact that $\{,\}_{\lambda}$ induces
a perfect pairing between $\mathcal{P}_{x,\bar{\tau}}$ and $\mathcal{H}_{x,\tau}/\mathcal{P}_{x,\tau}$
for every $\tau\in\mathscr{I},$ imply that $\varphi_{\bar{\tau}}$
is determined by $\varphi_{\tau}$, and that for a given choice of
representatives~$\tau$ for the pairs $\{\tau,\bar{\tau}\}\in\mathscr{I}^{+}$
the $\varphi_{\tau}$ may be chosen arbitrarily.
\item $\psi\in\mathrm{Hom}_{\mathcal{O}_{K}}(\mathcal{N}_{y},\mathcal{H}_{x}[V]/\mathcal{N}_{y})$
satisfies $\varphi|_{\mathcal{N}_{y}}=\psi\mod\underline{\omega}_{x}.$
\end{itemize}
The second condition means that $\psi=(\psi_{\tau})$ where the $\tau$
range over the same set of representatives for $\{\tau,\bar{\tau}\}\in\Sigma$
as above, $\psi_{\tau}\in\mathrm{Hom}(\mathcal{N}_{y,\tau},\mathcal{H}_{x,\tau}[V]/\mathcal{N}_{y,\tau})$
and $\varphi_{\tau}|_{\mathcal{N}_{y,\tau}}=\psi_{\tau}\mod\mathcal{P}_{x,\tau}$.

The tangent space to the fiber $\beta^{-1}(x)$ is the space of $(\varphi,$$\psi)$
with $\varphi=0,$ or, alternatively, the space
\[
\{(\psi_{\tau})|\,\psi_{\tau}\in\mathrm{Hom}(\mathcal{N}_{y,\tau},\mathcal{P}_{x,\tau}[V]/\mathcal{N}_{y,\tau})\}.
\]

\begin{thm}
There exists a unique smooth $p$-foliation $\mathscr{F}^{\Sigma}$
on $M^{\Sigma}$ , characterized by the fact that at any geometric
point $y\in M^{\Sigma}(k)$ as above, $\mathscr{F}_{y}^{\Sigma}$
is the subspace
\[
\mathscr{F}_{y}^{\Sigma}=\{(\varphi,\psi)|\,\psi=0\}.
\]
The foliation $\mathscr{F}^{\Sigma}$ agrees with $\mathscr{F}_{\Sigma}$
on $\beta^{-1}(M_{\Sigma})\simeq M_{\Sigma}$ and, in general, is
transversal to the fibers of $\beta$.
\end{thm}

\begin{proof}
A straightforward adaptation of the proof of Proposition 20 in \cite{G-dS1}
§4.3.
\end{proof}
Intuitively, $\mathscr{F}_{y}^{\Sigma}$ are the directions in $\mathcal{T}M_{y}^{\Sigma}$
in which $\mathcal{N}_{y}$ does not undergo any infinitesimal deformation.
By the alluded transversality, its projection to $\mathcal{T}M_{x}$
is injective, and identifies $\mathscr{F}_{y}^{\Sigma}$ with the
subspace of $\varphi=(\varphi_{\tau})_{\tau\in\mathscr{I}}$ such
that for every $\{\tau,\bar{\tau}\}\in\Sigma$ (and $\tau$ a representative
as above) $\varphi_{\tau}(\mathcal{N}_{y,\tau})=0.$ This subspace,
however, varies with $y\in\beta^{-1}(x).$

\end{document}